\documentclass[a4paper]{amsart}

\RequirePackage{amsmath} 
\RequirePackage{amssymb}
\usepackage{amscd,latexsym,amsthm,amsfonts,amssymb,amsmath,amsxtra}
\usepackage[colorlinks=true,urlcolor=blue,citecolor=blue]{hyperref}
\usepackage{color}
\usepackage[all]{xy}
\usepackage[OT2,T1]{fontenc}
\usepackage{bm}
\usepackage{mathtools}
\usepackage{ mathrsfs }
\usepackage{xcolor}
\usepackage{comment}
\usepackage{marginnote}
\usepackage{enumitem}
\usepackage{thmtools, thm-restate}

\DeclareSymbolFont{cyrletters}{OT2}{wncyr}{m}{n}
\DeclareMathSymbol{\Sha}{\mathalpha}{cyrletters}{"58}

\let\Re\undefined

\DeclareMathOperator{\Re}{Re}

\DeclareMathOperator{\supp}{supp}

\DeclareMathOperator{\Spec}{Spec}
\DeclareMathOperator{\Hom}{Hom}

\newcommand{\Res}{\operatorname{Res}}

\newcommand{\Geo}{\operatorname{Geo}}

\newcommand{\K}{\operatorname{K}}

\newcommand{\Ad}{\operatorname{Ad}}

\newcommand{\diag}{\operatorname{diag}}
\newcommand{\Ind}{\operatorname{Ind}}

\newcommand{\fin}{\operatorname{fin}}

\newcommand*{\transp}[2][-3mu]{\ensuremath{\mskip1mu\prescript{\smash{\mathrm t\mkern#1}}{}{\mathstrut#2}}}

\newcommand{\Reg}{\operatorname{Reg}}

\newcommand{\Vol}{\operatorname{Vol}}

\newcommand{\sm}{\operatorname{Small}}
\newcommand{\du}{\operatorname{Dual}}
\newcommand{\bi}{\operatorname{Big}}

\newcommand{\Kuz}{\operatorname{Kuz}}

\newcommand{\ER}{\operatorname{ER}}

\newcommand{\gen}{\operatorname{gen}}

\newcommand{\Main}{\operatorname{Main}}

\newcommand{\dist}{\operatorname{dist}}

\newcommand{\Lie}{\operatorname{Lie}}
\newcommand{\RNum}[1]{\uppercase\expandafter{\romannumeral #1\relax}}

\begin{document}
	
	\theoremstyle{plain}
	\newtheorem{thm}{Theorem}[section]
	
	\newtheorem{cor}[thm]{Corollary}
	\newtheorem{thmy}{Theorem}
	\renewcommand{\thethmy}{\Alph{thmy}}
	\newenvironment{thmx}{\stepcounter{thm}\begin{thmy}}{\end{thmy}}
	\newtheorem{cory}{Corollary}
	\renewcommand{\thecory}{\Alph{cory}}
	\newenvironment{corx}{\stepcounter{thm}\begin{cory}}{\end{cory}}
	\newtheorem*{thma}{Theorem A}
	\newtheorem*{corb}{Corollary B}
	\newtheorem*{thmc}{Theorem C}
	\newtheorem{lemma}[thm]{Lemma}  
	\newtheorem{prop}[thm]{Proposition}
	\newtheorem{conj}[thm]{Conjecture}
	\newtheorem{fact}[thm]{Fact}
	\newtheorem{claim}[thm]{Claim}
	
	\theoremstyle{definition}
	\newtheorem{defn}[thm]{Definition}
	\newtheorem{example}[thm]{Example}
	\theoremstyle{remark}
	
	\newtheorem{remark}[thm]{Remark}	
	\numberwithin{equation}{section}

\title[Rankin-Selberg $L$-functions for $\mathrm{GL}(n+1)\times \mathrm{GL}(n)$]{Relative Trace Formula, Subconvexity and Quantitative Nonvanishing of Rankin-Selberg $L$-functions for $\mathrm{GL}(n+1)\times\mathrm{GL}(n)$}%
\author{Liyang Yang}
	
\address{Department of Mathematics, Texas A\&M University, College Station,  TX 77843, USA} 
\email{liyangy@tamu.edu} 
	
\begin{abstract}
We establish hybrid subconvexity bounds and quantitative simultaneous
nonvanishing for Rankin--Selberg $L$-functions on
$\mathrm{GL}(n+1)\times\mathrm{GL}(n)$ over $\mathbb Q$. For a fixed unitary
cuspidal representation $\pi'$ of $\mathrm{GL}(n)$ and any unitary pure
isobaric representation $\pi$ of $\mathrm{GL}(n+1)$, we prove the
$t$-aspect bound
\begin{align*}
L(1/2+it,\pi\times\pi')\ll_{\pi,\pi',\varepsilon}
(1+|t|)^{\frac{n(n+1)}{4}-\frac{1}{4(3n^2+n-1)}+\varepsilon}.
\end{align*}
When $\pi$ is tempered, the saving improves to $1/(4(3n-1))$; in particular,
the resulting bound for standard $L$-functions improves the general-rank
exponent of Nelson. We also obtain,
in both the spectral and level aspects, the first quantitative simultaneous
nonvanishing result in higher rank for central Rankin--Selberg values in
suitable cuspidal families. Our proofs use an amplified, regularized relative
trace formula whose spectral expansion retains the full generic spectrum and
is expressed directly in terms of central values. On the geometric side, we
estimate collective combinations of regular orbital integrals rather than
the individual orbits separately. This structure leads to an optimized count
of rational parameters, sharper amplification bounds, and a treatment of the
level aspect. In rank one, the resulting exponents agree with the Burgess
bounds.
\end{abstract}
	
\date{\today}%
\maketitle
\enlargethispage{2pt}
\tableofcontents

\section{Introduction}

In this paper, we use an amplified, regularized relative trace formula for
$\mathrm{GL}(n+1)\times\mathrm{GL}(n)$ to establish hybrid subconvexity bounds
and quantitative simultaneous nonvanishing for Rankin--Selberg $L$-functions.
Its spectral side is expressed directly in terms of central Rankin--Selberg
values across the full generic spectrum, while its geometric side is analyzed
through collective combinations of orbital integrals.

Fix a unitary cuspidal representation $\pi'$ of $\mathrm{GL}(n)/\mathbb{Q}$, and let
$\pi=\pi_1\boxplus\cdots\boxplus\pi_m$ be a unitary pure isobaric automorphic
representation of $\mathrm{GL}(n+1)/\mathbb{Q}$ (see \cite{Lan79}). Thus each
$\pi_j$ is a unitary cuspidal representation of
$\mathrm{GL}(n_j)/\mathbb{Q}$ and $n_1+\cdots+n_m=n+1$. The two extreme cases
are $m=1$, when $\pi$ is cuspidal, and $m=n+1$, when every $\pi_j$ is a Hecke
character and $\pi$ belongs to the minimal Eisenstein spectrum. The completed
Rankin--Selberg $L$-function factors as
\begin{equation}\label{1.1}
\Lambda(s,\pi\times\pi')=\Lambda(s,\pi_1\times\pi')\cdots \Lambda(s,\pi_m\times\pi'),
\end{equation}
where the factors on the right are the Rankin--Selberg $L$-functions of
Jacquet--Piatetski-Shapiro--Shalika \cite{JPSS83}. If
$\pi=\textbf{1}\boxplus\cdots\boxplus\textbf{1}$, then \eqref{1.1} reduces to
$\Lambda(s,\pi\times\pi')=\Lambda(s,\pi')^{n+1}$, with $\Lambda(s,\pi')$
the standard $L$-function of Godement--Jacquet \cite{GJ72}. Langlands
functoriality predicts more generally that Rankin--Selberg $L$-functions factor
into shifted products of standard $L$-functions \cite{Lan70}.

Throughout the paper, $L(s,\pi\times\pi')$ denotes the finite part of the
completed $L$-function. It is given by an Euler product over the rational primes
in $\Re(s)>1$ and by meromorphic continuation elsewhere.

\subsection{\texorpdfstring{Subconvexity of $L(s,\pi\times\pi')$ in the $t$-Aspect}
{Subconvexity of L-Functions in the t-Aspect}}\label{subc}
The subconvexity problem asks for bounds on the critical line that improve on
those obtained from the functional equation and the Phragm\'en--Lindel\"of
principle. We refer to \cite{Fri95,IS00,Mic07,Mun18a} for surveys and historical
background. 

In this paper we consider bounds for Rankin--Selberg $L$-functions of the form 
\begin{equation}\label{1.4}
 	L(1/2+it,\pi\times\pi')\ll_{\pi,\pi',\varepsilon} (1+|t|)^{\frac{(1-\delta)n(n+1)}{4}+\varepsilon},
 	\end{equation}
where $\delta\in[0,1]$ and the implied constant may depend on $\pi$, $\pi'$,
and $\varepsilon$. Convexity gives \eqref{1.4} with $\delta=0$, whereas the
generalized Lindel\"of hypothesis predicts $\delta=1$. Thus subconvexity in the
$t$-aspect amounts to establishing \eqref{1.4} for some $\delta>0$ depending at
most on $n$.

The cases $n=1$ and $n=2$ have been treated in, among other works,
\cite{Goo81,Goo82,Li09,Li11,MV10,Mun18b,BB20,LS21}. In arbitrary rank, Nelson
\cite{Nel21} proved $t$-aspect subconvexity for standard $L$-functions on
$\mathrm{GL}(n)/\mathbb{Q}$. This yields \eqref{1.4} unconditionally when
$\pi=\textbf{1}\boxplus\cdots\boxplus\textbf{1}$, with a saving
$\delta\asymp n^{-5}$. Subject to Langlands functoriality, the factorization
into standard $L$-functions would give the general case as well.

Our first main result proves \eqref{1.4} unconditionally for every unitary pure
isobaric $\pi$, with $\delta\asymp n^{-4}$; in the tempered case, this improves
to $\delta\asymp n^{-3}$.
\begin{thmx}\label{A}
Let $n\geq 1$. Let $\pi$ be a unitary pure isobaric automorphic representation
of $\mathrm{GL}(n+1)/\mathbb{Q}$, and let $\pi'$ be a unitary cuspidal
representation of $\mathrm{GL}(n)/\mathbb{Q}$. Then
\begin{equation}\label{1}
L(1/2+it,\pi\times\pi')\ll_{\pi,\pi',\varepsilon}(1+|t|)^{\frac{n(n+1)}{4}-\frac{1}{4(3n^2+n-1)}+\varepsilon},
\end{equation}
where the implied constant depends on $\varepsilon$ and the conductors of $\pi$
and $\pi'$. If $\pi$ is tempered, then
\begin{equation}\label{1.4*}
L(1/2+it,\pi\times\pi')\ll_{\pi,\pi',\varepsilon}(1+|t|)^{\frac{n(n+1)}{4}-\frac{1}{4(3n-1)}+\varepsilon}.
	\end{equation}
\end{thmx}

Applying the tempered estimate in Theorem \ref{A} to
$\pi=\textbf{1}\boxplus\cdots\boxplus\textbf{1}$ gives the following
subconvex bound for standard $L$-functions.
\begin{cor}[Standard $L$-functions]\label{cor1.2}
Let $n\geq 1$, let $\delta=\frac{1}{n(n+1)(3n-1)}$, and let $\pi'$ be a
unitary cuspidal representation of $\mathrm{GL}(n)/\mathbb{Q}$. Then
\begin{equation}\label{1.}
L(1/2+it, \pi')\ll_{\pi',\varepsilon}(1+|t|)^{\frac{(1-\delta)n}{4}+\varepsilon},
\end{equation} 
where the implied constant depends on $\varepsilon$ and the conductor of $\pi'$.
\end{cor}

\begin{remark}
	When $n=1$, the estimate \eqref{1.4*} and Corollary \ref{cor1.2} give the
	Burgess exponents $3/8$ and $3/16$, respectively, in agreement with
	\cite{Yan23b}. This rank-one compatibility reflects the sharp dependence on
	the Hecke normalization retained in the regular orbital contribution. For
	$n\geq2$, the assertions follow from Theorem \ref{B} below.
\end{remark}

\begin{remark}
	For $n\geq2$, 
	Corollary \ref{cor1.2} improves the exponent in \cite[Theorem 1.1]{Nel21}
	from $\delta=\frac{2}{n^2(3n^3-2n^2-1)}\asymp n^{-5}$ to
	$\delta=\frac{1}{n(n+1)(3n-1)}\asymp n^{-3}$.
\end{remark}

%\marginpar{the bound $O(n^{-4})$ is optimal in the current approach (amplification in \cite{Iwa92}); see \textsection?}

\begin{cor}[The case of $\mathrm{GL}(n)\times\mathrm{GL}(m)$]
Let $n\geq 1$, let $\pi'$ be a unitary cuspidal representation of
$\mathrm{GL}(n)/\mathbb{Q}$, and let $\sigma$ be a unitary cuspidal
representation of $\mathrm{GL}(m)/\mathbb{Q}$ with $m\mid(n+1)$. Then
	\begin{align*}
		L(1/2+it, \sigma\times \pi')\ll_{\sigma,\pi',\varepsilon}(1+|t|)^{\frac{(1-\delta)nm}{4}+\varepsilon}
	\end{align*}
	for $\delta=\frac{1}{n(n+1)(3n^2+n-1)},$ 
	where the implied constant depends on $\pi'$, $\sigma$, and $\varepsilon$.
\end{cor}
%\begin{remark}
%	The approach of Nelson \cite{Nel21} should also yield
%	$$
%	L(1/2+it, \sigma\times \pi')\ll_{\pi',\varepsilon}(1+|t|)^{\frac{(1-\delta)nm}{4}+\varepsilon}
%	$$
%	with $\delta=\frac{1}{3n^5-2n^4-2n^2},$ for a unitary cuspidal representation $\sigma$ of $\mathrm{GL}(m)/\mathbb{Q}$ with $n\equiv 1\pmod{m}$.
%\end{remark}

\begin{comment}
	\begin{remark}
		When $n=2$ and $F=\mathbb{Q},$ \eqref{1} improves \cite{Mun18b}, where Munshi proves $L(1/2+it,\pi\times f)\ll_{\pi,f,\varepsilon} (1+|t|)^{3/2-1/51+\varepsilon}$ for a Hecke-Maass form $\pi$ for $\mathrm{SL}(3,\mathbb{Z})$ and a holomorphic (or Maass) Hecke-cusp form $f$ for $\mathrm{SL}(2,\mathbb{Z})$.
	\end{remark}
\end{comment}

\subsection{Uniform Parameter Growth}\label{subc.}
Write the archimedean $L$-factor of $\pi$ as
\begin{align*}
\Lambda_{\infty}(s,\pi)=\prod_{1\leq j\leq n+1}\Gamma_{\mathbb{R}}(s+\lambda_{\pi_\infty,j}).
\end{align*}
Here $\lambda_{\pi_\infty,j}\in\mathbb{C}$ are the Langlands parameters. Given
$T\geq1$ and $C_{\infty}>c_{\infty}>0$, we say that $\pi$ has
\textit{uniform parameter growth of size $(T;c_{\infty},C_{\infty})$} if
\begin{equation}\label{1.6}
	c_{\infty}T\leq |\lambda_{\pi_\infty,j}|\leq C_{\infty}T,\ \ 1\leq j\leq n+1.
\end{equation}

For $n\geq2$, applying the following hybrid estimate to
$\pi\otimes|\det|^{it}$ yields Theorem \ref{A}. Its rank-one specialization
agrees with the relative trace formula analysis in \cite{Yan23b}.
\begin{thmx}\label{B} 
Let $n\geq 2$, and let $\pi=\pi_1\boxplus\cdots\boxplus\pi_m$ be a unitary
pure isobaric automorphic representation of $\mathrm{GL}(n+1)/\mathbb{Q}$.
Suppose that $\pi$ has arithmetic conductor $M$ and uniform parameter growth of
size $(T;c_{\infty},C_{\infty})$. Let $\pi'$ be a unitary cuspidal
representation of $\mathrm{GL}(n)/\mathbb{Q}$, with conductor $M'$. Set
\begin{align*}
M^\dagger:=\prod_{p\mid M,\ p\nmid M'}p^{e_p(M)},
\qquad R:=\frac{M}{M^\dagger},
\end{align*}
so that $R$ is the part of $M$ supported on the primes dividing $M'$. Define
\begin{align*}
\mathbf{L}:=
\begin{cases}
(M^\dagger)^{\frac{n}{2(3n^2+n-1)}}T^{\frac{1}{4(3n^2+n-1)}},\ \ & \text{if $(M^\dagger)^{n^2-n-1}\leq T^{n+1}$,}\\
(M^\dagger)^{\frac{1}{2(n^2+n+1)}}T^{\frac{1}{4(n^2+n+1)}},\ \ & \text{if $(M^\dagger)^{n^2-n-1}>T^{n+1}$.}
\end{cases}
\end{align*} 
Then
\begin{equation}\label{bur}
L(1/2,\pi\times\pi')\ll T^{\frac{n(n+1)}{4}+\varepsilon}M^{\frac{n}{2}+\varepsilon}R^{\frac n2}\mathbf{L}^{-1}\prod_{j=1}^m\sqrt{L(1,\pi_j,\Ad)},
\end{equation}
where the implied constant depends on $\varepsilon$, $c_{\infty}$,
$C_{\infty}$, and the conductor of $\pi'$. In particular,
\begin{equation}\label{b.}
L(1/2,\pi\times\pi')\ll T^{\frac{n(n+1)}{4}-\frac{1}{4(3n^2+n-1)}+\varepsilon},
\end{equation}
where the implied constant depends on $\varepsilon$, $c_{\infty}$,
$C_{\infty}$, and the conductors of $\pi_{\fin}$ and $\pi'$.
\end{thmx}

\begin{remark}
Nelson \cite{Nel23} proves a subconvex bound for
$L(1/2,\sigma_E\times\sigma_E')$, where $\sigma$ and $\sigma'$ are tempered
cuspidal representations of $\mathrm{U}(n+1)$ and an anisotropic
$\mathrm{U}(n)$, respectively, and their parameters grow uniformly; here the
subscript $E$ denotes quadratic base change. Estimate \eqref{b.} is a general
linear analogue with $\pi'$ fixed, and it does not require $\pi$ to be
cuspidal.
\end{remark}

\begin{remark}
The rank-one hybrid bound has a different form; see \cite{Yan23b}.
\end{remark}

%\begin{remark}
%The uniform parameter growth assumption is crucial in the proof of Theorem \ref{B}, as  in \cite{Nel23} and \cite{Nel21}. Using an amplified relative trace formula with a different choice of the test function, \cite{MRY22} gives a subconvex bound for Rankin-Selberg $L$-functions for $\mathrm{U}(3)\times \mathrm{U}(2)$ in certain  `conductor dropping' case, i.e., failing the uniform parameter growth condition.  
%\end{remark}

\begin{comment}
Since $L(1,\pi_j,\Ad)\ll (MT)^{o(1)},$ $1\leq j\leq m,$  then Theorem \ref{B} yields the following:
\begin{cor}[Hybrid Subconvexity]
Let notation be as before. Let $\alpha\in (0,1)$. Suppose that $T\gg_{\varepsilon} M^{\frac{4n^3+4n^2-n+1}{1-\alpha}+\varepsilon}$. Then   
\begin{align*}
	L(1/2,\pi\times\pi')\ll C(\pi\times\pi')^{\frac{1}{4}-\frac{\alpha}{4n(n+1)(2n+1)^2}},
\end{align*} 	
where $C(\pi\times\pi')$ is the analytic conductor, and the implied constant depends on $\varepsilon,$ $c_0,$ $C_0,$ $\pi',$ and $\pi_p$ at $p\mid M'$.  
\end{cor}
\end{comment}

\subsection{\texorpdfstring{Quantitative Nonvanishing of Rankin--Selberg
$L$-Functions}{Quantitative Nonvanishing of Rankin--Selberg L-Functions}}\label{nonv}

Nonvanishing of Rankin--Selberg $L$-functions is connected with several major
themes, including Landau--Siegel zeros \cite{IS00a}, Langlands functoriality
\cite{GJR04}, the Gan--Gross--Prasad conjecture \cite{Zha14a,Zha14b},
Whittaker periods \cite{GH16}, and the Bloch--Kato conjecture \cite{LTXZZ22}.

%\cite{BLM19}, 

Quantitative nonvanishing is known in several low-rank settings; see, for
example, \cite{Duk95,IS00a,MRY22,ST22,WYZ24}. In higher rank, simultaneous
nonvanishing has been established in certain families
\cite{Li09,Tsu21,JN21,Yan22,MRY24}, but the known results are qualitative. Our second
main result gives, for $n\geq2$, the first quantitative simultaneous
nonvanishing theorem for central values of Rankin--Selberg $L$-functions on
$\mathrm{GL}(n+1)\times\mathrm{GL}(n)$.

\begin{thmx}\label{thmC}
Let $n\geq2$. For $j=1,2$, let $\pi_j'=\otimes_{p\leq\infty}\pi_{j,p}'$ be a unitary
cuspidal representation of $\mathrm{GL}(n,\mathbb{A})$ with arithmetic
conductor $M_j'$. Suppose that $M_1'M_2'>1$, $(M_1',M_2')=1$, and
$\pi_{1,\infty}'\simeq\pi_{2,\infty}'$. Let $p_*\nmid M_1'M_2'$ be a prime and
$T,M\geq1$ with $(M,M_1'M_2'p_*)=1$. Fix a unitary Hecke character
$\omega$ of $\mathbb Q^\times\backslash\mathbb A^\times$, trivial at infinity,
a unitary irreducible admissible
representation $\pi_{\infty}$ of $\mathrm{PGL}_{n+1}(\mathbb{R})$ with uniform
parameter growth of size $(T;c_{\infty},C_{\infty})$, and a supercuspidal
representation $\pi_{p_*}$ of $\mathrm{GL}(n+1,\mathbb{Q}_{p_*})$ with
central character $\omega_{p_*}$. Let
$\mathcal{A}_0(T,M;\pi_{\infty},\pi_{p_*},\pi_1',\pi_2')$ be the cuspidal
family of Definition \ref{defn11.1}, localized to a fixed-width spectral
neighborhood of $\pi_\infty$. Then, for every
$\varepsilon>0$,
\begin{align*} 
\sum_{\substack{\pi\in \mathcal{A}_0(T,M;\pi_{\infty},\pi_{p_*},\pi_1',\pi_2')\\ L(1/2,\pi\times\pi_1')L(1/2,\pi\times\pi_2')\neq 0}}1\gg (TM)^{-\varepsilon}\begin{cases}
M^{\frac{n}{3n^2+n-1}}T^{\frac{1}{2(3n^2+n-1)}}, & \text{if $M\leq T^{\frac{n+1}{n^2-n-1}}$}\\
M^{\frac{1}{n^2+n+1}}T^{\frac{1}{2(n^2+n+1)}}, & \text{if $M>T^{\frac{n+1}{n^2-n-1}}$}
\end{cases}, 
\end{align*}
where the implied constant depends on $\varepsilon$, $c_{\infty}$,
$C_{\infty}$, $\omega$, $\pi_{p_*}$, $\pi_1'$, and $\pi_2'$. In particular, as
$T+M\to\infty$,
\begin{align*}
	\sum_{\substack{\pi\in \mathcal{A}_0(T,M;\pi_{\infty},\pi_{p_*},\pi_1',\pi_2')\\ L(1/2,\pi\times\pi_1')L(1/2,\pi\times\pi_2')\neq 0}}1\gg |\mathcal{A}_0(T,M;\pi_{\infty},\pi_{p_*},\pi_1',\pi_2')|^{\frac{1}{n(n+1)(3n^2+n-1)}-\varepsilon},
\end{align*}
where the implied constant depends on $\varepsilon$, $c_{\infty}$,
$C_{\infty}$, $\omega$, $\pi_{p_*}$, $\pi_1'$, and $\pi_2'$.
\end{thmx}
\begin{remark}
When $n=1$, the set $\mathcal{A}_0(T,M;\pi_{\infty},\pi_{p_*},\pi_1',\pi_2')$ may be viewed as a family of Maass cusp forms on $\mathrm{GL}(2)$ with spectral parameter in a fixed-width window centered at $T$ and level $\asymp M$.	
\end{remark}
%\begin{remark}
%	Theorem \ref{thmC} 
%\end{remark}
%When $T$ is fixed, Theorem \ref{thmC} is a quantitative form of  Theorem 9.9 in \cite{Yan22} when the base field is $\mathbb{Q}$.

\begin{remark}
For every $\varepsilon>0$, the Weyl law gives
\begin{align*}
T^{\frac{n(n+1)}{2}-\varepsilon}M^n
\ll
|\mathcal{A}_0(T,M;\pi_{\infty},\pi_{p_*},\pi_1',\pi_2')|
\ll
T^{\frac{n(n+1)}{2}+\varepsilon}M^n,
\end{align*}
with constants depending on $\varepsilon$, $c_{\infty}$, $C_{\infty}$,
$\pi_{p_*}$, $\pi_1'$, and $\pi_2'$.
\end{remark}

\subsection{Discussion of the Proofs}

The proofs combine the relative trace formula developed in \cite{Yan22} with
Nelson's construction of archimedean test functions
\cite[Parts 2 and 3]{Nel21}, his volume estimates
\cite[\S\S 15--16]{Nel23}, and arithmetic amplification in
the spirit of Duke--Friedlander--Iwaniec \cite{DFI02}. The arguments for
Theorems \ref{A} and \ref{B} extend the rank-one methods of
\cite{Yan23b,Yan23c}, while the proof of Theorem \ref{thmC} is motivated in
part by the nonvanishing argument of \cite{MRY22} for
$\mathrm{U}(2,1)\times\mathrm{U}(1,1)$.

\subsubsection{The formal model}

To motivate the strategy, consider the formal pre-trace identity
\begin{equation}\label{39}
J_{\Spec}(f,\varphi')
=
\int_{[G']}\int_{[G']}
\K(\iota(x),\iota(y))
\varphi'(x)\overline{\varphi'(y)}
\,dx\,dy
=
J_{\Geo}(f,\varphi'),
\end{equation}
where
\begin{align*}
[G']:=\mathrm{GL}(n,\mathbb{Q})\backslash\mathrm{GL}(n,\mathbb{A}),
\qquad
\iota(x):=\diag(x,1),
\end{align*}
and $\K$ is the kernel associated with
$f\in C_c^\infty(\mathrm{GL}_{n+1}(\mathbb{A}))$. The spectral side formally
decomposes as
\begin{align*}
J_{\Spec}(f,\varphi')
=
J_0(f,\varphi')+J_{\ER}(f,\varphi'),
\end{align*}
where $J_0$ and $J_{\ER}$ denote the contributions of the cuspidal and
noncuspidal spectra, respectively. Rankin--Selberg convolution identifies the
cuspidal contribution with
\begin{align*}
J_0(f,\varphi')
=
\sum_{\sigma}
|L(1/2,\sigma\times\pi')|^2
\beta_{f,\varphi'}(\sigma),
\end{align*}
where $\sigma$ ranges over cuspidal representations of
$\mathrm{GL}(n+1)/\mathbb{Q}$, $\varphi'$ belongs to $\pi'$, and
$\beta_{f,\varphi'}(\sigma)$ is a spectral weight.

For a fixed representation $\pi$ with uniform parameter growth of size
$(T;c_\infty,C_\infty)$, one would like to choose $f$ and $\varphi'$ so that
the spectral weights are nonnegative,
$\beta_{f,\varphi'}(\pi)\gg T^{-\alpha}$, and
$J_{\Geo}(f,\varphi')\ll T^\beta$. If the noncuspidal contribution is also
nonnegative, these estimates give
\begin{align*}
|L(1/2,\pi\times\pi')|^2
\leq
\frac{J_0(f,\varphi')}{\beta_{f,\varphi'}(\pi)}
\ll
T^\alpha J_{\Spec}(f,\varphi')
=
T^\alpha J_{\Geo}(f,\varphi')
\ll
T^{\alpha+\beta}.
\end{align*}
The convexity bound corresponds to
$\alpha+\beta=n(n+1)/2+\varepsilon$, while amplification yields the power
saving required for subconvexity. 

The difficulty is that \eqref{39} generally diverges when $\varphi'$ is a
genuine automorphic form. Nelson \cite[\textsection 5.3]{Nel21} avoids this
divergence by replacing $\varphi'$ with a wave packet and uses bilinear
estimates for the amplified geometric side to establish $t$-aspect
subconvexity for standard $L$-functions. This approach is effective for upper
bounds, but its spectral side averages integrals of $L$-functions along
vertical lines rather than values at the central point and is therefore not
applicable to the nonvanishing problem considered here.

\subsubsection{The relative trace formula}

We instead use the regularized relative trace formula of \cite{Yan22}; see
Theorem \ref{C}:
\begin{align*}
J_{\Spec}^{\Reg,\heartsuit}(f,\mathbf{0})
=
J_{\Geo}^{\Reg,\heartsuit}(f,\mathbf{0}).
\end{align*}
Its spectral side has the expansion
\begin{equation}\label{1.11}
J_{\Spec}^{\Reg,\heartsuit}(f,\mathbf{0})
=
\int_{\widehat{G(\mathbb{A})}_{\gen}}
|L(1/2,\sigma\times\pi')|^2
\beta_{f,\varphi'}(\sigma)
\,d\mu_\sigma.
\end{equation}
Two features of \eqref{1.11} are essential. First, the spectral expansion
encompasses the full generic spectrum, allowing us to treat pure isobaric
representations in Theorems \ref{A} and \ref{B}. Second, it is expressed
directly in terms of central Rankin--Selberg values, making it suitable for the
quantitative nonvanishing problem addressed in Theorem \ref{thmC}.

The main analytic task is therefore to estimate the geometric side. We retain
the Hecke structure of the small-cell and dual terms throughout the amplified
sum and exploit the collective combinations of regular orbital integrals in
the resulting rational-point count. Write
$M^\dagger:=\prod_{p\mid M,\,p\nmid M'}p^{e_p(M)}$ and
$R:=M/M^\dagger$. For the test function of
\textsection\ref{testfunction}, this gives
\begin{multline}\label{1.13}
\frac{\mathcal J_{\Geo}^{\heartsuit}
(\boldsymbol\alpha,\boldsymbol\ell)}
{\langle\varphi',\varphi'\rangle T^{n/2}}
\ll (TMM')^\varepsilon
\left(M^nR^nL
+\frac{M^{n-1}R^{n+1}L^{n(n+1)+2}}{T^{1/2}}\right)\\
+(TMM')^\varepsilon
\frac{R^{2n}L^{(3n-2)(n+1)+2}}{T^{1/2}}.
\end{multline}
The first term comes from the small-cell and dual contributions. At the
amplification primes, their local Hecke transforms are Satake polynomials
evaluated at the parameters of $\widetilde\pi'\boxplus(\omega\omega')$.
Estimating their amplified sum by Rankin--Selberg theory avoids any pointwise
loss from the Ramanujan bounds for $\pi'$. The remaining two terms come from
the regular orbital contribution.
For an individual component $f$, their dependence on the Hecke length is
$M^{n-1}R^{n+1}\mathcal N_f^n+R^{2n}\mathcal N_f^{3n-2}$, which is sharp for the
rational-point count permitted by the local support conditions. These
exponents grow linearly with $n$, whereas the corresponding bilinear estimates
in \cite[\S\S 5.9--5.11]{Nel21} have quadratic dependence
on $n$. This leads to the saving $\delta\asymp n^{-4}$ in the general pure
isobaric case.

For standard $L$-functions, the representation
$\mathbf{1}\boxplus\cdots\boxplus\mathbf{1}$ in the continuous spectrum is
tempered. Incorporating Rankin--Selberg estimates into the amplification
allows us to take $\mathcal{N}_f\asymp L$, rather than
$\mathcal{N}_f\ll L^{n+1}$ as in \cite{Nel21}. This gives the stronger saving
$\delta\asymp n^{-3}$ in Corollary \ref{cor1.2}.

Finally, the geometric estimates above yield power savings in both the spectral and level
aspects. For a suitable choice of local data, the noncuspidal contribution
vanishes, while the geometric side admits a matching lower bound. Combining
this lower bound with Theorem \ref{B} gives the quantitative nonvanishing
estimate in Theorem \ref{thmC}.

\subsection{Outline of the Paper}

\subsubsection{The relative trace formula and amplification}

In \textsection\ref{sec2}, we recall the relative trace formula developed in
\cite{Yan22}; see Theorem \ref{C}. In \textsection\ref{11.2}, we specify the
local and global data, construct the test functions, and introduce the
amplification parameters; see Theorem \ref{thmD'}. In
\textsection\ref{sec4.}, following \cite[Part 2]{Nel21}, we construct a suitable
vector $\phi'\in\pi'$ through its Kirillov model and establish the local
estimates needed below.

\subsubsection{The spectral side}

In \textsection\ref{sec5}, we analyze the amplified spectral side. Combining
this analysis with the local estimates from \textsection\ref{sec4.}, we obtain
a lower bound in terms of central $L$-values; see Proposition \ref{thm6}.

\subsubsection{The geometric side}

In \textsection\ref{sec5.1}, we describe the support of the local test function
at the amplification places in Iwasawa coordinates. This description is used
throughout the subsequent analysis of the geometric side.

In \textsection\ref{8.5.1}--\textsection\ref{sec10}, we decompose and estimate
$J_{\Geo}^{\Reg,\heartsuit}(f,\mathbf{0})$ as follows:
\begin{enumerate}
\item The small-cell orbital integral, one of the two main terms, is evaluated
in Proposition \ref{prop54}.
\item The dual orbital integral, obtained from the small-cell term by Poisson
summation, forms the second main term and is estimated in Proposition
\ref{8.1}.
\item The contribution from the regular orbital integrals, which constitutes
the main difficulty on the geometric side, is estimated in Theorem
\ref{thm9.1}. The proof exploits the geometric combinations of orbital integrals
arising from the relative trace formula and reduces their collective analysis
to an optimized rational-point counting problem. The inverse-transpose
argument used in this analysis appears in Lemmas \ref{lem9.14} and
\ref{lem9.15}.
\end{enumerate}

\subsubsection{Proofs of the main results}

In \textsection\ref{proof}, we choose the amplification data and combine the
spectral and geometric estimates. The resulting bounds, formulated in
Theorems \ref{E} and \ref{J}, imply Theorems \ref{A} and \ref{B}.

Finally, in \textsection\ref{sec12}, we choose local and global data adapted to
the nonvanishing problem. Proposition \ref{proposition12.2} gives the required
lower bound for $J_{\Geo}^{\Reg,\heartsuit}(f,\mathbf{0})$. Combining it with
Theorem \ref{E}, we prove Theorem \ref{G}, and hence Theorem \ref{thmC}.

\subsection{Further Directions}
The adelic formulation should extend to totally real fields with only local
modifications. Complex places require additional archimedean calculations. In
the level aspect, fixing the archimedean test function should similarly give a
quantitative simultaneous nonvanishing theorem over number fields, extending
\cite[Theorem C]{Yan22}.

The geometric analysis underlying \eqref{1.13} is sufficiently flexible to
accommodate other choices of local test functions and suggests applications to
subconvexity in finite-place aspects. Related results include Marshall's
depth-aspect bound for $\mathrm{U}(n+1)\times\mathrm{U}(n)$ \cite{Mar23},
the horizontal-aspect bounds of Hu--Nelson in the unitary setting
\cite{HN23}, and the level-aspect bound of Hu--Nelson for standard
$L$-functions on $\mathrm{GL}_n$ \cite{HN25}. It would be natural to
investigate analogous depth-aspect bounds, as well as sharper hybrid bounds in
the levels of $\pi$ and $\pi'$, using the present relative trace formula. 

%The case $n=1$ is treated in \cite{Yan23b,Yan23c}.

\subsection{Notation}\label{notation}
Let $\mathbb{A}$ be the adele group of $\mathbb{Q}$. For a rational prime $p,$ we denote by $\mathbb{Q}_p$ the corresponding local field and $\mathbb{Z}_p$ its ring of integers. Denote by $e_p(\cdot)$ the evaluation relative to $p$ normalized as $e_p(p)=1$. Denote by $\widehat{\mathbb{Z}}=\prod_{p<\infty}\mathbb{Z}_{p}$. Let $|\cdot|_p$ be the norm in $\mathbb{Q}_p$. Let $|\cdot|_{\fin}=\prod_{v<\infty}|\cdot|_p$. Let $|\cdot|_{\infty}$ be the norm in $\mathbb{R}$. Let $|\cdot|_{\mathbb{A}}=|\cdot|_{\infty}\otimes|\cdot|_{\fin}$. We will simply write $|\cdot|$ for $|\cdot|_{\mathbb{A}}$ in calculation over $\mathbb{A}^{\times}$ or its quotient by $\mathbb{Q}^{\times}$.   

Let $\psi=\otimes_{p\leq \infty}\psi_p$ be the additive character on $\mathbb{Q}\backslash \mathbb{A}$ such that $\psi(t)=\exp(2\pi it),$ for $t\in \mathbb{R}\hookrightarrow\mathbb{A}$. For $p\in \Sigma_{\mathbb{Q}},$ let $dt$ be the additive Haar measure on $\mathbb{Q}_p,$ self-dual relative to $\psi_p$. Then $dt=\prod_{p\leq\infty}dt$ is the standard Tamagawa measure on $\mathbb{A}$. Let $d^{\times}t=\zeta_{p}(1)dt/|t|_p,$ where $\zeta_{p}(\cdot)$ is the local Riemann zeta factor. In particular, $\Vol(\mathbb{Z}_p^{\times},d^{\times}t)=\Vol(\mathbb{Z}_p,dt)=1$ for every finite place $p$. Moreover, $\Vol(\mathbb{Q}\backslash\mathbb{A}; dt)=1$ and $\Vol(\mathbb{Q}\backslash\mathbb{A}^{(1)},d^{\times}t)=\underset{s=1}{\Res}\ \zeta(s)=1,$ where $\mathbb{A}^{(1)}$ is the subgroup of ideles $\mathbb{A}^{\times}$ with norm $1,$ and $\zeta(s)$ is the Riemann zeta function. More properties of $dt$ and $d^{\times}t$ can be found in \cite{Lan94}, Ch. XIV.

%Let $A_{G'}$ be  $Z'(\mathbb{R})^0$ embedded diagonally in the connected abelian Lie group $\prod_{v\mid\infty} Z'(\mathbb{Q}_p)^0$.	

Let $G=\mathrm{GL}(n+1)$ and $G'=\mathrm{GL}(n)$. Denote by $Z$ (resp. $Z'$) the center of $G$ (resp. $G'$). Let $G'^0(\mathbb{A})$ be the subgroup of $G'(\mathbb{A})$ consisting of $g'\in G'(\mathbb{A})$ with $|\det g'|=1$. Then the subgroup $Z'(\mathbb{A})G'^0(\mathbb{A})$ is open and has index $1$ (or $n-1$ if the base field is a function field) in $G'(\mathbb{A})$. So we may write $G'(\mathbb{A})=Z'(\mathbb{A})G'^0(\mathbb{A})$. Let $\overline{G}=Z\backslash G$ and $\overline{G'}=Z'\backslash G'$. We will identify $\overline{G'}$ with $G'^0$ as the conventional notation from Rankin-Selberg theory. For $x\in G(\mathbb{A})$ or $G'(\mathbb{A}),$ we denote by $x^{\vee}$ the transpose inverse of $x$. Fix the embedding from $G'$ to $G:$
\begin{align*}
	\iota:\ G'\longrightarrow G,\quad \gamma\mapsto \begin{pmatrix}
		\gamma&\\
		&1
	\end{pmatrix}.
\end{align*} 

For a matrix $g=(g_{i,j})\in G(\mathbb{A}),$ we denote by $E_{i,j}(g)=g_{i,j},$ the $(i,j)$-th entry of $g$. For a vector $\mathbf{v}=(v_1, \cdots, v_m)$ denote by $E_i(\mathbf{v})=v_i,$ the $i$-th component of $\mathbf{v}$. Let $m_1, m_2\in \mathbb{N}$. We write $M_{m_1, m_2}$ for the group of $m_1\times m_2$ matrices.

For an algebraic group $H$ over $\mathbb{Q}$, we will denote by $[H]:=H(\mathbb{Q})\backslash H(\mathbb{A})$. We equip measures on $H(\mathbb{A})$ as follows: for each unipotent group $U$ of $H,$ we equip $U(\mathbb{A})$ with the Haar measure such that, $U(\mathbb{Q})$ being equipped with the counting measure and the measure of $[U]$ is $1$. We equip the maximal compact subgroup $K$ of $H(\mathbb{A})$ with the Haar measure such that $K$ has total mass $1$. When $H$ is split, we also equip the maximal split torus of $H$ with Tamagawa measure induced from that of $\mathbb{A}^{\times}$. 

Let $\omega$ and $\omega'$ be unitary idele class characters on $\mathbb{A}^{\times}$. Denote by $\mathcal{A}_0\left([G],\omega\right)$ (resp. $\mathcal{A}_0\left([G'],\omega'\right)$) the set of cuspidal representations on $G(\mathbb{A})$ (resp. $G'(\mathbb{A})$) with central character $\omega$ (resp. $\omega'$). 

%Let $\tilde{f}\in C_c^{\infty}(G(\mathbb{A}))$. Denote by 
%\begin{equation}\label{test}
%\tilde{f}^{\omega}(g):=\int_{Z(\mathbb{A})}\tilde{f}(zg)\omega(z)d^{\times}z.
%\end{equation}

%\marginpar{\textcolor{red}{be careful of the notation $T$}}
Let $B$ (resp. $B'$) be the group of upper triangular matrices in $G$ (resp. $G'$). Let $T_B$ (resp. $T_{B'}$) be the diagonal  subgroup of $B$ (resp. $B'$). Let $A=Z\backslash T_B$ and $A'=Z'\backslash T_{B'}$. Let $N$ (resp. $N'$) be the unipotent radical of $B$ (resp. $B'$). Let $W_G$ be Weyl group of $G$ with respect to $(B,T_{B})$. Let $\Delta=\{\alpha_{1,2},\alpha_{2,3},\cdots,\alpha_{n,n+1}\}$ be the set of simple roots, and for each simple root $\alpha_{k,k+1},$ $1\leq k\leq n,$ denote by $w_k$ the corresponding reflection. Explicitly, for each $1\leq k\leq n,$
\begin{align*}
	w_k=\begin{pmatrix}
		I_{k-1} &\\
		& S&\\
		&&I_{n-k}
	\end{pmatrix},\ \text{where $S=\begin{pmatrix}
			&1\\
			1&
		\end{pmatrix}$}.
\end{align*}
For $1\leq k\leq n-1,$ denote by $w_k'$ the unique element in $G'$ such that $\iota(w_k')=w_k$. Then $w_k$'s generate the Weyl group $W_{G'}$ of $G'$ with respect to $(B',T_{B'})$. Let $\widetilde{w}_j'=w_j'w_{j+1}'\cdots w_{n-1}',$ $1\leq j\leq n-1$. Denote by $N_j'=N'/(\widetilde{w}_j 'N' \widetilde{w}_j'^{-1})$. Denote by $\widetilde{w}_n'=I_n$ and $N_n'=I_n$. Define the generic character $\theta$ on $[N]$ by setting $\theta(u)=\prod_{j=1}^{n}\psi(u_{j,j+1})$ for $u=(u_{i,j})_{1\leq i, j\leq n+1}\in N(\mathbb{A})$. Let $\theta'=\theta\mid_{[\iota(N')]}$ be the generic character on $[N']$.

Let $P$ (resp. $P'$) be the standard parabolic subgroup of $G$ (resp. $G'$) of type $(n,1)$ (resp. $(n-1,1)$). Let $P_0=Z\backslash P$ (resp. $P_0'=Z'\backslash P'$). We will denote by $Q$ a general parabolic subgroup of $G$. Denote by $N_P$ (resp. $N_Q$) the unipotent radical of $P$ (resp. $Q$).

For a function $h$ on $G(\mathbb{A}),$ we define $h^*$ by assigning $h^*(g)=\overline{h({g}^{-1})},$ $g\in G(\mathbb{A})$. Let $F_1(s), F_2(s)$ be two meromorphic functions. Denote by $A\asymp B$ for $A, B\in\mathbb{R}$ if there are absolute constants $c_1$ and $c_2$ such that $c_1B\leq A\leq c_2B$. Denote by $A\ll B$ for $A\in\mathbb{C}$ and $B\in\mathbb{R}_{>0}$ if there is an absolute constant $c$ such that $|A|\leq cB$. 

Throughout, we follow the $\varepsilon$-convention: that is, $\varepsilon$ will always be positive number which can be taken as small as we like, but may differ from one occurrence to another.

\textbf{Acknowledgements}
I am deeply grateful to Paul Nelson for his helpful discussions. I would also like to express my gratitude to Valentin Blomer, Yongxiao Lin, Simon Marshall, Philippe Michel, Dinakar Ramakrishnan, and Peter Sarnak for their precise comments and valuable suggestions.

\section{The Relative Trace Formula}\label{sec2}

We recall the form of the relative trace formula from \cite{Yan22} used in
this paper. For $i=1,2$, fix a unitary cuspidal representation $\pi_i'$ of
$G'(\mathbb{A})$ and a vector $\phi_i'\in\pi_i'$, with Whittaker function
$W_{\phi_i'}'$ relative to $\theta'$. Let $f$ be a continuous function on
$G(\mathbb{A})$, compactly supported modulo the center and transforming under
the center by a unitary character. Set
\begin{align*}
&\mathcal{R}^*:=\big\{\mathbf{s}=(s_1, s_2)\in\mathbb{C}^2:\ |\Re(s_1)|<1/(n+1),\ |\Re(s_2)|<1/(n+1)\big\},\\
&\mathcal{R}:=\big\{\mathbf{s}=(s_1,s_2)\in\mathbb{C}^2:\  \Re(s_1)>-1/(n+1),\ \Re(s_2)>-1/(n+1)\big\}.
\end{align*}

\subsection{\texorpdfstring{The Spectral Side
$J_{\Spec}^{\Reg,\heartsuit}(f,\mathbf{s})$}{The Spectral Side}}\label{sec2.1.}
For $\Re(s_1),\Re(s_2)\gg1$, consider
\begin{align*}
J^{\Reg}(f,\mathbf{s}):=\iint J_{\Kuz}\left(\begin{pmatrix}
		x&\\
		&1
	\end{pmatrix},\begin{pmatrix}
		y&\\
		&1
	\end{pmatrix}\right)W_{\phi_1'}'(x)\overline{W_{\phi_2'}'(y)}|\det x|^{s_1}|\det y|^{s_2}dxdy,
\end{align*}
where both integrals are over $N'(\mathbb{Q})\backslash G'(\mathbb{A})$ and
\begin{align*}
	J_{\Kuz}(g_1,g_2):=\int_{[N]}\int_{[N]}\K(u_1g_1,u_2g_2)\theta(u_1)\overline{\theta}(u_2)du_1du_2.
\end{align*}
Here $\K=\K^f$ is the kernel associated with $f$; we suppress the dependence
on $\phi_1'$ and $\phi_2'$ throughout. Expanding $\K$ spectrally and
interchanging the integrals gives the absolutely convergent expression
\begin{equation}\label{64}
	J_{\Spec}^{\Reg}(f,\mathbf{s})=\int_{\widehat{G(\mathbb{A})}_{\gen}}\sum_{\phi\in\mathfrak{B}_{\pi}}\Psi(s_1,\pi(f)W_{\phi},W_{\phi_1'}')\Psi(s_2,\widetilde{W}_{\phi},\widetilde{W}_{\phi_2'}')d\mu_{\pi},
\end{equation}
where $\widehat{G(\mathbb{A})}_{\gen}$ denotes the generic automorphic
spectrum with the prescribed central character, $d\mu_\pi$ is the automorphic
Plancherel measure, and
\begin{equation}\label{whittaker}
\Psi(s,W_{\phi},W_{\phi_i'}')=\int_{N'(\mathbb{A})\backslash G'(\mathbb{A})}W_{\phi}\left(\begin{pmatrix}
		x\\
		&1
	\end{pmatrix}\right)W_{\phi_i'}'(x)|\det x|^sdx,\ \ \Re(s)\gg 1.
\end{equation}
is the Rankin--Selberg period in Whittaker form. Here
$W_\phi(g):=\int_{[N]}\phi(ug)\overline{\theta(u)}\,du$, and
$\mathfrak{B}_\pi$ is a basis compatible with the spectral decomposition; see
\textsection\ref{sec5.2.1}. Absolute convergence in this right half-plane is
proved in \cite[\textsection 6.1]{Yan22}.

By Rankin--Selberg theory \cite{JPSS83,Jac09}, each period admits a
meromorphic continuation to $\mathbb{C}$ and is a holomorphic multiple of
$\Lambda(s+1/2,\pi\times\pi_i')$. Thus
\begin{align*}
P(s,\phi,\phi_i'):=\Lambda(s+1/2,\pi\times\pi_i')^{-1}\Psi(s,W_\phi,W_{\phi_i'}')
\end{align*}
is entire. Put
$\mathcal{P}(f,\mathbf{s},\phi):=
P(s_1,\pi(f)\phi,\phi_1')
\overline{P(\overline{s_2},\phi,\phi_2')}$ and define
\begin{align*}
J_{\Spec}^{\Reg,\heartsuit}(f,\mathbf{s})=\int_{\widehat{G(\mathbb{A})}_{\gen}}\sum_{\phi\in\mathfrak{B}_{\pi}}
\Lambda(s_1+1/2,\pi\times\pi_1')\Lambda(s_2+1/2,\widetilde{\pi}\times\widetilde{\pi}_2')\mathcal{P}(f,\mathbf{s},\phi)
	d\mu_{\pi}.
\end{align*}

In the initial region of convergence this agrees with
$J_{\Spec}^{\Reg}(f,\mathbf{s})$, with the Rankin--Selberg periods understood
by meromorphic continuation. Contour shifting gives a meromorphic continuation
$\widetilde{J}_{\Spec}^{\Reg}(f,\mathbf{s})$ to $\mathcal{R}$ and, on
$\mathcal{R}^*$, the relation
\begin{equation}\label{59}
\widetilde{J}_{\Spec}^{\Reg}(f,\mathbf{s})=J_{\Spec}^{\Reg,\heartsuit}(f,\mathbf{s})-\mathcal{G}_{\chi}(\mathbf{s},\phi_1',\phi_2')+\widetilde{\mathcal{G}}_{\chi}(\mathbf{s},\phi_1',\phi_2'),
\end{equation}
where $\mathbf{s}=(s_1,s_2)\in \mathcal{R}^*\subset \mathcal{R},$ and 
\begin{equation}\label{183}
	\mathcal{G}_{\chi}(\mathbf{s},\phi_1',\phi_2'):=\sum_{\phi\in\mathfrak{B}_{Q,\chi}}\underset{\lambda=n(s_1-1/2)}{\Res}\Psi(s_1,\mathcal{I}(\lambda,f)W_{\phi_{\lambda}},W'_{\phi_1'})\overline{\Psi(\overline{s_2},W_{\phi_{\lambda}},W'_{\phi_2'})},
\end{equation}
with $\phi_{\lambda}:=\phi\otimes|\cdot|^{\lambda}$, and similarly
\begin{equation}\label{198}
	\widetilde{\mathcal{G}}_{\chi}(\mathbf{s},\phi_1',\phi_2'):=\sum_{\phi\in\mathfrak{B}_{Q,\chi}}\underset{\lambda=n(1/2-s_2)}{\Res}\Psi(s_1,\mathcal{I}(\lambda,f)W_{\phi_{\lambda}},W'_{\phi_1'})\overline{\Psi(\overline{s_2},W_{\phi_{\lambda}},W'_{\phi_2'})}.
\end{equation}

The right-hand side of \eqref{59} converges absolutely on $\mathcal{R}^*$;
see \cite[\textsection 6]{Yan22}. Equivalently, \eqref{64}, with each
Rankin--Selberg period interpreted by meromorphic continuation, converges
absolutely there. The correction terms in \eqref{59} account for the residues
encountered when the spectral contours are shifted. Since our applications
take $\mathbf{s}=\mathbf{0}$, we record their geometric descriptions.

Let $\eta=(0,\ldots,0,1)\in\mathbb{Q}^n$. Define
$\widehat{E}_1^{(2)}(x,s;y)$ by
\begin{equation}\label{276}
	\sum_{\delta \in P'(\mathbb{Q})\backslash G'(\mathbb{Q})}\int_{Z'(\mathbb{A})}\int_{N_P(\mathbb{A})}\int_{\mathbb{A}^n}f\left(u(\mathbf{x})n\iota(y)\right)\psi(\eta \delta zx\mathfrak{u})d\mathbf{x}dn|\det zx|^{s}d^{\times}z,
\end{equation}
where we write $n=\begin{pmatrix}
	I_n&\mathfrak{u}\\
	&1
\end{pmatrix}\in N_P(\mathbb{A})$ and $u(\mathbf{x})=\begin{pmatrix}
	I_n&\\
	\mathbf{x}&1
\end{pmatrix},$ $	\mathbf{x}\in M_{1,n}(\mathbb{A})\simeq \mathbb{A}^n$.

This series converges absolutely for $\Re(s)>1$, satisfies a functional
equation, and continues meromorphically to $\mathbb{C}$. The same holds for
$\widehat{E}_1^{(1)}(x,s;y)$, defined by
\begin{equation}\label{60}
\sum_{\delta \in P'(\mathbb{Q})\backslash G'(\mathbb{Q})}\int_{Z'(\mathbb{A})}\int_{N_P(\mathbb{A})}\int_{\mathbb{A}^n}f\left(nu(\mathbf{x})\iota(y)\right)\psi(\eta \delta zx\mathfrak{u})d\mathbf{x}dn|\det zx|^{s}d^{\times}z
\end{equation}
admits the same analytic properties. 

\begin{prop}[{\cite[Theorem 8.1]{Yan22}}]\label{thm49}
For $\Re(s_1),\, \Re(s_2)\gg 1$, the residual terms have the geometric expressions
\begin{align*}
&\mathcal{G}_{\chi}(\mathbf{s},\phi_1',\phi_2')=
-\int_{G'(\mathbb{A})}\int_{[\overline{G'}]}\phi_1'(x)
\overline{\phi_2'(xy)}\widehat{E}_1^{(2)}(x,s_1+s_2;y)
\,dx\,|\det y|^{s_2}dy,\\
&\widetilde{\mathcal{G}}_{\chi}(\mathbf{s},\phi_1',\phi_2')=
\int_{G'(\mathbb{A})}\int_{[\overline{G'}]}\phi_1'(x)
\overline{\phi_2'(xy)}\widehat{E}_1^{(1)}(x,s_1+s_2;y)
\,dx\,|\det y|^{s_2}dy.
\end{align*}
Both sides continue meromorphically to $\mathcal{R}$.
\end{prop}

\subsection{\texorpdfstring{The Geometric Side
$J^{\Reg,\heartsuit}_{\Geo}(f,\mathbf{s})$}{The Geometric Side}}
Inserting the geometric expansion
\begin{align*}
\K(g_1,g_2)=\sum_{\gamma\in Z(\mathbb{Q})\backslash G(\mathbb{Q})}f(g_1^{-1}\gamma g_2),\ \ g_1,\, g_2\in G(\mathbb{A})
\end{align*}
into $J^{\Reg}(f,\mathbf{s})$ gives an elementary geometric expansion in a
right half-plane. Its direct form involves higher-rank Kloosterman sums.
Instead, \cite[\textsection3--\textsection5]{Yan22} reorganizes it as a
linear combination of automorphic periods admitting a meromorphic
continuation $\widetilde{J}^{\Reg}_{\Geo}(f,\mathbf{s})$ to $\mathcal{R}$.
The meromorphic relative trace formula is
\begin{equation}\label{rtf2.8}
	\widetilde{J}^{\Reg}_{\Spec}(f,\mathbf{s})=\widetilde{J}^{\Reg}_{\Geo}(f,\mathbf{s})
\end{equation}
	as an identity of meromorphic functions. In view of \eqref{59}, define
\begin{align*}
	J^{\Reg,\heartsuit}_{\Geo}(f,\mathbf{s}):=\widetilde{J}^{\Reg}_{\Geo}(f,\mathbf{s})+\mathcal{G}_{\chi}(\mathbf{s},\phi_1',\phi_2')-\widetilde{\mathcal{G}}_{\chi}(\mathbf{s},\phi_1',\phi_2'),\ \ \mathbf{s}\in\mathcal{R}^*.
\end{align*}

Although the correction terms originate spectrally, Proposition \ref{thm49}
expresses them geometrically, giving the decomposition
\begin{equation}\label{63}
J^{\Reg,\heartsuit}_{\Geo}(f,\mathbf{s})
=J^{\Reg}_{\Geo,\sm}(f,\mathbf{s})
+J^{\Reg}_{\Geo,\du}(f,\mathbf{s})
-J^{\Reg,\RNum{1}}_{\Geo,\bi}(f,\mathbf{s})
+J^{\Reg,\RNum{2}}_{\Geo,\bi}(f,\mathbf{s}).
\end{equation}
We describe its four terms below.

\subsubsection{The small-cell and dual terms}\label{2.2.1}
For $\Re(s)>1$, define
\begin{align*}
&E(s,x;\check{f}_P,y):=
	\sum_{\delta\in P'(\mathbb{Q})\backslash G'(\mathbb{Q})}\int_{Z'(\mathbb{A})}\int_{\mathbb{A}^n}f\left(\begin{pmatrix}
		y&\mathfrak{u}\\
		&1
	\end{pmatrix}\right)\psi(\eta z_1\delta x\mathbf{u})|\det z_1x|^{s}d\mathbf{u}d^{\times}z_1,\\
&E(s,x;f,y):=\sum_{\delta\in P'(\mathbb{Q})\backslash G'(\mathbb{Q})}\int_{Z'(\mathbb{A})}f(u(\eta z_1\delta x) \iota(y)) |\det z_1x|^{s}d^{\times}z_1,
\end{align*}
Both series continue meromorphically to $\mathbb{C}$ and satisfy functional
equations. For $\Re(s_1),\Re(s_2)>1/2$, set
\begin{align*}
	J^{\Reg}_{\Geo,\sm}(f,\mathbf{s}):=&\int_{G'(\mathbb{A})}\int_{[\overline{G'}]}\phi'(x)\overline{\phi'(xy)}E(1+s_1+s_2,x;\check{f}_P,y)dx|\det y|^{s_2}dy,\\
	J^{\Reg}_{\Geo,\du}(f,\mathbf{s}):=&\int_{G'(\mathbb{A})}\int_{[\overline{G'}]}\phi'(x)\overline{\phi'(xy)}E(s_1+s_2,x;f,y)dx|\det y|^{s_2}dy.
\end{align*}

The first term arises from the small Bruhat cell; the second is its dual under
Poisson summation. They continue meromorphically, with possible simple poles
at $s_1+s_2\in\{-1,0\}$ and $s_1+s_2\in\{0,1\}$, respectively. %Their residues at $s_1+s_2=0$ cancel:

\subsubsection{The first big-cell term}\label{2.2.2}
For $\Re(s)>1$, define
\begin{align*}
&E_1^{(2)}(x,s;y):=
	\sum_{\delta\in P'(\mathbb{Q})\backslash G'(\mathbb{Q})}\int_{Z'(\mathbb{A})}\int_{\mathbb{A}^n}f\left(\begin{pmatrix}
		I_n&\\
		\eta \delta x&1
	\end{pmatrix}\begin{pmatrix}
		y&u\\
		&1
	\end{pmatrix}\right)|\det zx|^{s}dud^{\times}z,\\
&E_1^{(1)}(x,s;y):=\sum_{\delta\in P'(\mathbb{Q})\backslash G'(\mathbb{Q})}\int_{Z'(\mathbb{A})}\int_{\mathbb{A}^n}f\left(\begin{pmatrix}
		I_n&u\\
		&1
	\end{pmatrix}\begin{pmatrix}
		y&\\
		\eta \delta xy&1
	\end{pmatrix}\right)|\det zx|^{s}dud^{\times}z.
\end{align*}

These series continue meromorphically to $\mathbb{C}$ and satisfy functional
equations relating them, after transposing $f$, to the series in \eqref{276}
and \eqref{60}.

Let $\Re(s_1), \Re(s_2)>0$. Write $s'=s_1+s_2+1$. Define 
\begin{align*}
	&\mathcal{F}_{1}^{(2)}J^{\bi}_{\Geo}(f,\mathbf{s})=\int_{[\overline{G'}]}\int_{{G'}(\mathbb{A})}\phi_1'(x)\overline{\phi_2'(xy)}E_1^{(2)}(x,s';y)|\det x|^{s'}|\det y|^{s_2}dydx,\\
	&\mathcal{F}_{1}^{(1)}J^{\bi}_{\Geo}(f,\mathbf{s})=\int_{[\overline{G'}]}\int_{{G'}(\mathbb{A})}\phi_1'(x)\overline{\phi_2'(xy)}E_1^{(1)}(x,s';y)|\det x|^{s'}|\det y|^{s_2}dydx.
\end{align*}

Define the first big-cell contribution by
\begin{equation}\label{62} 
J^{\Reg,\RNum{1}}_{\Geo,\bi}(f,\mathbf{s})=\sum_{j=1}^2\mathcal{F}_{1}^{(j)}J^{\bi}_{\Geo}(f,\mathbf{s})-\mathcal{G}_{\chi}(\mathbf{s},\phi_1',\phi_2')+\widetilde{\mathcal{G}}_{\chi}(\mathbf{s},\phi_1',\phi_2').
\end{equation}

The first two terms converge absolutely for $\Re(s_1+s_2)>0$ and continue
meromorphically to $\mathcal{R}$ \cite[\textsection 5]{Yan22}. Proposition
\ref{thm49} therefore gives the same continuation for
$J^{\Reg,\RNum{1}}_{\Geo,\bi}(f,\mathbf{s})$.

\subsubsection{The second big-cell term}\label{2.2.3}
For $\Re(s_1),\Re(s_2)>0$, define
\begin{align*}
	J^{\Reg,\RNum{2}}_{\Geo,\bi}(f,\mathbf{s})=&\sum_{\substack{(t,\boldsymbol{\xi})\in \mathbb{Q}\oplus \mathbb{Q}^{n-1}\\ (t,\boldsymbol{\xi})\neq (0,\mathbf{0}),\ t\neq1}}\int_{{G'}(\mathbb{A})}\int_{P_0'(\mathbb{Q})\backslash {G'}(\mathbb{A})}f\left(\iota(x)^{-1}\begin{pmatrix}
		I_{n-1}&&\boldsymbol{\xi}\\
		&1&t\\
		&1&1
	\end{pmatrix}\iota(xy)\right)\\
	&\qquad \qquad \qquad\qquad \phi'_1(x)\overline{\phi'_2(xy)}|\det x|^{s_1+s_2}|\det y|^{s_2}dxdy.
\end{align*}

The excluded zero orbit gives $J^{\Reg}_{\Geo,\du}(f,\mathbf{s})$. The
remaining sum is entire:

\begin{prop}[{\cite[Theorem 6.6]{Yan22}}]
The integral $J^{\Reg,\RNum{2}}_{\Geo,\bi}(f,\mathbf{s})$ converges
absolutely for every $\mathbf{s}\in\mathbb{C}^2$.
\end{prop}

\subsection{The Relative Trace Formula}
We summarize the form of the relative trace formula needed below.
\begin{thmx}\label{C}
Let $\mathbf{s}=(s_1,s_2)\in\mathcal{R}^*$. Then the spectral expansion
\begin{equation}\label{67}
		J_{\Spec}^{\Reg,\heartsuit}(f,\mathbf{s})=\int_{\widehat{G(\mathbb{A})}_{\gen}}\sum_{\phi\in\mathfrak{B}_{\pi}}\Psi(s_1,\pi(f)W_{\phi},W_{\phi_1'}')\Psi(s_2,\widetilde{W}_{\phi},\widetilde{W}_{\phi_2'}')d\mu_{\pi},
	\end{equation}
converges absolutely, and
	\begin{equation}\label{70}
		J_{\Spec}^{\Reg,\heartsuit}(f,\mathbf{s})=J_{\Geo}^{\Reg,\heartsuit}(f,\mathbf{s}),
	\end{equation}
holds as an identity of holomorphic functions on $\mathcal{R}^*$, with the
geometric side given by \eqref{63}.
\end{thmx}
\begin{remark}
In \eqref{67}, each Rankin--Selberg period is understood via its meromorphic
continuation
$\Psi(s,W_\phi,W_{\phi_i'}')=
\Lambda(s+1/2,\pi\times\pi_i')P(s,\phi,\phi_i')$.
\end{remark}

We apply \eqref{70} at $\mathbf{s}=\mathbf{0}$ and amplify it through the
choice of $f$; see \textsection\ref{sec3.7}.

\section{Test Functions and Amplification}\label{11.2} 
\subsection{Intrinsic Data}\label{3.1}
We now fix the data and notation used to construct the amplified test function
and prove Theorem \ref{B}.

\subsubsection{The automorphic weight}\label{sec3.1.2}
Let $\pi'=\otimes_{p\leq \infty}\pi_p'$ be a fixed unitary cuspidal representation of $G'(\mathbb{A})$ with central character $\omega'=\otimes_{p\leq\infty}\omega_p'$ and level $M'$. After twisting by a Dirichlet character (see \textsection\ref{seclevel}), we may assume that $\omega'\neq \textbf{1}$. Let $M''>1$ be the arithmetic conductor of $\omega'$. Then $M''\mid M'$.  
\begin{comment}
The assumption $\omega'\neq \textbf{1}$ is introduced to facilitate the construction of a local test function at $p\mid M''$ that satisfies $\mathcal{F}_{1}^{(\ell)}J^{\bi}_{\Geo}(f,\mathbf{s})\equiv 0$ for $\ell\in \{1,2\}$ and for all $\mathbf{s}\in\mathbb{C}^2$. We will also introduce another local test function at a fixed place to annahinalte $\mathcal{G}_{\chi}(\mathbf{s},\phi_1',\phi_2')$ and $\widetilde{\mathcal{G}}_{\chi}(\mathbf{s},\phi_1',\phi_2')$ simultaneously. Then by \eqref{62} we will have $J^{\Reg,\RNum{1}}_{\Geo,\bi}(f,\mathbf{s})\equiv 0$. 

While it is true that the contribution from $J^{\Reg,\RNum{1}}_{\Geo,\bi}(f,\mathbf{s})$ is overshadowed by that of $J^{\Reg,\RNum{2}}_{\Geo,\bi}(f,\mathbf{s})$, employing this particular choice of test function simplifies the calculations without compromising the main results.
\end{comment}

\subsubsection{Ramanujan parameters}\label{sec3.1.3}
Let $\vartheta_p$ be an admissible Ramanujan bound for the Satake parameters of $\pi_p'$. By \cite{KS03}, one has $\vartheta_p\leq \frac{1}{2}-\frac{2}{n(n+1)+2}$ for $p\leq \infty$.

\subsubsection{Uniform parameter growth}\label{sec3.14}
Let $\pi=\otimes_{p\leq \infty}\pi_p$ be a fixed unitary automorphic representation of $G(\mathbb{A})$ with central character $\omega$ and level $M$. By the Langlands classification, $\pi$ is the Langlands quotient of $\Ind_{Q}^{G}\pi_1\otimes\cdots\otimes\pi_m,$ where $Q$ is a parabolic subgroup of $G$ and $\pi_1\otimes\cdots\otimes\pi_m$ is a discrete automorphic representation of its Levi subgroup. In particular, if $\pi$ is cuspidal, then $Q=G,$ $m=1$ and $\pi_1=\pi$.

Let $T\geq 1$ and $C_{\infty}>c_{\infty}>0$. We assume that $\pi$ has uniform parameter growth of size $(T;c_{\infty},C_{\infty})$ in the sense of \textsection\ref{subc.}; that is,
\begin{equation*}
	c_{\infty}T\leq |\lambda_{\pi_\infty,j}|\leq C_{\infty}T,\ \ 1\leq j\leq n+1.\tag{\ref{1.6}}
\end{equation*}

\subsubsection{Level structure}\label{seclevel}
For $p<\infty$, define the Hecke congruence subgroup
\begin{equation}\label{equ3.4}
K_p(M):=\bigg\{\begin{pmatrix}
A'&\mathfrak{b}\\
\mathfrak{c}& d
\end{pmatrix}\in G(\mathbb{Z}_p):\ \mathfrak{c}\in p^{e_p(M)}M_{1,n}(\mathbb{Z}_p)\bigg\},
\end{equation}
where $e_p(M)$ is the valuation of $M$ at $p$. Thus $K_p(M)=K_p:=G(\mathbb{Z}_p)$ whenever $p\nmid M$. Similarly, define
\begin{equation}\label{3.4.}
K_p'(M'):=\bigg\{\begin{pmatrix}
A'&\mathfrak{b}\\
\mathfrak{c}& d
\end{pmatrix}\in G'(\mathbb{Z}_p):\ \mathfrak{c}\in p^{e_p(M')}M_{1,n-1}(\mathbb{Z}_p)\bigg\}. 
\end{equation}

We do not assume that $(M,M')=1$. Indeed, if
$\omega'=\mathbf{1}$, choose a fixed prime $q>n+1$ with $q\nmid MM'$ and
a Dirichlet character $\chi$ modulo $q$ of order $q-1$, and replace $(\pi,\pi')$ by $(\pi\otimes\chi^{-1},\,\pi'\otimes\chi)$. This leaves $L(s,\pi\times\pi')$ unchanged. Moreover, $\chi^n$ is nontrivial,
and hence primitive modulo $q$, so the twisted representations have levels
$Mq^{n+1}$ and $M'q^n$, respectively. Thus the twist makes the central
character of $\pi'$ nontrivial while introducing only a fixed common
ramified prime.

\subsubsection{The auxiliary vanishing prime}
Fix a sufficiently large prime $p_{\circ}\nmid MM'$, called the
\emph{auxiliary vanishing prime}, at which
we choose a local test function annihilating both residual terms
$\mathcal{G}_{\chi}(\mathbf{s},\phi_1',\phi_2')$ and
$\widetilde{\mathcal{G}}_{\chi}(\mathbf{s},\phi_1',\phi_2')$.
Since $p_{\circ}$ is fixed throughout, all implied constants may depend on it.

\subsubsection{Deformation parameters}\label{specpara}
If $m=1$, set $\mu=0$. If $m\geq2$, let
$\mathbf c=(c_1,\ldots,c_m)$, where $c_j=2^{-j}$ for $j<m$ and
$c_m=-\sum_{j<m}c_j$, and set
\begin{align*}
\mu=it_\mu\mathbf c,
\qquad
\pi_\mu:=\pi_1|\cdot|^{\mu_1}\boxplus\cdots\boxplus
\pi_m|\cdot|^{\mu_m}.
\end{align*}
Here $t_\mu=\exp(-O(\sqrt{\log T}))$ is chosen in
\textsection\ref{sec5.2} by a zero-avoidance argument. In particular,
$t_\mu=T^{-o(1)}$ and the coordinates $\mu_j$ are distinct.

This deformation is needed when $\pi$ belongs to the noncuspidal spectrum; see
Corollary \ref{cor1.2} in \textsection\ref{subc}.

%Denote by $\omega=\omega_0\widetilde{\chi}^{n+1}=\otimes_{p\leq \infty}\omega_p,$ where $\widetilde{\chi}$ is the idele class character induced from $\chi$. 

%Define $I_p(q)$ to be the subgroup 
%$$
%\Big\{g=(g_{i,j})\in G(\mathbb{Z}_p):\ e_p(\det g -1)\geq e_p(q),\ e_p(g_{i,j})\geq e_p(q),\ 1\leq j<i\leq n+1\Big\}.
%$$ 
%Let $I_p'(q)\subseteq G'(\mathbb{Z}_p)$ such that $\iota(I_p'(q))=I_p(q)\cap \iota(G'(\mathbb{Z}_p))$.
\subsubsection{Amplification parameters}\label{sec3.1.6}

Let $L\gg1$ satisfy $L\ll(TM)^{O(1)}$, and let
\begin{align*}
\mathcal{L}
\subseteq
\left\{p\ \text{prime}: L<p\leq 2L,\quad p\nmid p_{\circ}MM'\right\}.
\end{align*}
For each $p\in\mathcal{L}$, choose an integer
$1\leq\ell_p\leq n+1$ and a complex number $\alpha_p$, and write
\begin{align*}
\boldsymbol{\ell}:=(\ell_p)_{p\in\mathcal{L}},
\qquad
\boldsymbol{\alpha}:=(\alpha_p)_{p\in\mathcal{L}}.
\end{align*}

\subsubsection{Other notation}\label{3.1.6}
For a function $h$ on $G(\mathbb{A})$ or $G(\mathbb{Q}_p),$ $p\leq \infty,$ set
\begin{align*}
h^*(g)=\overline{h(g^{-1})},\quad (h*h^*)(g)=\int h(gg'^{-1})h^*(g')dg',
\end{align*} 
where $g'$ ranges over the domain of $h$. For a set $X$, let $\textbf{1}_X$ denote its characteristic function.

\subsection{Construction of Test Functions: The Archimedean Place}\label{11.1.1}

%\subsubsection{Analytic Test Functions at $p=\infty$}
We use the test functions ${f}_{\infty}$ constructed in \cite[\textsection 1.5.2 and \textsection 14]{Nel23}, following the formulation of \cite[\textsection 1.10]{NV21}; see also \cite[Part 2]{Nel21}.

\subsubsection{Construction of ${f}_{\infty}$}\label{3.2.1}
Let $\mathfrak{g}$ (resp. $\mathfrak{g}'$) be the Lie algebra of $G(\mathbb{R})$ (resp. $G'(\mathbb{R})$), with imaginary dual $\hat{\mathfrak{g}}$ (resp. $\hat{\mathfrak{g}}'$). We may choose $\tau\in\hat{\mathfrak{g}}$ such that, for $\tau':=\tau|_{\mathfrak g'}\in\hat{\mathfrak g}'$, the elements $T\tau$ and $T\tau'$ lie in the coadjoint orbits $\mathcal{O}_{\pi_\infty}$ and $\mathcal{O}_{\pi_\infty'}$, respectively. Thus $\tau$ remains in a fixed compact set as $T$ varies. Let $\tilde{f}^{\wedge}_{\infty}:\hat{\mathfrak{g}}\rightarrow\mathbb{C}$ be a smooth bump function concentrated on $\{T\tau+(\xi,\xi^{\bot}):\ \xi\ll T^{\frac{1}{2}+\varepsilon},\ \xi^{\bot}\ll T^{\varepsilon}\},$ where $\xi$ lies in the tangent space to the coadjoint orbit at $T\tau$ and $\xi^{\bot}$ lies in the normal space. Let $\tilde{f}_{\infty}^{\sharp}\in C_c^{\infty}(G(\mathbb{R}))$ be the pushforward of the Fourier transform of $\tilde{f}_{\infty}^{\wedge}$, and let $\tilde{f}_{\infty}^{\dag}$ be its truncation to the essential support:
\begin{equation}\label{245}
	\supp \tilde{f}_{\infty}^{\dagger}\subseteq \big\{g\in G(\mathbb{R}):\ g=I_{n+1}+O(T^{-\varepsilon}),\ \Ad^*(g)\tau=\tau+O(T^{-\frac{1}{2}+\varepsilon})\big\}.
\end{equation}

%Let $f_{\infty}(\cdot; \omega_{\infty})$ be defined by $f_{\infty}(g_{\infty}; \omega_{\infty})=\int_{Z(\mathbb{R})}\tilde{f}_{\infty}(g_{\infty}z_{\infty})\omega_{\infty}(z)d^{\times}z_{\infty}$. Let $f_{\infty}=f_{\infty}(\cdot; \omega_{\infty})*f_{\infty}(\cdot; \omega_{\infty})^*$. 

In the sense of \cite[\textsection 2.5]{NV21}, the operator $\pi_{\infty}(\tilde{f}_{\infty}^{\dagger})$ is approximately a rank-one projector whose image is spanned by a unit vector microlocalized at $\tau$. Define
\begin{equation}\label{eq3.2}
\tilde{f}_{\infty}(g):=\int_{Z(\mathbb{R})}\tilde{f}_{\infty}^{\dagger}(zg)\omega_{\infty}(z)d^{\times}z,\ \ g\in G(\mathbb{R}).
\end{equation} 

The Archimedean test function is the convolution $f_{\infty}=\tilde{f}_{\infty}*\tilde{f}_{\infty}^*$. %where $\tilde{f}_{\infty}^*(g):=\overline{\tilde{f}_{\infty}(g^{-1})}$; see \textsection\ref{3.1.6}.

\subsubsection{Transversality}
By construction (see \cite[(14.13)]{Nel21}),
\begin{equation}\label{250}
	\|\tilde{f}_{\infty}\|_{\infty}\ll_{\varepsilon} T^{\frac{n(n+1)}{2}+\varepsilon},
\end{equation}
where $\|\cdot \|_{\infty}$ is the sup-norm. For $g\in\overline{G}(\mathbb{R}),$ we may write 
\begin{align*}
g=\begin{pmatrix}
A&\mathfrak{b}\\
\transp{\mathfrak{c}}&d
\end{pmatrix},\ \ \ g^{-1}=\begin{pmatrix}
A'&\mathfrak{b}'\\
\transp{\mathfrak{c}}'&d'
\end{pmatrix},
\end{align*}
where $A, A'\in M_{n,n}(\mathbb{R}),$ $\mathfrak{b}, \mathfrak{b}', \mathfrak{c}, \mathfrak{c}'\in M_{n,1}(\mathbb{R}),$ and $d, d'\in \mathbb{R}$. Define 
\begin{equation}\label{dg}	
d_{G'}(g):=\begin{cases}
\min\big\{1, \|d^{-1}\mathfrak{b}\|+ \|d^{-1}\mathfrak{c}\|+ \|d'^{-1}\mathfrak{b}'\|+\|d'^{-1}\mathfrak{c}'\| \big\},&\ \text{if $dd'\neq 0$}\\
1,&\ \text{if $dd'=0$}.
\end{cases}
\end{equation}

\begin{prop}[Theorem 15.1 of \cite{Nel23}]\label{prop3.1}
Retain the preceding notation. There exists a fixed neighborhood $\mathcal{Z}$ of the identity in $Z'(\mathbb{R})$ such that, for $g$ in a sufficiently small neighborhood of $I_{n+1}$ in $\overline{G}(\mathbb{R})$ and sufficiently small $r>0$,
\begin{align*}
\Vol\left(\big\{z\in\mathcal{Z}:\ \dist(gz\tau, G'(\mathbb{R})\tau)\leq r\big\}\right)\ll \frac{r}{d_{G'}(g)}. 
\end{align*} 
Here $\dist(\cdots)$ is the infimum of $\|gz\tau-g'\tau\|$ over $g'\in G'(\mathbb{R})$, where $\|\cdot\|$ is a fixed norm on $\hat{\mathfrak{g}}$.
\end{prop}

We use Proposition \ref{prop3.1}, with $r=T^{-1/2+\varepsilon}$, to exploit the constraint $\Ad^*(g)\tau=\tau+O(T^{-\frac{1}{2}+\varepsilon})$ on the support of $\tilde{f}_{\infty}$.

\subsection{Construction of Test Functions: The Auxiliary Vanishing Prime}
At the auxiliary prime $p=p_{\circ}$, define
\begin{equation}\label{e3.7}
\tilde{f}_p(g):=\int_{Z(\mathbb{Q}_p)}\bigg[
\sum_{\beta\in(\mathbb{Z}/p\mathbb{Z})^\times}
\psi_p(\beta p^{-1})\mathbf{1}_{u_\beta^*K_p}(zg)
+\mathbf{1}_{\mathbf{d}K_p}(zg)\bigg]\omega_p(z)\,d^\times z,
\end{equation}
where $g\in G(\mathbb{Q}_p)$, and 
\begin{align*}
u_\beta^*
:=
\begin{pmatrix}
I_{n-1}&&\\
&1&\beta p^{-1}\\
&&1
\end{pmatrix},
\qquad
\mathbf{d}
:=
\begin{pmatrix}
I_{n-1}&&\\
&p^{-1}&\\
&&p
\end{pmatrix}.
\end{align*}

Set $f_p:=\tilde{f}_p*\tilde{f}_p^*$. Explicitly,
\begin{equation}\label{e3.8}
f_p(g)
=\sum_{i,j\in\{1,2\}}\int_{Z(\mathbb{Q}_p)}h_p^{(i,j)}(zg)\omega_p(z)\,d^\times z,
\end{equation}
where
\begin{align*}
h_p^{(1,1)}
&:=
\sum_{\beta,\beta'\in(\mathbb{Z}/p\mathbb{Z})^\times}
\psi_p(\beta p^{-1})
\overline{\psi_p(\beta' p^{-1})}
\mathbf{1}_{u_\beta^*K_pu_{-\beta'}^*},\quad h_p^{(2,2)}:=
\mathbf{1}_{\mathbf{d}K_p\mathbf{d}^{-1}},\\
h_p^{(1,2)}
&:=\sum_{\beta\in(\mathbb{Z}/p\mathbb{Z})^\times}
\psi_p(\beta p^{-1})
\mathbf{1}_{u_\beta^*K_p\mathbf{d}^{-1}},\quad 
h_p^{(2,1)}:=
\sum_{\beta'\in(\mathbb{Z}/p\mathbb{Z})^\times}
\overline{\psi_p(\beta' p^{-1})}
\mathbf{1}_{\mathbf{d}K_pu_{-\beta'}^*}.
\end{align*}

\subsection{Construction of Test Functions: Ramified Places}\label{sec3.3.}
\subsubsection{Test functions at $p\mid M$ and $p\nmid M'$}\label{11.1.3}
Define the local test function $f_p$ on $G(\mathbb{Q}_p),$ supported on $Z(\mathbb{Q}_p)\backslash K_p(M),$ by
\begin{align*}
f_p(zk)=\Vol(\overline{K_p(M)})^{-1}\omega_p(z)^{-1}\omega_p(E_{n+1,n+1}(k))^{-1},\ \ z\in Z(\mathbb{Q}_p),\ k\in K_p(M),
\end{align*}
where $\overline{K_p(M)}$ is the image of $K_p(M)$ in $\overline{G}(\mathbb{Q}_p),$ and $E_{n+1,n+1}(k)$ is the $(n+1,n+1)$-th entry of $k\in K_p(M)$.

\subsubsection{Test functions at $p\mid M'$}\label{11.1.4}
Set $m':=e_p(M')$ and $m'':=e_p(M'')$, and define the Gauss sum
\begin{align*}
G(\omega_p',\psi_p):=\sum_{\substack{\alpha_n\in (\mathbb{Z}/p^{m''}\mathbb{Z})^{\times}}}\psi_p(\alpha_np^{-m''})\omega_p'(\alpha_n).
\end{align*}
%Then $|G(\omega_p',\psi_p)|=p^{m''/2}$.

Let $\boldsymbol{\alpha}=(\alpha_1,\cdots,\alpha_{n-1},\alpha_n),$ where $\alpha_j\in \mathbb{Z}_p/p^{m'}\mathbb{Z}_p\simeq \mathbb{Z}/p^{m'}\mathbb{Z},$ $1\leq j<n$, and $\alpha_n\in (\mathbb{Z}_p/p^{m''}\mathbb{Z}_p)^{\times}\simeq (\mathbb{Z}/p^{m''}\mathbb{Z})^{\times}$. Set
\begin{align*}
u_{\boldsymbol{\alpha}}:=\begin{pmatrix}
I_n&\mathfrak{b}_{\boldsymbol{\alpha}}\\
&1
\end{pmatrix},\quad \mathfrak{b}_{\boldsymbol{\alpha}}=\transp{(\alpha_1p^{-m'},\cdots, \alpha_{n-1}p^{-m'},\alpha_n p^{-m''})}.
\end{align*}
Define the local test function by
\begin{equation}\label{3.4}
f_p(g):=(\tilde{f}_p*\tilde{f}_p^*)(g),\ \ \ g\in G(\mathbb{Q}_p),
\end{equation}
where $\tilde{f}_p^*(g):=\overline{\tilde{f}_p(g^{-1})}$ (see \textsection\ref{3.1.6}), and  
\begin{equation}\label{3.5}
\tilde{f}_p(g):=v_p\sum_{\substack{\alpha_j\in \mathbb{Z}/p^{m'}\mathbb{Z}\\ 1\leq j<n}}\sum_{\alpha_n\in (\mathbb{Z}/p^{m''}\mathbb{Z})^{\times}}\omega_p'(\alpha_n)\int_{Z(\mathbb{Q}_p)}\textbf{1}_{u_{\boldsymbol{\alpha}}K_p(M)}(zg)\omega_p(z)d^{\times}z.
\end{equation}
Here $v_p:=\Vol(K_p'(M'))^{-1}\Vol(K_p(M))^{-1}p^{-(n-1)m'}G(\omega_p',\psi_p)^{-1}$. Note that $|v_p|\asymp p^{ne_p(M)-m''/2}$.

\subsection{Construction of Test Functions: Amplification}\label{11.1.5}
Let $p\in\mathcal{L}$. For $\ell\geq 0,$ let $T_{p^{\ell}}=p^{-\frac{n\ell}{2}}\textbf{1}_{K_p\diag(p^\ell,I_n)K_p}$ be the Hecke operator at $p,$ where $K_p=G(\mathbb{Z}_p)$. Its adjoint is $T^*_{p^\ell}=p^{-\frac{n\ell}{2}}\textbf{1}_{K_p\diag(I_n,p^{-\ell})K_p}$.

Recall that $\boldsymbol{\ell}=(\ell_p)_{p\in\mathcal{L}}$. Define
\begin{align*}
f_p(g):=\int_{Z(\mathbb{Q}_p)}T_{p^{\ell_p}}(zg)\omega_p(z)d^{\times}z,\ \ g\in G(\mathbb{Q}_p).
\end{align*}
Its adjoint is given by
\begin{align*}
f_p^*(g):=\int_{Z(\mathbb{Q}_p)}T_{p^{\ell_p}}^*(zg)\omega_p(z)d^{\times}z. 
\end{align*}

For $p\in \mathcal{L}$ and $0\leq i\leq \ell_p,$ define
\begin{align*}
f_{p,i}(g):=p^{-ni}\int_{Z(\mathbb{Q}_p)}\textbf{1}_{K_{p}\diag(p^i,I_{n-1},p^{-i})K_{p}}(zg)\omega_p(z)d^{\times}z,\ \ g\in G(\mathbb{Q}_p).
\end{align*}
This normalization is the one occurring in the Hecke relation \eqref{3.66}.

\subsection{Construction of Test Functions: Remaining Places}\label{11.1.6}
For a prime $p\notin\mathcal{L}$ with $p\nmid MM'$, set
\begin{align*}
f_p(g):=\int_{Z(\mathbb{Q}_p)}\textbf{1}_{K_p}(zg)\omega_p(z)d^{\times}z,\ \ g\in G(\mathbb{Q}_p).
\end{align*}

\subsection{The Space of Test Functions}\label{testfunction}

\subsubsection{Decomposition of the test function}\label{sec3.6.1}

Retain the notation of
\textsection\ref{11.1.1}--\textsection\ref{11.1.6}.
For $p\in\mathcal{L}$, \cite[Lemma 4.4]{BM15} gives constants
$c_{p,i}\ll 1$ such that
\begin{equation}\label{3.66}
T_{p^{\ell_p}}*T_{p^{\ell_p}}^*=\sum_{i=0}^{\ell_p}
c_{p,i}p^{-ni}
\mathbf{1}_{K_p\diag(p^i,I_{n-1},p^{-i})K_p}.
\end{equation}

We now define the diagonal and off-diagonal components of the global test
function. For $p_0\in\mathcal{L}$ and $0\leq i\leq\ell_{p_0}$, let
\begin{align*}
f_{p_0,i}^{\diag}
:=
(f_\infty*f_\infty^*)
\otimes f_{p_0,i}
\otimes\bigotimes_{q<\infty:\ q\neq p_0}f_q.
\end{align*}
For distinct $p_1,p_2\in\mathcal{L}$, let
\begin{align*}
f_{p_1,p_2}^{\mathrm{off}}
:=
(f_\infty*f_\infty^*)
\otimes f_{p_1}
\otimes f_{p_2}^*
\otimes
\bigotimes_{q<\infty:\ q\nmid p_1p_2}f_q.
\end{align*}
The dependence of these functions on the amplification data
$\boldsymbol{\ell}=(\ell_p)_{p\in\mathcal{L}}$ is suppressed from the notation.

The amplified test function is
\begin{equation}\label{73}
f^{\mathrm{amp}}:=
\sum_{\substack{p_1,p_2\in\mathcal{L}\\p_1\neq p_2}}\alpha_{p_1}\overline{\alpha_{p_2}}\,
f_{p_1,p_2}^{\mathrm{off}}+\sum_{p_0\in\mathcal{L}}\sum_{i=0}^{\ell_{p_0}}
c_{p_0,i}|\alpha_{p_0}|^2f_{p_0,i}^{\diag}.
\end{equation}
In \textsection\ref{8.5.1}--\textsection\ref{sec10}, we estimate the orbital
integrals associated with the individual components
\begin{align*}
f\in
\left\{
f_{p_0,i}^{\diag},\,
f_{p_1,p_2}^{\mathrm{off}}
\right\}.
\end{align*}
These estimates are combined in \textsection\ref{proof} to bound the amplified
geometric side of Theorem \ref{C} for the test function \eqref{73}.

\subsubsection{Auxiliary notation}\label{3.6.2}

For an individual component $f$ of the amplified test function, define
\begin{equation}\label{61}
(\nu(f),\mathcal{N}_f)
:=
\begin{cases}
\big(p_1p_2,\,
p_1^{\ell_{p_1}/2}p_2^{\ell_{p_2}/2}\big),
& \text{if } f=f_{p_1,p_2}^{\mathrm{off}},\\[2mm]
(p_0,\, p_0^i),
& \text{if } f=f_{p_0,i}^{\diag}.
\end{cases}
\end{equation}
We suppress the dependence of $f$ on the levels, the auxiliary and
amplification primes, and the archimedean data.

Recall that $\tilde{f}_p$ was defined at $p=\infty$, $p=p_{\circ}$, and
$p\mid M'$ in \eqref{eq3.2}, \eqref{e3.7}, and \eqref{3.5}, respectively.
At the remaining finite primes, define
\begin{align*}
\tilde{f}_p(g):=
\begin{cases}
f_p(g), & p\nmid\nu(f),\\[1mm]
\displaystyle
\int_{Z(\mathbb{Q}_p)}
\mathbf{1}_{K_p}(zg)\omega_p(z)\,d^\times z,
& p\mid\nu(f).
\end{cases}
\end{align*}
Set $\tilde{f}:=\bigotimes_{p\leq\infty}\tilde{f}_p$. Then $\tilde{f}*\tilde{f}^{*}$ and $f$ have the same local components
away from the amplification primes $p\mid\nu(f)$. Thus
$\tilde{f}*\tilde{f}^{*}$ is the unamplified counterpart of $f$; it
will be used in the proof of Proposition \ref{thm6} in
\textsection\ref{sec5.2}.

\subsection{\texorpdfstring{Vanishing of
$J^{\Reg,\RNum{1}}_{\Geo,\bi}(f,\mathbf{s})$}{Vanishing of a Correction Term}}
\begin{lemma}\label{lem3.2}
Retain the notation of \textsection\ref{2.2.2}, and let $f$ be as in
\textsection\ref{testfunction}. Then
\begin{equation}\label{3.8}
J^{\Reg,\RNum{1}}_{\Geo,\bi}(f,\mathbf{s})\equiv 0,
\quad \mathbf{s}\in\mathcal{R}.
\end{equation}
\end{lemma}

\begin{proof}
Recall from \eqref{62} that
\begin{align*}
J^{\Reg,\RNum{1}}_{\Geo,\bi}(f,\mathbf{s})
=
\sum_{j=1}^2
\mathcal{F}_{1}^{(j)}J^{\bi}_{\Geo}(f,\mathbf{s})
-\mathcal{G}_{\chi}(\mathbf{s},\phi_1',\phi_2')
+\widetilde{\mathcal{G}}_{\chi}(\mathbf{s},\phi_1',\phi_2').
\end{align*}
We show that all four terms on the right-hand side vanish.

First choose a prime $p\mid M''$. By \eqref{3.4}, we have
\begin{equation}\label{3.13}
f_p(g)
=
v_p^2
\sum_{\boldsymbol{\alpha}}
\sum_{\boldsymbol{\beta}}
\omega_p'(\alpha_n)\omega_p'(\beta_n)
\int_{Z(\mathbb{Q}_p)}
\mathbf{1}_{u_{\boldsymbol{\alpha}}K_p
u_{\boldsymbol{\beta}}^{-1}}(zg)
\omega_p(z)\,d^\times z,
\end{equation}
where the ranges of $\boldsymbol{\alpha}$ and $\boldsymbol{\beta}$ are as in
\eqref{3.5}. Since $\omega_p'$ is ramified, character orthogonality gives
\begin{align*}
\sum_{\beta_n\in
(\mathbb{Z}/p^{m''}\mathbb{Z})^\times}
\omega_p'(\beta_n)=0.
\end{align*}
Consequently, for every $\mathbf{x}\in M_{1,n}(\mathbb{Q}_p)$ and
$y\in G'(\mathbb{Q}_p)$,
\begin{equation}\label{eq3.9}
\int_{N_P(\mathbb{Q}_p)}
f_p\bigl(u(\mathbf{x})n\iota(y)\bigr)\,dn
=
\int_{N_P(\mathbb{Q}_p)}
f_p\bigl(nu(\mathbf{x})\iota(y)\bigr)\,dn
=0.
\end{equation}
Inserting \eqref{eq3.9} into the defining series gives
\begin{align*}
E_1^{(1)}(x,s;y)=E_1^{(2)}(x,s;y)=0
\end{align*}
in their domain of absolute convergence. Unfolding and meromorphic
continuation therefore yield
\begin{align*}
\mathcal{F}_{1}^{(j)}J^{\bi}_{\Geo}(f,\mathbf{s})
\equiv 0,
\quad j=1,2,\quad \mathbf{s}\in\mathcal{R}.
\end{align*}

It remains to treat the residual terms at $p=p_{\circ}$. We use the
four-term decomposition \eqref{e3.8}. The identity
\begin{equation}\label{f2.22}
\begin{pmatrix}
1&\beta p^{-1}\\
&1
\end{pmatrix}
=
\begin{pmatrix}
1&\\
\beta^{-1}p&1
\end{pmatrix}
\begin{pmatrix}
p^{-1}&\\
&p
\end{pmatrix}
\begin{pmatrix}
p&\beta\\
-\beta^{-1}&0
\end{pmatrix},
\qquad
\beta\in(\mathbb{Z}/p\mathbb{Z})^\times,
\end{equation}
in the lower-right $2\times2$ block relates the $(i,1)$- and
$(i,2)$-terms. After a unipotent change of variables, the Ramanujan-sum
identity
\begin{align*}
\sum_{\beta\in(\mathbb{Z}/p\mathbb{Z})^\times}
\psi_p(\beta p^{-1})=-1
\end{align*}
therefore gives
\begin{align*}
\int_{M_{1,n}(\mathbb{Q}_p)}
\left[
h_p^{(i,1)}\bigl(nu(\mathbf{x})\iota(y)\bigr)
+
h_p^{(i,2)}\bigl(nu(\mathbf{x})\iota(y)\bigr)
\right]d\mathbf{x}=0,
\qquad i=1,2.
\end{align*}
The same argument on the other side gives
\begin{align*}
\int_{M_{1,n}(\mathbb{Q}_p)}
\left[
h_p^{(1,j)}\bigl(u(\mathbf{x})n\iota(y)\bigr)
+
h_p^{(2,j)}\bigl(u(\mathbf{x})n\iota(y)\bigr)
\right]d\mathbf{x}=0,
\qquad j=1,2.
\end{align*}
This is the same local calculation as in \cite[\S 3.4]{Yan23b}; the
remaining $n-1$ coordinates lie in $K_p$ and do not affect the argument.
Summing over $i,j$ and integrating over the center, we find that the local
integrals defining $\widehat E_1^{(1)}$ and $\widehat E_1^{(2)}$ vanish.
Consequently,
\begin{align*}
\mathcal{G}_{\chi}(\mathbf{s},\phi_1',\phi_2')
=
\widetilde{\mathcal{G}}_{\chi}(\mathbf{s},\phi_1',\phi_2')
=0
\end{align*}
in the domain of absolute convergence, and hence identically on
$\mathcal{R}$ by meromorphic continuation.

All terms on the right-hand side of \eqref{62} vanish identically, proving
\eqref{3.8}.
\end{proof}

\subsection{The Amplified Relative Trace Formula}\label{sec3.7}
Retain the local and global data chosen in
\textsection\ref{3.1}--\textsection\ref{testfunction}, including
$\mathcal{L},$ $\boldsymbol{\alpha},$ $\boldsymbol{\ell},$ and
$f\in \big\{f_{p_0,i}^{\diag},\,f_{p_1,p_2}^{\mathrm{off}}\big\}$.

\subsubsection{The spectral side}\label{ampsp}
In accordance with \eqref{73}, define
\begin{align*}
\mathcal{J}_{\Spec}^{\heartsuit}(\boldsymbol{\alpha},\boldsymbol{\ell})
:=\sum_{p_1\neq p_2\in\mathcal{L}}\alpha_{p_1}\overline{\alpha_{p_2}}
J_{\Spec}^{\Reg,\heartsuit}(f_{p_1,p_2}^{\mathrm{off}},\mathbf{0})+\sum_{p_0\in\mathcal{L}}\sum_{i=0}^{\ell_{p_0}}c_{p_0,i}|\alpha_{p_0}|^2
J_{\Spec}^{\Reg,\heartsuit}(f_{p_0,i}^{\diag},\mathbf{0}).
\end{align*}

\subsubsection{The geometric side}\label{sec3.7.2}
Similarly, define
\begin{align*}
\mathcal{J}_{\Geo}^{\heartsuit}(\boldsymbol{\alpha},\boldsymbol{\ell})
:=\sum_{p_1\neq p_2\in\mathcal{L}}\alpha_{p_1}\overline{\alpha_{p_2}}
J_{\Geo}^{\Reg,\heartsuit}(f_{p_1,p_2}^{\mathrm{off}},\mathbf{0})+\sum_{p_0\in\mathcal{L}}\sum_{i=0}^{\ell_{p_0}}c_{p_0,i}|\alpha_{p_0}|^2
J_{\Geo}^{\Reg,\heartsuit}(f_{p_0,i}^{\diag},\mathbf{0}),
\end{align*}
where $c_{p_0,i}$ is determined by the Hecke relation \eqref{3.66}.

\subsubsection{The amplified relative trace formula}\label{3.7.3}

Applying Theorem \ref{C} with $\mathbf{s}=\mathbf{0}$ and
$f=f^{\mathrm{amp}}$ as defined in \eqref{73}, we obtain the following
identity.

\begin{thmx}\label{thmD'}
Retain the preceding notation. Then
\begin{equation}\label{3.7}
\mathcal{J}_{\Spec}^{\heartsuit}
(\boldsymbol{\alpha},\boldsymbol{\ell})
=
\mathcal{J}_{\Geo}^{\heartsuit}
(\boldsymbol{\alpha},\boldsymbol{\ell}).
\end{equation}
\end{thmx}

\section{\texorpdfstring{Construction of Test Vectors in $\pi'$}
{Construction of Test Vectors}}\label{sec4.}
Recall that $\pi'=\otimes_{p\leq\infty}\pi_p'$ is a fixed unitary cuspidal
representation of $G'(\mathbb{A})$, with central character $\omega'$ and level
$M'$. We construct a $T$-dependent factorizable vector $\phi'\in\pi'$ for the
relative trace formula. Its Whittaker function is
$W'=\otimes_{p\leq\infty}W_p'$, where
\begin{itemize}
\item $W_p'$ is the normalized new Whittaker vector at every finite prime $p$;
\item $W_\infty'$ is the microlocalized Whittaker vector constructed below.
\end{itemize}

\subsection{The Archimedean Test Vector}\label{spec}

Although $\pi_\infty'$ is fixed, the vector $W_\infty'$ must vary with $T$:
it is microlocalized at the frequency determined by the spectral parameter of
$\pi_\infty$. This matching is essential for the required lower bound for the
local Rankin--Selberg integral, whereas a fixed smooth vector would have
negligible mass at the relevant frequency. We therefore adapt the construction
of \cite[\S\S 10--11]{Nel21} to the fixed generic representation
$\pi_\infty'$.

By the Langlands classification and genericity, followed by Casselman's
subrepresentation theorem, $\pi_\infty'$ embeds into a principal series
$\Ind_{B'(\mathbb{R})}^{G'(\mathbb{R})}\chi$. Write
\begin{align*}
\chi(\diag(a_1,\ldots,a_n))
=\prod_{j=1}^n\chi_j(a_j)|a_j|^{\lambda_j},
\end{align*}
where each $\chi_j$ is unitary and the exponents $\lambda_j$ are fixed. In
particular, $|\lambda_j|=o(T)$, which is the only uniformity required below.

\subsubsection{Functions on $N'\backslash G'$}\label{4.1.2}
Let $\overline{N'}$ be the unipotent subgroup opposite to $N'$, and let
$\theta$ denote the relevant nondegenerate unitary characters on $N'$ and
$\overline{N'}$, according to context. Fix
$\beta\in C_c^{\infty}(A'(\mathbb{R}))$ such that
$\int_{A'(\mathbb{R})}\beta(a)\overline{\chi}(a)\,d^{\times}a=1$, and define
\begin{align*}
\mathcal{J}[\theta,\beta](g'):=\begin{cases}
\delta_{B'}^{\frac{1}{2}}(a)\beta(a)\theta(u),&\  \text{if $g'=au\in A'(\mathbb{R})\overline{N'}(\mathbb{R})$,}\\
0, &\ \text{if $g'\notin A'(\mathbb{R})\overline{N'}(\mathbb{R})$}.
\end{cases}
\end{align*}

Choose $h_0\in C_c^{\infty}(\mathfrak{g}')$ as in
\cite[Corollary 10.25 and the proof of Proposition 10.26]{Nel21}, normalized
so that $\int h_0=1$, and set
$\theta_T(\exp X):=\theta(\exp(TX))$. Denote its logarithmic differential,
extended by zero to a complementary subspace of $\mathfrak g'$, by
$\theta^\dagger\in i(\mathfrak g')^*$.

Define a smooth measure $\gamma_0$ on $G'(\mathbb{R})$ by
\begin{equation}\label{75}
\int_{G'(\mathbb{R})}\Phi(g)\,d\gamma_0(g)
:=\int_{\mathfrak{g}'}h_0(X)\Phi(\exp(T^{-1/2}X))
\exp(-\langle T\theta^{\dagger},T^{-1/2}X\rangle)\,dX.
\end{equation}

Let $\gamma=\gamma_0^* * \gamma_0$. Define
\begin{equation}\label{65}
\mathcal{J}[\theta_T,\beta,\gamma](g'):=\int_{G'(\mathbb{R})}\gamma(g'')\mathcal{J}[\theta_T,\beta](g'g'')dg''
\end{equation}
and define $\mathcal{J}_T[\theta,\beta,\gamma]$ by replacing $\beta$ in
\eqref{65} with
\begin{align*}
L(T^{\rho^{\vee}})\beta(a)
:=\delta_{B'}^{-1/2}(T^{\rho^{\vee}})\beta(T^{\rho^{\vee}}a),
\end{align*}
where $\rho^\vee$ denotes the half-sum of the positive coroots of
$(G',B')$.

Define the Mellin component
\begin{align*}
\mathcal{M}_T[\theta,\beta,\gamma;\chi](g'):=\int_{A'(\mathbb{R})}\delta_{B'}^{\frac{1}{2}}(a)\mathcal{J}_T[\theta,\beta,\gamma](a^{-1}g')\chi(a)d^{\times}a,
\end{align*}
which is a smooth section of
$\Ind_{B'(\mathbb{R})}^{G'(\mathbb{R})}\chi$, and set
\begin{align*}
\mathcal{W}_T(g'):=\chi(T^{-\rho^{\vee}})
\int_{\overline{N'}(\mathbb{R})}
\mathcal{M}_T[\theta,\beta,\gamma;\chi](ug')
\overline{\theta}_T(u)du.
\end{align*}

Changing variables gives
\begin{align*}
\mathcal{M}_T[\theta,\beta,\gamma;\chi](g')
=\chi(T^{\rho^{\vee}})
\mathcal{M}[\theta_T,\beta,\gamma;\chi](g'),
\end{align*}
where
\begin{align*}
\mathcal{M}[\theta_T,\beta,\gamma;\chi](g'):=\int_{A'(\mathbb{R})}\delta_{B'}^{\frac{1}{2}}(a)\mathcal{J}[\theta_T,\beta,\gamma](a^{-1}g')\chi(a)d^{\times}a.
\end{align*}

By the normalization of $\beta$, unfolding the convolution
$\gamma=\gamma_0^**\gamma_0$ gives
\begin{align*}
\mathcal{W}_T(g')
=\int_{\overline{N'}(\mathbb{R})}
\mathcal{M}[\theta_T,\beta,\gamma_0;\chi](u)
\overline{\mathcal{M}[\theta_T,\beta,\gamma_0*g'^{-1};
\overline{\chi}^{-1}](u)}du.
\end{align*}

\subsubsection{Construction of $W_{\infty}'$ and $\phi'$}\label{sec4.1.3}
\begin{defn}\label{def4.2}
Let $w_{n-1}$ be the long Weyl element associated with $(G',B')$. For
$g'\in G'(\mathbb{R})$, set
\begin{align*}
W_{\infty}'(g'):=\chi(T^{-\rho^{\vee}})T^{\frac{n(n-1)}{8}}\int_{N'(\mathbb{R})}\mathcal{M}_T[\theta,\beta,\gamma;\chi](w_{n-1}ug')\overline{\theta}(u)du.
\end{align*}
Let $\phi'\in\pi'$ be the cusp form whose Whittaker function is
$\otimes_{p\leq\infty}W_p'$, where $W_p'$ is the normalized local new
Whittaker function for every $p<\infty$.
\end{defn}

The auxiliary function $\mathcal{W}_T$ is written in the
$\overline{\theta}_T$-Whittaker model, whereas $W_\infty'$ belongs to the
ordinary $\overline{\theta}$-Whittaker model. Weyl conjugation and the
rescaling identity of \cite[Lemma 10.3]{Nel21} pass between these two
realizations. The factor $T^{n(n-1)/8}$ is the normalization used when
returning to the ordinary Whittaker model.

\subsection{Archimedean Estimates}
The normalization of $\beta$ gives
\begin{equation}\label{74}
\mathcal{M}[\theta_T,\beta,\gamma_0;\chi](u)
=\int_{G'(\mathbb{R})}\gamma_0(g)
\delta_{B'}^{1/2}((ug)_{A'})
\theta_T((ug)_{\overline{N'}})
\chi((ug)_{A'})\,dg.
\end{equation}

\begin{lemma}\label{lem4.3}
Let $u\in\overline{N'}(\mathbb{R})$ satisfy
$u=I_n+O(T^{-1/2})$. Then
\begin{equation}\label{78}
\mathcal{M}[\theta_T,\beta,\gamma_0;\chi](u)
=\theta_T(u)+O(T^{-\varepsilon})
\end{equation}
for some fixed $\varepsilon>0$.
\end{lemma}
\begin{proof}
Insert \eqref{75} into \eqref{74} and use the local coordinates
$u\exp(T^{-1/2}X)$. The oscillatory factor in \eqref{75} cancels the linear
term of $\theta_T$, while the $A'$-component is $I_n+O(T^{-1/2})$.
Since the inducing exponents of $\chi$ are fixed, the remaining factors equal
$1+O(T^{-\varepsilon})$. The normalization $\int h_0=1$ now gives
\eqref{78}.
\end{proof}

\begin{lemma}\label{lem4.4}
Let $g'\in B'(\mathbb{R})$ satisfy
$g'=I_n+O(T^{-1/2-\varepsilon})$. Then
\begin{equation}\label{79}
\int_{\overline{N'}(\mathbb{R})}
\big|\mathcal{M}[\theta_T,\beta,
\gamma_0*g'^{-1}-\gamma_0;\overline{\chi}^{-1}](u)\big|^2du
\ll T^{-\frac{n(n-1)}4-2\varepsilon+o(1)}.
\end{equation}
\end{lemma}
\begin{proof}
The fundamental-theorem-of-calculus argument of
\cite[Proposition 6.53]{Nel21} contributes the factor $T^{-\varepsilon}$
arising from the translation by $g'^{-1}$. Applying this in each copy and
then using the $L^2$ estimate of \cite[Lemma 10.20]{Nel21} gives
\eqref{79}; the exponent $n(n-1)/4$ is half the dimension of
$\overline{N'}$.
\end{proof}

\begin{lemma}\label{lem4.5}
Let $g'\in B'(\mathbb{R})$ satisfy
$g'=I_n+O(T^{-1/2-\varepsilon})$. Then
\begin{equation}\label{81}
\mathcal{W}_T(g')
=\int_{\overline{N'}(\mathbb{R})}
\big|\mathcal{M}[\theta_T,\beta,\gamma_0;\chi](u)\big|^2du
+O\left(T^{-\frac{n(n-1)}4-2\varepsilon+o(1)}\right).
\end{equation}
In particular,
\begin{align*}
T^{-\frac{n(n-1)}4}
\ll \Re(\mathcal{W}_T(g'))
\ll T^{-\frac{n(n-1)}4+\varepsilon}.
\end{align*}
\end{lemma}
\begin{proof}
Unfolding the convolution $\gamma=\gamma_0^* * \gamma_0$ expresses
$\mathcal{W}_T(g')$ as the inner product of the Mellin transforms
associated with $\gamma_0$ and $\gamma_0*g'^{-1}$. Lemma \ref{lem4.4}
replaces the latter by the former with the error in \eqref{81}, and Lemma
\ref{lem4.3} gives the asserted bounds.
\end{proof}

\subsection{Norms}\label{sec4.5}
Let $H=\diag(\mathrm{GL}(n-1),1)$ and $N_H=N'\cap H$. Define
\begin{equation}\label{W_inf}
\|W_{\infty}'\|_2
:=\left(
\int_{N_H(\mathbb{R})\backslash H(\mathbb{R})}
|W_{\infty}'(h)|^2\,dh
\right)^{1/2}.
\end{equation} 

\begin{lemma}\label{lem4.6}
With the preceding choices,
\begin{equation}\label{83}
1\ll_{\pi_\infty'}\|W_{\infty}'\|_2
\ll_{\pi_\infty'}T^{o(1)},
\qquad
\langle \phi',\phi'\rangle
\asymp_{M'}L(1,\pi',\Ad)\|W_{\infty}'\|_2^2,
\end{equation}
\end{lemma}
\begin{proof}
The Kirillov-model Plancherel formula identifies the square of the first norm
with the $K_\infty'$-norm of the Mellin component, multiplied by
$T^{n(n-1)/4}$. Lemma 10.20 of \cite{Nel21} bounds the latter by
$T^{-n(n-1)/4+o(1)}$. Conversely, Lemma \ref{lem4.5}, evaluated at the
identity, gives the matching lower bound. This proves the first assertion.
The second follows from the standard Rankin--Selberg norm identity and the
normalization $W_p'(I_n)=1$ at every finite prime.
\end{proof}

\subsection{Siegel Domains and Bounds for Cusp Forms}\label{sec4.4}
Set
\begin{align*}
A_H(\mathbb{R})
:=\big\{\diag(a_1,\ldots,a_{n-1},1)\in A'(\mathbb{R}):
a_1\geq\cdots\geq a_{n-1}>0\big\}.
\end{align*}

By reduction theory \cite{Bor62}, there exists a compact subset $\Omega$ of
$P_0'(\mathbb{R})$ such that every nonnegative integrable function $h$ on
$P_0'(\mathbb{Q})\backslash \overline{G'}(\mathbb{A})$ satisfies
\begin{equation}\label{87}
\int_{P_0'(\mathbb{Q})\backslash \overline{G'}(\mathbb{A})}h(g)dg\leq \int_{A_H(\mathbb{R})}\int_{\Omega K_{\fin}'}h(ag)dg\delta_{B'}^{-1}(a)d^{\times}a,
\end{equation}
This is the parabolic analogue of \cite[Proposition 2.25]{Nel21}.
\begin{lemma}\label{lem4.7'}
There exists an absolute constant $C_0>0$ such that
\begin{equation}\label{88}
	\int_{\Omega K_{\fin}'}h(ag)dg\ll \delta_{B'}(a)\int_{P_0'(\mathbb{Q})\backslash \overline{G'}(\mathbb{A})}h(g)\textbf{1}_{C_0|\det a|_{\infty}\leq |\det g|\leq 2C_0|\det a|_{\infty}}dg,
\end{equation}
where the implied constant is absolute. 
\end{lemma}
\begin{proof}
The proof of \cite[Proposition 2.26]{Nel21}, applied up to its final
enlargement of the torus domain, gives
\begin{equation}\label{90}
\int_{\Omega K_{\fin}'}h(ag)dg\ll \delta_{B'}(a)\int_{[N']}\int_{a\Omega_{A'}}\int_{K'}h(ua'k)dk\delta_{B'}^{-1}(a')d^{\times}a'du,
\end{equation}
where $\Omega_{A'}\subset A_H(\mathbb{R})$ is the set of semisimple
components of elements of $\Omega$. Retaining the restriction
$a'\in a\Omega_{A'}$ gives the determinant constraint, and the right-hand side
of \eqref{90} is bounded by
\begin{align*}
\delta_{B'}(a)\int_{[N']}\int_{A_H(\mathbb{R})}\int_{K'}h(ua'k)\textbf{1}_{C_0|\det a|_{\infty}\leq |\det a'|\leq 2C_0|\det a|_{\infty}}dk\delta_{B'}^{-1}(a')d^{\times}a'du
\end{align*}
for a constant $C_0$ depending only on $\Omega$. This proves \eqref{88}.
\end{proof}

For $\mathfrak{c}=(c_1,\cdots,c_n)\in\mathbb{A}^n,$ let $\|\mathfrak{c}_{\infty}\|_{\infty}:=\sqrt{|c_{1,\infty}|_{\infty}^2+\cdots+|c_{n,\infty}|_{\infty}^2}$ and 
\begin{align*}
\|\mathfrak{c}\|:=\|\mathfrak{c}_{\infty}\|_{\infty}\cdot \prod_{p<\infty}\max\{|c_{1,p}|_p,\cdots,|c_{n,p}|_p\}.
\end{align*}

Fix $h_{\infty}\in C_c^{\infty}(\mathbb{R})$ such that
$0\leq h_{\infty}(t)\leq1$, with $h_{\infty}(t)=1$ for $1\leq t\leq2$
and $h_{\infty}(t)=0$ for $t\leq1/2$ or $t\geq3$.

Let $h^{\dagger}(\mathfrak{c}):=h_{\infty}^{\dagger}(\mathfrak{c}_{\infty})\prod_{p<\infty}h_p^{\dagger}(\mathfrak{c}_p),$ where $h_{\infty}^{\dagger}(\mathfrak{c}_{\infty})=h_{\infty}(\|\mathfrak{c}_{\infty}\|_{\infty}),$ and  $h_p^{\dagger}(\mathfrak{c}):=\mathbf{1}_{\mathbb{Z}_p}(c_{1,p})\cdots \mathbf{1}_{\mathbb{Z}_p}(c_{n,p})\textbf{1}_{\max\{|c_{1,p}|_p,\cdots,|c_{n,p}|_p\}=1}$ for $p<\infty$. 

Then $h^{\dagger}$ is a bounded, compactly supported Borel function on
$\mathbb{A}^n$. In particular, $h^{\dagger}(\mathfrak{c})=0$ unless
$1/2\leq \|\mathfrak{c}\|\leq 3$.

\begin{lemma}\label{lem4.8.}
Let $C_0>0$ be an absolute constant. Define
\begin{equation}\label{4.14.}
\mathcal{I}(a,y):=\int_{N'(\mathbb{A})\backslash G'(\mathbb{A})}\big|W'(gy)\big|^2\textbf{1}_{C_0|\det a|_{\infty}\leq |\det g|\leq 2^{n+1}C_0|\det a|_{\infty}}h^{\dagger}(\eta g)dg,
\end{equation}
where $\eta=(0,\cdots,0,1)\in \mathbb{Q}^n$. Then 
\begin{equation}\label{4.15.}
\mathcal{I}(a,y)\ll T^{\varepsilon}\cdot |\det a|_{\infty}^{-1}\cdot d(y)^2\cdot \langle\phi',\phi'\rangle,
\end{equation}
where the implied constant depends at most on $\varepsilon$, $n$, and $\pi'$,
and
\begin{equation}\label{dy}
d(y):=|\det y|_{\mathbb{A}}^{-1/2}\cdot \bigg[\int_{ K'}\|\eta ky^{-1}\|^{-n}dk\bigg]^{1/2}.
\end{equation}
\end{lemma}
\begin{proof}
After a change of variable $g\mapsto gy^{-1},$ we obtain 
\begin{align*}
\mathcal{I}(a,y)=\int_{N'(\mathbb{A})\backslash G'(\mathbb{A})}\big|W'(g)\big|^2\textbf{1}_{C_0\leq \frac{|\det g|}{|\det a|_{\infty}|\det y|_{\mathbb{A}}}\leq 2^{n+1}C_0}h^{\dagger}(\eta gy^{-1})dg.
\end{align*}

Write $g=zbk$ according to the Iwasawa decomposition
$G'(\mathbb{A})=Z'(\mathbb{A})P_0'(\mathbb{A})K'$, and set
$\mathfrak{c}:=\eta ky^{-1}$. Then
\begin{align*}
\mathcal{I}(a,y)\ll \int\int\big|W'(bk)\big|^2\textbf{1}_{\frac{|\det b|}{|\det a|_{\infty}|\det y|_{\mathbb{A}}\|\eta ky^{-1}\|^n}\asymp C_0}\int_{Z'(\mathbb{A})}h^{\dagger}(\eta zky^{-1})d^{\times}zdbdk,
\end{align*}
where $b$ and $k$ range over $N'(\mathbb{A})\backslash P_0'(\mathbb{A})$
and $K'$, respectively. Identifying a central matrix with its scalar entry, we
have
\begin{align*}
\int_{Z'(\mathbb{A})}h^{\dagger}(\eta zky^{-1})d^{\times}z=\int_{\mathbb{R}^{\times}}h_{\infty}^{\dagger}(z_{\infty}\mathfrak{c}_{\infty})d^{\times}z_{\infty}\prod_{p<\infty}\int_{\mathbb{Q}_p^{\times}}h_{p}^{\dagger}(z_{p}\mathfrak{c}_p)d^{\times}z_{p}.
\end{align*}

For $p<\infty$, set $m_p:=\min_j e_p(c_{j,p})$. By definition,
$h_p^{\dagger}(z_p\mathfrak{c}_p)$ is nonzero precisely when
$e_p(z_p)=-m_p$. Thus
\begin{align*}
\int_{\mathbb{Q}_p^{\times}}h_{p}^{\dagger}(z_{p}\mathfrak{c}_p)d^{\times}z_{p}=1
\end{align*}
for all $p<\infty$. At the Archimedean place, scaling by
$\|\mathfrak{c}_\infty\|_\infty^{-1}$ gives
\begin{align*}
\int_{\mathbb{R}^{\times}}h_{\infty}^{\dagger}(z_{\infty}\mathfrak{c}_{\infty})d^{\times}z_{\infty}=\int_{\mathbb{R}^{\times}}h_{\infty}(|z_{\infty}|)d^{\times}z_{\infty}=2\int_{0}^{\infty}\frac{h_{\infty}(t)}{t}dt\asymp 1.
\end{align*}

Therefore,
$\int_{Z'(\mathbb{A})}h^{\dagger}(\eta zky^{-1})d^{\times}z\asymp1$
uniformly in $k$ and $y$. Consequently,
\begin{align*}
\mathcal{I}(a,y)\ll \int_{K'}\frac{J(a,y;k)}{|\det a|_{\infty}|\det y|_{\mathbb{A}}\|\mathfrak{c}\|^n}\int_{N'(\mathbb{R})\backslash P_0'(\mathbb{R})}\big|W_{\infty}'(b_{\infty}k_{\infty})\big|^2|\det b_{\infty}|_{\infty}db_{\infty}dk,
\end{align*}
where  
\begin{align*}
J(a,y;k):=\int_{N'(\mathbb{A}_{\fin})\backslash P_0'(\mathbb{A}_{\fin})}\big|W_{\fin}'(b_{\fin}k_{\fin})\big|^2|\det b_{\fin}|_{\fin}\textbf{1}_{|\det b_{\fin}|_{\fin}\asymp X}db_{\fin},
\end{align*}
with $X:=|\det a|_{\infty}|\det y|_{\mathbb{A}}\|\eta ky^{-1}\|^nC_0|\det b_{\infty}|_{\infty}^{-1}$. 

The support properties of the normalized essential Whittaker vectors imply
that $J(a,y;k)=0$ unless $X\ll 1$. If
$\lambda_{\pi'\times\widetilde{\pi}'}(m)$ denotes the $m$-th Dirichlet
coefficient of $L(s,\pi'\times\widetilde{\pi}')$, the standard dyadic
Rankin--Selberg estimate gives
\begin{align*}
J(a,y;k)\ll \sum_{m\asymp X^{-1}}\frac{\lambda_{\pi'\times\widetilde{\pi}'}(m)}{m}\ll_{\pi'}1.
\end{align*}

Let $\delta_{H\cap B'}$ be the modular character on $H\cap B'$. Then 
\begin{equation}\label{89}
\delta_{B'}^{-1}(a')=|a_n'|^{n}|\det a'|^{-1}\delta_{H\cap B'}^{-1}(a'),\ a'=\diag(a_1',\cdots,a_n').
\end{equation}

Combining \eqref{89} with the Kirillov model and local Rankin--Selberg theory
gives
\begin{align*}
\mathcal{I}(a,y)\ll \frac{L(1,\pi',\Ad)}{|\det a|_{\infty}|\det y|_{\mathbb{A}}}\int_{K'}\frac{1}{\|\mathfrak{c}\|^n}dk\int_{N_H(\mathbb{R})\backslash H(\mathbb{R})}\big|W_{\infty}'(h)\big|^2dh,
\end{align*}
where the implied constant depends at most on $\pi'$. 

Equation \eqref{4.15.} now follows from Lemma \ref{lem4.6}.
\end{proof}

\begin{lemma}\label{norm}
Let $\phi'\in\pi'$ be defined in \textsection\ref{sec4.1.3}. For
$a=\diag(a_1,\ldots,a_{n-1},1)\in A_H(\mathbb{R})$ and
$y\in G'(\mathbb{A})$, one has
\begin{align*}
\int_{\Omega K_{\fin}'}\big|\phi'(ag)\phi'(agy)\big|dg\ll T^{\varepsilon}\langle\phi',\phi'\rangle\delta_{B'}(a)\cdot \min\Bigg\{\frac{|\det a|_{\infty}d(\transp{y}^{-1})}{|a_1|_{\infty}^n},\frac{d(y)}{|\det a|_{\infty}}\Bigg\},
\end{align*}
where the implied constant depends at most on $\varepsilon$, $n$, and $\pi'$,
and $d(\cdot)$ is defined by \eqref{dy}.
\end{lemma}
\begin{proof}
We first prove
\begin{equation}\label{91}
\int_{\Omega K_{\fin}'}\big|\phi'(ag)\phi'(agy)\big|dg
\ll T^{\varepsilon}\langle\phi',\phi'\rangle\delta_{B'}(a)
\frac{d(y)}{|\det a|_{\infty}}.
\end{equation}

Apply Lemma \ref{lem4.7'} to
$h(g)=|\phi'(g)\phi'(gy)|$. It remains to bound
\begin{equation}\label{94}
\int_{P_0'(\mathbb{Q})\backslash \overline{G'}(\mathbb{A})}
|\phi'(g)\phi'(gy)|
\textbf{1}_{C_0|\det a|_{\infty}\leq|\det g|\leq
2C_0|\det a|_{\infty}}\,dg.
\end{equation}
By Cauchy--Schwarz, \eqref{94} is
\begin{equation}\label{4.18}
\ll \sqrt{\mathcal{J}(a,I_n)\mathcal{J}(a,y)}
\end{equation}
where
\begin{align*}
\mathcal{J}(a,y):=\int_{P_0'(\mathbb{Q})\backslash \overline{G'}(\mathbb{A})}\big|\phi'(gy)\big|^2\textbf{1}_{\substack{|\det g|\geq C_0|\det a|_{\infty}\\|\det g|\leq 2C_0|\det a|_{\infty}}}dg.
\end{align*}

Unfolding the Fourier--Whittaker expansion gives
\begin{align*}
\mathcal{J}(a,y)=\int_{N'(\mathbb{A})\backslash \overline{G'}(\mathbb{A})}\big|W'(gy)\big|^2\textbf{1}_{C_0|\det a|_{\infty}\leq |\det g|\leq 2C_0|\det a|_{\infty}}dg.
\end{align*}

For $z\in Z'(\mathbb{A})$ and $k\in K'$, the definition of $h^\dagger$
gives
\begin{align*}
h^{\dagger}(\eta zk)
=h_{\infty}(|\det z_{\infty}|_{\infty}^{1/n})
\textbf{1}_{z_{\fin}\in\widehat{\mathbb{Z}}^{\times}I_n}.
\end{align*}
Consequently,
\begin{align*}
\mathcal{J}(a,y)\ll \int_{N'(\mathbb{A})\backslash G'(\mathbb{A})}
\big|W'(gy)\big|^2
\textbf{1}_{C_0|\det a|_{\infty}\leq |\det g|\leq
2^{n+1}C_0|\det a|_{\infty}}h^{\dagger}(\eta g)dg=\mathcal{I}(a,y).
\end{align*}
Lemma \ref{lem4.8.} therefore yields
\begin{align*}
\mathcal{J}(a,y)
\ll T^{\varepsilon}|\det a|_{\infty}^{-1}d(y)^2
\langle\phi',\phi'\rangle.
\end{align*}
Since $d(I_n)=1$, substitution into \eqref{4.18} proves \eqref{91}.

It remains to obtain the other bound. Let $w_{n-1}$ be the long Weyl element
of $G'$ and set
\begin{align*}
a^\vee:=a_1w_{n-1}\transp{a}^{-1}w_{n-1}^{-1},
\qquad
y^\vee:=w_{n-1}\transp{y}^{-1}w_{n-1}^{-1}.
\end{align*}
The element $a^\vee$ lies in the image of $A_H(\mathbb{R})$ under this
involution, and
\begin{align*}
\delta_{B'}(a^\vee)=\delta_{B'}(a),
\qquad
|\det a^\vee|_{\infty}
=\frac{|a_1|_{\infty}^n}{{|\det a|_{\infty}}},
\qquad
d(y^\vee)=d(\transp{y}^{-1}).
\end{align*}
The inverse-transpose map, followed by conjugation by $w_{n-1}$, carries
$\phi'$ to a vector in $\widetilde{\pi}'$ with the same Petersson norm. It
also carries the preceding Siegel domain to one associated with the
transformed Weyl chamber. Repeating the proof of \eqref{91} on that domain,
with $(a,y)$ replaced by $(a^\vee,y^\vee)$, gives
\begin{align*}
\int_{\Omega K_{\fin}'}|\phi'(ag)\phi'(agy)|\,dg
\ll T^{\varepsilon}\langle\phi',\phi'\rangle\delta_{B'}(a)
\frac{|\det a|_{\infty}d(\transp{y}^{-1})}{|a_1|_{\infty}^n}.
\end{align*}
Taking the minimum of this estimate and \eqref{91} proves the lemma.
\end{proof}

\begin{lemma}\label{lem4.9}
For every $c_0\geq0$, there is a $c_1\geq 0$ such that for each fixed
$c_2\in\mathbb{R}$, we have, for all $a\in A_H(\mathbb{R})$ and
$g\in \Omega K_{\fin}'$, that
\begin{equation}\label{97}
\big|\phi'(ag)\big|\ll \Bigg[\prod_{j=1}^{n-1}\max\{a_j,a_j^{-1}\}\Bigg]^{-c_0}T^{c_1}|\det a|^{c_2}.
\end{equation}
\end{lemma}
\begin{proof}
This is the cuspidal analogue of the crude-growth estimate in
\cite[Lemma 2.23]{Nel21}. Applying the argument there to the cuspidal
realization of $W_\infty'$ reduces \eqref{97} to bounds for suitable Sobolev
seminorms of the underlying test vector. These bounds follow from
\cite[Lemma 21.27]{Nel21}. Crucially, the exponent $c_1$ may be chosen
independently of the fixed parameter $c_2$; only the auxiliary seminorms and
the implied constant may depend on $c_2$. This proves \eqref{97}.
\end{proof}

\section{The Spectral Side}\label{sec5}
Let $\pi=\pi_1\boxplus\cdots\boxplus\pi_m\in\mathcal{A}([G],\omega)$ be
the pure isobaric representation fixed in \textsection\ref{3.1}. Thus $\pi$ is
the Langlands quotient induced from the unitary cuspidal representations
$\pi_j$ on the Levi subgroup of a standard parabolic $Q$ of type
$(n_1,\ldots,n_m)$, where $n_1+\cdots+n_m=n+1$. We assume that $\pi$ has
uniform parameter growth of size $(T;c_\infty,C_\infty)$ as in \eqref{1.6}.
Recall the deformation parameter
$\mu\in i\mathfrak a_Q^*/i\mathfrak a_G^*\simeq(i\mathbb R)^{m-1}$:
\begin{equation}\label{5.1.}
\mu=\begin{cases}
	0,\ & \text{if $m=1$,}\\
	(\mu_1,\cdots,\mu_m),\ & \text{if $m\geq 2$,}
\end{cases}
\end{equation}
with $\mu=it_\mu\mathbf c$ when $m\geq2$, as in
\textsection\ref{specpara}. Define
\begin{equation}\label{mu}
\boldsymbol{L}(\mu,\pi):=\prod_{j=1}^mL(1,\pi_j,\Ad)\prod_{1\leq j_1<j_2\leq m}
\big|L(1+\mu_{j_1}-\mu_{j_2},
\pi_{j_1}\otimes\widetilde{\pi}_{j_2})\big|^2.
\end{equation}
By the nonvanishing theorem for Rankin--Selberg $L$-functions on
$\Re(s)=1$, we have $\boldsymbol{L}(\mu,\pi)\neq0$.
	
For $\lambda=(\lambda_1,\ldots,\lambda_m)\in\mathbb{C}^m$, set
\begin{equation}\label{eq5.3}
\pi_{\lambda}:=\pi_1|\cdot|^{\lambda_1}\boxplus \cdots\boxplus \pi_m|\cdot|^{\lambda_m}.
\end{equation}
This agrees with the notation $\pi_\mu$ introduced in
\textsection\ref{specpara}.

Recall from \textsection\ref{ampsp} that
\begin{align*}
\mathcal{J}_{\Spec}^{\heartsuit}(\boldsymbol{\alpha},\boldsymbol{\ell})=\sum_{p_1\neq p_2\in\mathcal{L}}\alpha_{p_1}\overline{\alpha_{p_2}}J_{\Spec}^{\Reg,\heartsuit}(f_{p_1,p_2}^{\mathrm{off}},\mathbf{0})+\sum_{p_0\in\mathcal{L}}\sum_{i=0}^{\ell_{p_0}}c_{p_0,i}|\alpha_{p_0}|^2J_{\Spec}^{\Reg,\heartsuit}(f_{p_0,i}^{\diag},\mathbf{0}).
\end{align*}

The main result of this section is the following spectral lower bound.
\begin{restatable}[]{thm}{thmf} \label{prop11.2}
Let notation be as above. Suppose that
$|L(1/2,\pi\times\pi')|\geq2$. Then
\begin{align*}
\mathcal{J}_{\Spec}^{\heartsuit}(\boldsymbol{\alpha},\boldsymbol{\ell})\gg  T^{-\frac{n^2}{2}-\varepsilon}\cdot  \big|L(1/2,\pi\times\pi')\big|^2\cdot \bigg|\sum_{p\in \mathcal{L}}\alpha_p\lambda_{\pi_{\mu}}(p^{\ell_p})\bigg|^2\cdot \prod_{j=1}^m\frac{1}{L(1,\pi_j,\Ad)},
\end{align*}
where the implied constant depends on $\varepsilon$, $c_\infty$, $C_\infty$,
and $\pi_\infty'$.
\end{restatable}

\subsection{The Archimedean Period}

\begin{lemma}\label{prop4.11}
	There exists a vector $W_{\mu,\infty}$ in the Whittaker model of
	$\pi_{\mu,\infty}$ such that
	\begin{align*}
	\int_{N'(\mathbb{R})\backslash G'(\mathbb{R})}\Bigg|W_{\mu,\infty}\left(\begin{pmatrix}
			g'\\
			&1
		\end{pmatrix}\right)\Bigg|^2dg'\leq 1,
	\end{align*}
	and 
	\begin{align*}
		T^{-\frac{n^2}{4}-\varepsilon}\ll \Bigg|\int_{N'(\mathbb{R})\backslash G'(\mathbb{R})}W_{\mu,\infty}\left(\begin{pmatrix}
			g'\\
			&1
		\end{pmatrix}\right)W_{\infty}'(g')dg'\Bigg|\ll  T^{-\frac{n^2}{4}+\varepsilon},
	\end{align*}
	where the implied constants depend on $\pi_\infty'$ and $\varepsilon$.
\end{lemma}
\begin{proof}
Apply the rescaled test-vector construction of
\cite[Lemmas 8.5--8.6 and \S 11.4]{Nel21} to $\pi_{\mu,\infty}$ and normalize
the resulting vector to have the stated $L^2$-norm. Lemma
\ref{lem4.5}, together with the rescaling identity
\cite[Lemma 10.3]{Nel21}, gives the lower bound for the local zeta integral;
the corresponding $L^2$ estimate \cite[Lemma 10.20]{Nel21} gives the upper
bound. The normalization in \cite[\S 11.4]{Nel21} contributes
$T^{-\dim G'/4}=T^{-n^2/4}$. Since the inducing exponents of $\chi$ are fixed
and $\mu$ is purely imaginary with $|\mu_j|\leq1$, all remaining losses are
absorbed into $T^\varepsilon$.
\end{proof}

\subsection{The Auxiliary Local Period}
Write $\Psi_p(s,W_p,W_p')$ for the local Rankin--Selberg integral associated
with \eqref{whittaker}.

Let $p=p_{\circ}$. For
$\lambda\in i\mathfrak a_Q^*/i\mathfrak a_G^*$, let $W_{\lambda,p}$ be the
normalized spherical Whittaker function of $\pi_{\lambda,p}$, and let $W_p'$
be that of $\pi_p'$.
\begin{lemma}\label{lem5.aux}
Uniformly for $\lambda$ in a sufficiently small fixed neighborhood of $0$ in
$i\mathfrak a_Q^*/i\mathfrak a_G^*$,
\begin{align*}
\big|\Psi_p(0,\pi_{\lambda,p}(\tilde f_p)W_{\lambda,p},W_p')\big|
\gg
\big|L_p(1/2,\pi_{\lambda,p}\times\pi_p')\big|
\gg 1.
\end{align*}
The implied constants may depend on $n$, $p_{\circ}$, and $\pi_{p_{\circ}}'$.
\end{lemma}
\begin{proof}
The argument of \cite[Lemma 3.6]{Yan23b}, applied in the lower-right
$2\times2$ block, gives the first inequality once $p_{\circ}$ is chosen
sufficiently large; spherical invariance in the remaining coordinates leaves
the calculation unchanged. The second inequality follows from the standard
bounds for unramified Satake parameters. Since
$f_p=\tilde f_p*\tilde f_p^*$, the corresponding contribution to the spectral
side is the square of this local period and is therefore nonnegative.
\end{proof}

\subsection{Ramified Local Periods}
Let $p\mid M''$, and recall the function $\tilde f_p$ from
\textsection\ref{11.1.4}:
\begin{equation*}
\tilde{f}_p(g):=v_p\sum_{\boldsymbol{\alpha}}\omega_p'(\alpha_n)\int_{Z(\mathbb{Q}_p)}\textbf{1}_{u_{\boldsymbol{\alpha}}K_p}(zg)\omega_p(z)d^{\times}z,\tag{\ref{3.5}}
\end{equation*}
where $\boldsymbol{\alpha}=(\alpha_1,\cdots,\alpha_{n-1},\alpha_n),$ $\alpha_j\in  \mathbb{Z}/p^{m'}\mathbb{Z},$ $1\leq j<n,$ $\alpha_n\in (\mathbb{Z}/p^{m''}\mathbb{Z})^{\times},$ and   $u_{\boldsymbol{\alpha}}=\begin{pmatrix}
I_n&\mathfrak{b}_{\boldsymbol{\alpha}}\\
&1
\end{pmatrix},$ with $\mathfrak{b}_{\boldsymbol{\alpha}}=\transp{(\alpha_1p^{-m'},\cdots, \alpha_{n-1}p^{-m'},\alpha_n p^{-m''})},$ and 
\begin{align*}
v_p=\Vol(K_p'(M'))^{-1}\Vol(K_p(M))^{-1}p^{-(n-1)m'}G(\omega_p',\psi_p)^{-1}.
\end{align*} 

\begin{lemma}\label{lem5.3}
Let $\lambda\in i\mathfrak a_Q^*/i\mathfrak a_G^*$, let $W_{\lambda,p}$ be
an essential Whittaker vector of $\pi_{\lambda,p}$, and let $W_p'$ be the
normalized new vector of $\pi_p'$. Then, as meromorphic functions of $s$,
\begin{align*}
\Psi_p(s,\pi_{\lambda,p}(\tilde f_p)W_{\lambda,p},W_p')
=W_{\lambda,p}(I_{n+1})W_p'(I_n)
L_p(s+1/2,\pi_{\lambda,p}\times\pi_p').
\end{align*}
\end{lemma}
\begin{proof}
Unfold \eqref{3.5} and write
$x=ak$, with $a=\diag(p^{r_1},\ldots,p^{r_n})$ and $k\in K_p'$.
Orthogonality in $\alpha_1,\ldots,\alpha_{n-1}$ restricts $k$ to
$K_p'(p^{m'-r_n})$. The Gauss sum in $\alpha_n$, together with the newvector
transformation law for $W_p'$, then forces $r_n\leq0$, whereas the support of
the essential vector $W_{\lambda,p}$ forces
$r_1\geq\cdots\geq r_n\geq0$. Thus $r_n=0$, and the normalizing constants
cancel by the definition of $v_p$. The resulting sum is the usual toral
expansion of the local Rankin--Selberg integral; see
\cite[Theorem 4.1]{Miy14}. This proves the identity.
\end{proof}

\subsection{Lower Bound of the Amplified Spectral Side}\label{sec5.2}
\subsubsection{Spectral decomposition}\label{sec5.2.1}
Let $R=M_RN_R$ be a standard parabolic subgroup of type
$(n_1',\ldots,n_l')$. Write
\begin{align*}
\mathfrak a_R=\Hom_{\mathbb Z}(X(M_R)_{\mathbb Q},\mathbb R),
\qquad
\mathfrak a_R^*=X(M_R)_{\mathbb Q}\otimes_{\mathbb Z}\mathbb R,
\end{align*}
and, for $\mathbf m=\diag(m_1,\ldots,m_l)\in M_R(\mathbb A)$, set
\begin{align*}
H_R(\mathbf m)
=\left(\frac{\log|\det m_1|}{n_1'},\ldots,
\frac{\log|\det m_l|}{n_l'}\right)\in\mathfrak a_R.
\end{align*}
Arthur's spectral decomposition \cite[pp. 256, 263]{Art79} gives
\begin{equation}\label{49}
L^2(R):=L^2\left(Z_G(\mathbb A)N_R(\mathbb A)M_R(F)\backslash G(\mathbb A)\right)
=\bigoplus_{\mathfrak X}L^2(R)_{\mathfrak X},
\end{equation}
where $\mathfrak X$ ranges over cusp data. For each $\mathfrak X$, choose an
orthonormal basis $\mathfrak B_{R,\mathfrak X}$ consisting of $K$-finite pure
tensors.

For $\lambda\in\mathfrak a_{R,\mathbb{C}}^*$, let
$\mathcal{I}_R(\cdot,\lambda)$ denote the normalized induced representation
attached to $\mathfrak X$, and let
\begin{align*}
E(x,\phi,\lambda)
:=\sum_{\delta\in R(F)\backslash G(F)}
\phi(\delta x)e^{\langle\lambda+\rho_R,H_R(\delta x)\rangle}
\end{align*}
be the associated Eisenstein series. Its Whittaker function is
\begin{align*}
W_{\phi,\lambda}(x)
:=\int_{[N]}E(ux,\phi,\lambda)\overline{\theta(u)}\,du.
\end{align*}

\begin{prop}\label{thm6}
Let notation be as in \textsection\ref{11.2}.  
\begin{enumerate}
\item Suppose $m=1,$ i.e., $\pi$ is cuspidal. Then 
\begin{equation}\label{5.3'}
\mathcal{J}_{\Spec}^{\heartsuit}(\boldsymbol{\alpha},\boldsymbol{\ell})\gg_{\varepsilon,\pi'} T^{-\frac{n^2}{2}-\varepsilon}\cdot  \frac{\big|L(1/2,\pi\times\pi')\big|^2}{L(1,\pi,\Ad)}\cdot \bigg|\sum_{p\in \mathcal{L}}\alpha_p\lambda_{\pi}(p^{\ell_p})\bigg|^2,
\end{equation}
where the implied constant depends on $\varepsilon,$ $\pi',$ and $\pi_p,$ $p\mid M'$.
\item Suppose $m\geq2$, and let $\mu=it_\mu\mathbf c$ be as in
\textsection\ref{specpara}. Then
\begin{equation}\label{5.3''}
\mathcal{J}_{\Spec}^{\heartsuit}(\boldsymbol{\alpha},\boldsymbol{\ell})\gg_{\varepsilon,\pi'} T^{-\frac{n^2}{2}-\varepsilon}\cdot \frac{\big|L(1/2,\pi_{\mu}\times\pi')\big|^2}{\boldsymbol{L}(\mu,\pi)}\cdot \bigg|\sum_{p\in \mathcal{L}}\alpha_p\lambda_{\pi_{\mu}}(p^{\ell_p})\bigg|^2,
\end{equation}
where the implied constant depends on $\varepsilon$ and $\pi_p,$ $p\mid M'$. Here $\boldsymbol{L}(\mu,\pi)$ is defined by \eqref{mu}. 
\end{enumerate}
\end{prop}

\begin{proof}
For $f\in\{f_{p_0,i}^{\diag},f_{p_1,p_2}^{\mathrm{off}}\}$, let
$\tilde f$ be as in \textsection\ref{3.6.2}. The spectral expansion gives
\begin{equation}\label{eqspec}
J_{\Spec}^{\Reg,\heartsuit}(f,\mathbf{0})=\sum_{\mathfrak{X}}\sum_{R}\frac{1}{n_R}\left(\frac{1}{2\pi i}\right)^{\dim A_R/A_G}\int_{i\mathfrak{a}_R^*/i\mathfrak{a}_G^*}\Psi^{\dagger}(\lambda)d\lambda,
\end{equation}
where
\begin{align*}
\Psi^{\dagger}(\lambda):=\sum_{\phi\in\mathfrak{B}_{R,\mathfrak{X}}}\Psi(0,
\mathcal{I}_R(\lambda,f)W_{\phi,\lambda},W_{\phi'}')\overline{\Psi(0,
W_{\phi,\lambda},W_{\phi'}')}.
\end{align*}
Since the amplification primes are coprime to $M$, the Hecke action factors as
\begin{align*}
\mathcal{I}_R(\lambda,f)W_{\phi,\lambda}=\boldsymbol{c}_{\pi_{\lambda}}(f)\mathcal{I}_R(\lambda,\tilde{f}*\tilde{f}^*)W_{\phi,\lambda},
\end{align*}
where $\boldsymbol c_{\pi_\lambda}(f)$ is the corresponding eigenvalue. Thus
\begin{equation}\label{eq5.10}
\Psi^{\dagger}(\lambda)=\sum_{\phi\in\mathfrak{B}_{R,\mathfrak{X}}}\boldsymbol{c}_{\pi_{\lambda}}(f)\big|\Psi(0,
\mathcal{I}_R(\lambda,\tilde{f})W_{\phi,\lambda},W_{\phi'}')\big|^2.
\end{equation}

The Hecke relation \eqref{3.66} gives
\begin{equation}\label{5.9}
\sum_{p_1\neq p_2\in\mathcal{L}}\alpha_{p_1}\overline{\alpha_{p_2}}\boldsymbol{c}_{\pi_{\lambda}}^{p_1,p_2}+\sum_{p_0\in\mathcal{L}}\sum_{i=0}^{\ell_{p_0}}c_{p_0,i}|\alpha_{p_0}|^2\boldsymbol{c}_{\pi_{\lambda}}^{i;p_0}=\bigg|\sum_{p\in \mathcal{L}}\alpha_p\lambda_{\pi_{\lambda}}(p^{\ell_p})\bigg|^2
\end{equation}
where $\boldsymbol{c}_{\pi_{\lambda}}^{p_1,p_2}=\boldsymbol{c}_{\pi_{\lambda}}(f_{p_1,p_2}^{\mathrm{off}})$ and $\boldsymbol{c}_{\pi_{\lambda}}^{i;p_0}:=\boldsymbol{c}_{\pi_{\lambda}}(f_{p_0,i}^{\diag})$.

We now form the amplified combination before discarding any spectral terms.
By \eqref{5.9}, every resulting spectral weight is nonnegative. We may
therefore retain only the cusp datum
$\mathfrak X=\{(M_Q,\pi_1\otimes\cdots\otimes\pi_m)\}$ and $R=Q$, obtaining
\begin{equation}\label{5.3.}	
\mathcal J_{\Spec}^{\heartsuit}(\boldsymbol\alpha,\boldsymbol\ell)
\gg \int_{i\mathfrak{a}_{Q}^*/i\mathfrak{a}_G^*}
\bigg|\sum_{p\in\mathcal L}\alpha_p
\lambda_{\pi_\lambda}(p^{\ell_p})\bigg|^2
\sum_{\phi\in\mathfrak B_{Q,\mathfrak X}}
\big|\Psi(0,\mathcal{I}_Q(\lambda,\tilde f)W_{\phi,\lambda},
W_{\phi'}')\big|^2d\lambda.
\end{equation}

Suppose first that $m\geq2$. The mean-value inequality for holomorphic
functions, applied on a $T^{-o(1)}$-neighborhood of $\mu$ as in
\cite[Theorem 1.2.1 and its corollary]{Ram95}, bounds \eqref{5.3.} from below
by $T^{-\varepsilon}$ times its integrand at $\lambda=\mu$. Rankin--Selberg
factorization \cite{JPSS81,JPSS83}, Lemmas \ref{lem5.aux} and \ref{lem5.3},
and the normalization $W_p'(I_n)=1$ give
\begin{align*}
\big|\Psi(0,\mathcal{I}_{Q}(\mu,\tilde{f})W_{\phi,\mu},W_{\phi'}')\big|^2
\gg
\big|\Psi_{\infty}(0,\mathcal{I}_{Q}(\mu,\tilde{f})W_{\phi,\mu},W_{\phi'}')\big|^2
\frac{|L(1/2,\pi_{\mu}\times\pi')|^2}
{\boldsymbol{L}(\mu,\pi)},
\end{align*}
where $\boldsymbol L(\mu,\pi)$ is defined by \eqref{mu}. Consequently,
\begin{align*}
\mathcal{J}_{\Spec}^{\heartsuit}(\boldsymbol{\alpha},\boldsymbol{\ell})
\gg T^{-\varepsilon}\mathcal{Q}_{\boldsymbol{\alpha},\boldsymbol{\ell}}(\pi,\pi';\mu)
\sum_{\phi\in\mathfrak{B}_{Q,\mathfrak{X}}}
\big|\Psi_{\infty}(0,\mathcal{I}_{Q}(\mu,\tilde{f})W_{\phi,\mu},W_{\phi'}')\big|^2,
\end{align*}
where
\begin{align*}
	\mathcal{Q}_{\boldsymbol{\alpha},\boldsymbol{\ell}}(\pi,\pi';\mu):=\frac{\big|L(1/2,\pi_{\mu}\times\pi')\big|^2}{\boldsymbol{L}(\mu,\pi)}\cdot \bigg|\sum_{p\in \mathcal{L}}\alpha_p\lambda_{\pi_{\mu}}(p^{\ell_p})\bigg|^2.
	\end{align*}
Lemma \ref{prop4.11} implies that
\begin{align*}
\sum_{\phi\in\mathfrak{B}_{Q,\mathfrak{X}}}\big|\Psi_{\infty}(0,\mathcal{I}_{Q}(\mu,\tilde{f})W_{\phi,\mu},W_{\phi'}')\big|^2\gg T^{-\frac{n^2}{2}-\varepsilon}.
\end{align*}
This proves \eqref{5.3''}.

If $m=1$, then $Q=G$ and the same argument has no continuous parameter.
Taking $\mu=0$ and applying Lemma \ref{prop4.11} directly proves
\eqref{5.3'}.
\end{proof}

\subsection{Proof of the Spectral Lower Bound}
For the noncuspidal case, we use the following zero-avoidance lemma, which
may be of independent interest. 
\begin{lemma}\label{lem11.1}
Let $\sigma$ be an automorphic representation of
	$\mathrm{GL}(r,\mathbb A)$, and let $t_0\in\mathbb R$. Suppose that
	$|L(1/2+it_0,\sigma\times\pi')|\geq2$, and let
	$\textbf{C}:=C(\sigma\times\pi'\otimes|\cdot|^{it_0})\geq10^4$.
	For every integer $l\geq1$, there exists
\begin{align*}
d_l'\in[2^{-l}e^{-3\sqrt{\log\textbf C}},
	2^{1-l}e^{-3\sqrt{\log\textbf C}}]
\end{align*}	
such that, whenever $|t-t_0|=d_l'$, one has
		\begin{equation}\label{11.1}
		|L(1/2+it_0,\sigma\times\pi')|\ll \exp(\log^{3/4}\textbf{C})\cdot |L(1/2+it,\sigma\times\pi')|,
	\end{equation}
	and 
	\begin{equation}\label{5.8}
	|L(1/2+it,\sigma\times\pi')|\ll \exp(\log^{3/4}\textbf{C})\cdot |L(1/2+it_0,\sigma\times\pi')|,	
	\end{equation}
		where the implied constants depend only on $l$.
\end{lemma}
\begin{proof}
Denote by $\mathcal{D}:=\big\{\rho\in\mathbb{C}:\ L(\rho,\sigma\times\pi')=0,\ |\rho-1/2-it_0|\leq \exp(-2\sqrt{\log \textbf{C}})\big\}$. By Jensen's formula and the convexity bound for $L(s,\sigma\times\pi'),$ we have 
 	\begin{align*}
 		\#\mathcal{D}\cdot \log (\exp(\sqrt{\log \textbf{C}}))\leq \max_{|s-1/2-it_0|\leq \exp(-\sqrt{\log \textbf{C}})}\log \frac{|L(s,\sigma\times\pi')|}{|L(1/2+it_0,\sigma\times\pi')|},
 	\end{align*} 	
 	which is $\ll \log \textbf{C},$ and the implied constant is absolute. Consequently, 
	\begin{align*}
	\#\mathcal{D}\ll  \sqrt{\log \textbf{C}},
	\end{align*}
 	where the implied constant is absolute. 
	By the pigeonhole principle, there exists $d_l\in [2^{-l}\exp(-2\sqrt{\log \textbf{C}}),2^{1-l}\exp(-2\sqrt{\log \textbf{C}})]$ such that 
 	\begin{equation}\label{11.3.}
 	|s-\rho|\geq \frac{\exp(-2\sqrt{\log \textbf{C}})}{3\cdot 2^l\cdot (1+\#\mathcal{D})},
 	\end{equation}
 	for all $|s-1/2-it_0|=d_l$ and for all $\rho\in\mathcal{D}$. 
 	
 	Likewise, there exists some $d_l'\in [2^{-l}\exp(-3\sqrt{\log \textbf{C}}),2^{1-l}\exp(-3\sqrt{\log \textbf{C}})]$ such that for all $|s-1/2-it_0|=d_l'$ and for all $\rho\in\mathcal{D},$ we have 
\begin{equation}\label{5.11}
 	|s-\rho|\geq \frac{\exp(-3\sqrt{\log \textbf{C}})}{3\cdot 2^{l}\cdot (1+\#\mathcal{D})},
\end{equation}
 	
 	 Let $\mathcal{B}_l:=\{s:\ |s-1/2-it_0|\leq d_l\}$ and $\mathcal{B}_l':=\{s:\ |s-1/2-it_0|\leq d_l'\}$.

 	Set $
 		\mathbf{P}(s):=\prod_{\rho\in \mathcal{D}}(s-\rho)$. 	Then $\mathbf{P}$ is a polynomial of degree $\deg\mathbf{P}=\#\mathcal{D}$. Define
	\begin{align*}
	\mathbb{L}(s):=\log\frac{L(s,\sigma\times\pi')}{\mathbf{P}(s)},\ \ s\in\mathcal{B}_l.
	\end{align*}
 	Then $\mathbb{L}(s)$ is holomorphic in $\mathcal{B}_l$. In particular, it is holomorphic in $\mathcal{B}_l'$. Therefore, by Borel-Carath\'{e}odory theorem, we have, for all $s\in \mathcal{B}_l',$ that 
\begin{align*}
\big|\mathbb{L}(1/2+it_0)-\mathbb{L}(s)\big|
\leq \frac{10\exp(-3\sqrt{\log \textbf{C}})\cdot \big[\max_{s'\in \mathcal{B}_l}\Re\mathbb{L}(s')-
		 \Re\mathbb{L}(1/2+it_0)\big]}
{\exp(-2\sqrt{\log \textbf{C}})-\exp(-3\sqrt{\log \textbf{C}})}.
\end{align*}
 	
 	Let $s_0\in\partial \mathcal{B}_l$ be such that $\Re(\mathbb{L}(s_0))=\max_{s'\in \mathcal{B}_l}\Re(\mathbb{L}(s')),$ where $\partial\mathcal{B}_l$ is the boundary of $\mathcal{B}_l$.  Then 
\begin{align*}
\max_{s'\in \mathcal{B}_l}\Re\mathbb{L}(s')-
	 	\Re\mathbb{L}(1/2+it_0)
=\log\Big|\frac{L(s_0,\sigma\times\pi')}
{L(1/2+it_0,\sigma\times\pi')}\Big|+\log\Big|\frac{\mathbf{P}(1/2+it_0)}{\mathbf{P}(s_0)}\Big|.
\end{align*}
 	
 	Using the trivial bound $|\mathbf{P}(1/2+it_0)|\leq 1$ and \eqref{11.3.}, 
 	\begin{align*}
	 	\Big|\frac{\mathbf{P}(1/2+it_0)}{\mathbf{P}(s_0)}\Big|\ll \frac{1}{|\mathbf{P}(s_0)|}\ll \Big[3\cdot 2^l\cdot (1+\#\mathcal{D})\exp(2\sqrt{\log \textbf{C}})\Big]^{\deg \mathbf{P}}.
 	\end{align*}
 	Hence, we have, by $\#\mathcal{D}\ll \sqrt{\log \textbf{C}},$ that
	\begin{align*}
	\log \Big|\frac{\mathbf{P}(1/2+it_0)}{\mathbf{P}(s_0)}\Big|\ll \deg \mathbf{P}\cdot 
	\log \Big[3\cdot 2^l\cdot (1+\#\mathcal{D})\exp(2\sqrt{\log \textbf{C}})\Big]\ll \log \textbf{C},
	\end{align*}
	where the implied constant depends on $l$. Moreover, by the convexity bound, 
\begin{equation}\label{5.13}
 	\log\Big|\frac{L(s_0,\sigma\times\pi')}{L(1/2+it_0,\sigma\times\pi')}\Big|\ll \log |L(s_0,\sigma\times\pi')|\ll  \log \textbf{C}.
\end{equation}

Therefore, uniformly for $s\in \mathcal{B}_l',$ we have
 \begin{equation}\label{5.10}
 \big|\mathbb{L}(1/2+it_0)-\mathbb{L}(s)\big|\ll \exp(-\sqrt{\log \textbf{C}})\cdot \log \textbf{C}\ll 1.
 \end{equation}
 	
 Take $s=1/2+it\in \partial\mathcal{B}_l',$ then \eqref{5.10} leads to 
 		\begin{align*}
 L(1/2+it_0,\sigma\times\pi')\ll  \big|L(1/2+it,\sigma\times\pi')\big|\cdot \Big|\frac{\mathbf{P}(1/2+it_0)}{\mathbf{P}(1/2+it)}\Big|,
\end{align*}
where the implied constant is independent of $t$. Moreover, by \eqref{5.11}, $\mathbf{P}(1/2+it)\neq 0$.
 	
 	For $\rho\in\mathcal{D}$ and $s\in \partial \mathcal{B}_l',$ by \eqref{5.11} there exists some  constant $C$ (which depends at most on $l$ and the analytic conductors of $\sigma$ and $\pi'$) such that for all $\rho\in\mathcal{D},$ 
\begin{equation}\label{5.12}
 	\frac{|1/2+it_0-\rho|}{|s-\rho|}\leq 
 	\begin{cases}
 		1,&\ \text{if $|s-1/2-it_0|\leq d_l'/2,$}\\
 		C,&\ \text{if $d_l'/2\leq |s-1/2-it_0|\leq 3d_l'/4,$}\\
 		C\cdot (1+\#\mathcal{D}),&\ \text{if $3d_l'/4\leq |s-1/2-it_0|\leq 5d_l'/4,$}\\
 		C,&\ \text{if $5d_l'/4\leq |s-1/2-it_0|\leq d_l$.}
 		\end{cases}
\end{equation}
Together with 
$\deg\mathbf{P}=\#\mathcal{D}\leq C_2\sqrt{\log \textbf{C}}$ 
for some absolute constant $C_2,$ we have
\begin{equation}\label{5.14}
	\Big|\frac{\mathbf{P}(1/2+it_0)}{\mathbf{P}(1/2+it)}\Big|\leq (C\sqrt{\log \textbf{C}})^{C_2\sqrt{\log \textbf{C}}}\ll \exp(\log^{3/4}\textbf{C}),
\end{equation}
where the implied constant depends on $C,$ and $C_2$. So \eqref{11.1} follows. 

The treatment of \eqref{5.8} is  different since $|1/2+it_0-\rho|$ may be quite tiny. So we cannot bound ${|s-\rho|}/{|1/2+it_0-\rho|}$ as \eqref{5.12}. By Borel-Carath\'{e}odory theorem, 
\begin{align*}
\big|\mathbb{L}(1/2+it_0)-\mathbb{L}(s)\big|\leq 20\exp(-\sqrt{\log \textbf{C}})\log\Big|\frac{L(s_0,\sigma\times\pi')}{L(1/2+it_0,\sigma\times\pi')}\cdot\frac{\mathbf{P}(1/2+it_0)}{\mathbf{P}(s_0)}\Big|,
\end{align*}
where $s_0$ is defined as before. Employing \eqref{5.13} and the assumption that $\textbf{C}\geq 10^4$, 
\begin{align*}
\log \Big|\frac{L(1/2+it,\sigma\times\pi')}{L(1/2+it_0,\sigma\times\pi')}\cdot \frac{\mathbf{P}(1/2+it_0)}{\mathbf{P}(1/2+it)}\Big|\leq 20+\log\Big|\frac{\mathbf{P}(1/2+it_0)}{\mathbf{P}(s_0)}\Big|,
\end{align*}
where we use the fact that $\exp(-\sqrt{\log \textbf{C}})\log \textbf{C}\leq 1$ and $20\exp(-\sqrt{\log \textbf{C}})\leq 1$. So
\begin{equation}\label{5.15}
\log \Big|\frac{L(1/2+it,\sigma\times\pi')}{L(1/2+it_0,\sigma\times\pi')}\Big|\leq 20+\log\Big|\frac{\mathbf{P}(1/2+it)}{\mathbf{P}(s_0)}\Big|.
\end{equation}

By \eqref{11.3.}, a similar analysis of \eqref{5.12} works for $|1/2+it-\rho|/|s_0-\rho|,$ $\rho\in\mathcal{D}$. As a consequence, 
we have an analogue of \eqref{5.14}:  
\begin{equation}\label{5.16}
	\Big|\frac{\mathbf{P}(1/2+it)}{\mathbf{P}(s_0)}\Big|\leq (C'\sqrt{\log \textbf{C}})^{C_2\sqrt{\log \textbf{C}}}\ll \exp(\log^{3/4}\textbf{C}),
\end{equation}
where $C'$ is some absolute constant, and thus the implied constant is absolute. Therefore, \eqref{5.8} follows from \eqref{5.15} and \eqref{5.16}.
\end{proof}
 
\begin{proof}[Proof of Theorem \ref{prop11.2}]
If $m=1$, the result is \eqref{5.3'}. Suppose that $m\geq2$, and let $\mathbf c$ be as in
\textsection\ref{specpara}. Apply the zero-avoidance argument in the proof of Lemma \ref{lem11.1}
to the holomorphic function
\begin{align*}
z\longmapsto
L(1/2,\pi_{iz\mathbf c}\times\pi')
=
\prod_{j=1}^m
L(1/2+ic_jz,\pi_j\times\pi'),\quad |z|<1/10.
\end{align*}
Indeed, the required growth estimate follows by applying the convexity bound
to each factor. Since $|L(1/2,\pi\times\pi')|\geq2$, there exists $t_\mu=\exp\bigl(-O(\sqrt{\log(MT)})\bigr)$ such that, for $\mu=it_\mu\mathbf c$,
\begin{align*}
|L(1/2,\pi\times\pi')|
\ll
(MT)^\varepsilon
|L(1/2,\pi_\mu\times\pi')|.
\end{align*}

Proposition \ref{thm6} yields
\begin{equation}\label{11.4.}
\mathcal{J}_{\Spec}^{\heartsuit}
(\boldsymbol{\alpha},\boldsymbol{\ell})
\gg_{\varepsilon,\pi'}
T^{-\frac{n^2}{2}-\varepsilon}
\frac{|L(1/2,\pi_\mu\times\pi')|^2}
{\boldsymbol{L}(\mu,\pi)}
\bigg|
\sum_{p\in\mathcal{L}}
\alpha_p\lambda_{\pi_\mu}(p^{\ell_p})\bigg|^2.
\end{equation}
Since the coordinates of $\mu$ are distinct and
$|\mu_{j_1}-\mu_{j_2}|^{-1}=(MT)^{o(1)}$, the standard bounds on
$\Re(s)=1$ give
\begin{align*}
\boldsymbol{L}(\mu,\pi)
\ll
\prod_{1\leq j_1<j_2\leq m}
|\mu_{j_1}-\mu_{j_2}|^{-2}
\prod_{j=1}^m L(1,\pi_j,\Ad)
\ll
(MT)^\varepsilon
\prod_{j=1}^m L(1,\pi_j,\Ad).
\end{align*}
Substituting these estimates into \eqref{11.4.} proves the theorem.
\end{proof}

\section{Classification of the Amplification Support}\label{sec5.1}

At each amplification prime $p\mid\nu(f)$, the local test function is supported
on a Cartan double coset. This section describes that support in the mirabolic
coordinates appearing in the orbital integrals of \eqref{63}. The resulting
classification will be used throughout the analysis of the amplified geometric
side.

\subsection{Notation}
Let $W'$ be the Weyl group of $G'=\mathrm{GL}(n)$, let $B'=T'N'$ be the
standard Borel subgroup, and let $\delta_{B'}$ denote its modular character.
Recall from \textsection\ref{sec3.1.6} that $\mathcal L$ is a subset of
\begin{align*}
\{p\text{ prime}:L<p\leq 2L,\ p\nmid p_\circ MM'\},
\end{align*}
and that $1\leq\ell_p\leq n+1$ for each $p\in\mathcal L$. Fix
$p\in\mathcal L$ and suppress $p$ from the local notation when no confusion
can arise.

\subsection{Off-Diagonal Amplification Support}\label{sec6.2}

Let $K_p'=G'(\mathbb Z_p)$. By $p$-adic Gaussian elimination, for each
$y\in G'(\mathbb Q_p)$ one can choose $w\in W'$ and write
\begin{equation}\label{dec}
	wy=tuk,\qquad
	t\in \diag(p^{\mathbb Z},\ldots,p^{\mathbb Z}),\quad
	u\in N'(\mathbb Q_p),\quad k\in K_p',
\end{equation}
with $tut^{-1}\in K_p'$. In particular, \eqref{dec} refines the Cartan
decomposition:
\begin{align*}
y=(w^{-1}tut^{-1})tk\in K_p'tK_p'.
\end{align*}
When $n=2$, all blocks and vectors of size $n-2$ below are omitted.

\begin{lemma}\label{lem7}
Let $f=f_{p_1,p_2}^{\mathrm{off}}$, set $p=p_2$, and write
$\ell=\ell_{p_2}$. Let $y\in G'(\mathbb Q_p)$ and
$\mathfrak u\in M_{n,1}(\mathbb Q_p)\simeq\mathbb Q_p^n$. Suppose that the
minimum valuations of the entries of $y$ and $\mathfrak u$ are $-e$ and
$-e'$, respectively, with $e,e'\geq0$, and write $wy=tuk$ as in \eqref{dec}.
Write $w\mathfrak u=\transp{(\mathfrak u',\alpha)}$, where
$\mathfrak u'\in\mathbb Q_p^{n-1}$ and $\alpha\in\mathbb Q_p$.
If
$f_p\left(\begin{pmatrix}
			y&\mathfrak{u}\\
			&1
\end{pmatrix}\right)\neq 0$, then $\delta_{B'}^{-1}(t)\leq p^{(n-1)\ell}$, and
the following alternatives hold.
\begin{enumerate}
\item[(A)] If $e\geq e'$, then $e=\ell$ and
\begin{align*}
wy=\begin{pmatrix}
I_{n-1}&\\
&p^{-\ell}
\end{pmatrix}\begin{pmatrix}
I_{n-1}&\mathfrak{u}''\\
&1
\end{pmatrix}k,\quad \mathfrak{u}''\in p^{-\ell}\mathbb{Z}_p^{n-1},\ k\in K_p',
\end{align*}
and $\mathfrak{u}'-\alpha p^\ell\mathfrak{u}''\in \mathbb{Z}_p^{n-1}$, with
$\alpha\in p^{-\ell}\mathbb{Z}_p$.

\item[(B)] If $e<e'$, then $2e'-e=\ell$,
$\mathfrak{u}'\in p^{-e'}\mathbb{Z}_p^{n-1}$, and
$\alpha\in p^{-e'}\mathbb{Z}_p^{\times}$. Moreover,
\begin{align*}
wy=\begin{pmatrix}
p^{\frac{\ell-e}{2}}I_{n-1}&\\
				&p^{-e}
\end{pmatrix}\begin{pmatrix}
I_{n-1}&\mathfrak{u}''\\
&1
\end{pmatrix}k,\quad \mathfrak{u}''-p^{-e'}\alpha^{-1}\mathfrak{u}'\in \mathbb{Z}_p^{n-1},\ k\in K_p'.
		\end{align*}
			In particular, $\ell>e'>e\geq 0$.
	\end{enumerate}
\end{lemma}
\begin{proof}
Write $t=\diag(p^{r_1-e},\ldots,p^{r_{n-1}-e},p^{-e})$. By the
choice of \eqref{dec},
\begin{align*}
r_1\geq\cdots\geq r_{n-1}\geq0.
\end{align*}
Write $u=\begin{pmatrix}
			u'&\mathfrak{u}''\\
			&1
		\end{pmatrix}$, with $\mathfrak{u}''\in\mathbb{Q}_p^{n-1}$, and set
\begin{align*}
w\mathfrak{u}=\transp{(\alpha_1,\ldots,\alpha_n)},\qquad
t'=\diag(p^{r_1},\ldots,p^{r_{n-1}}),\qquad
\mathfrak{u}'=\transp{(\alpha_1,\ldots,\alpha_{n-1})}.
\end{align*}
	
The support condition is equivalent to
\begin{equation}\label{2}
p^j\begin{pmatrix}
			tu&w\mathfrak{u}\\
			&1
\end{pmatrix}\in K_p\begin{pmatrix}
			I_n\\
			&p^{-\ell}
\end{pmatrix}K_p
\end{equation}
for some $j\in\mathbb Z$. Necessarily $-\ell\leq j\leq0$.
\begin{enumerate}
\item[(1)] If $e\geq e'$, then $j-e=-\ell$. Row and column operations reduce
\eqref{2} to
\begin{align*}
p^j\begin{pmatrix}
p^{-e}t'u'&&\mathfrak{u}'-\alpha_n t'\mathfrak{u}''\\
				&p^{-e}&\alpha_n\\
				&&1
\end{pmatrix}\in K_p\begin{pmatrix}
				I_n\\
				&p^{-\ell}
\end{pmatrix}K_p.
\end{align*}
It follows that $\alpha_n\in p^{-e}\mathbb{Z}_p$ and $p^j\begin{pmatrix}
p^{-e}t'u'&\mathfrak{u}'-\alpha_n t'\mathfrak{u}''\\
			&1
		\end{pmatrix}\in K_p'$. This gives case (A), together with
$\delta_{B'}^{-1}(t)\leq p^{(n-1)\ell}$.
		
\item[(2)] If $e<e'$, then $j-e'=-\ell$. Let $m$ be the largest
index such that $e_p(\alpha_m)=-e'$.
\begin{itemize}
\item If $m=n$, row and column operations give
\begin{align*}
p^j\begin{pmatrix}
p^{-e}t'u'&p^{-e}t'\mathfrak{u}''-p^{-e}\alpha_n^{-1}\mathfrak{u}'&\\
&p^{-e}&\alpha_n\\
&&1
\end{pmatrix}\in K_p\begin{pmatrix}
I_n\\
&p^{-\ell}
\end{pmatrix}K_p.
\end{align*}
Hence $2e'-e=\ell$ and $w\mathfrak{u}=\begin{pmatrix}
\mathfrak{u}'\\
\alpha_n\end{pmatrix}$, where $\mathfrak{u}'\in p^{-e'}\mathbb{Z}_p^{n-1}$ and
$\alpha_n\in p^{-e'}\mathbb{Z}_p^{\times}$. Moreover,
\begin{align*}
wy=\begin{pmatrix}
p^{\frac{\ell-e}{2}}I_{n-1}&\\
&p^{-e}
\end{pmatrix}\begin{pmatrix}
I_{n-1}&\mathfrak{u}''\\
&1
\end{pmatrix}k,\ \mathfrak{u}''-p^{j-e}\alpha_n^{-1}\mathfrak{u}'\in \mathbb{Z}_p^{n-1},\ k\in K_p'.
\end{align*}
			
Since $2e'-e=\ell$ and $e<e'$, we have $e'<\ell$; hence
$\delta_{B'}^{-1}(t)=p^{(n-1)(\ell+e)/2}<p^{(n-1)\ell}$.
			
\item If $m<n$, the support condition implies $j-e=0$ and
$j-e+r_m+e'=0$, so that $r_m=-e'<0$, contradicting $r_m\geq0$.
\end{itemize}
Thus case (B) holds when $e<e'$.
\end{enumerate}
\end{proof}

\subsection{Diagonal Amplification Support}\label{sec6.2.}
\begin{lemma}\label{lem7.}
Let $p=p_0$ and $\ell=i$, where $0\leq i\leq\ell_{p_0}$. Let
$y\in G'(\mathbb Q_p)$ and
$\mathfrak u\in M_{n,1}(\mathbb Q_p)\simeq\mathbb Q_p^n$. Suppose that the
minimum valuations of the entries of $y$ and $\mathfrak u$ are $-e$ and
$-e'$, respectively, with $e,e'\geq0$, and write $wy=tuk$ as in \eqref{dec}.
Write $w\mathfrak u=\transp{(\mathfrak u',\alpha)}$, where
$\mathfrak u'\in\mathbb Q_p^{n-1}$ and $\alpha\in\mathbb Q_p$.
If
$f_p\left(\begin{pmatrix}
y&\mathfrak{u}\\
			&1
\end{pmatrix}\right)\neq 0$, then $\delta_{B'}^{-1}(t)\leq p^{2(n-1)\ell}$,
and the following alternatives hold.
\begin{enumerate}
\item[(A)] If $e\geq e'$, one of the following holds:
\begin{enumerate}
\item[(A.1)] $e=\ell$ and 
\begin{align*}
wy=\begin{pmatrix}
p^\ell\\
&I_{n-2}&\\
&&p^{-\ell}
\end{pmatrix}\begin{pmatrix}
1&\mathfrak{c}&\mathfrak{u}_1''\\
&I_{n-2}&\mathfrak{u}_2''\\
&&1
\end{pmatrix}k,\quad
\begin{pmatrix}
\mathfrak{u}_1''\\
\mathfrak{u}_2''
\end{pmatrix}\in \begin{pmatrix}
p^{-2\ell}\mathbb{Z}_p\\
p^{-\ell}\mathbb{Z}_p^{n-2} \end{pmatrix},
\end{align*}
$\mathfrak{c}\in p^{-\ell}\mathbb{Z}_p^{n-2}$, $k\in K_p'$, and
$\mathfrak{u}'-\alpha p^\ell\begin{pmatrix}
p^\ell\mathfrak{u}_{1}''\\
\mathfrak{u}_{2}''
\end{pmatrix}\in \mathbb{Z}_p^{n-1}$ and $\alpha \in p^{-\ell}\mathbb{Z}_p$. 
			
\item[(A.2)] $\ell<e\leq 2\ell$. Set $r_1=3\ell-e$. Then
\begin{align*}
wy=\begin{pmatrix}
p^{-e+r_1}\\
	&p^{\ell-e}I_{n-2}&\\
	&&p^{-e}
\end{pmatrix}\begin{pmatrix}
1&\mathfrak{c}&\mathfrak{u}_{1}''\\
&I_{n-2}&\mathfrak{u}_{2}''\\
&&1
\end{pmatrix}k,\quad k\in K_p',
\end{align*}
$\mathfrak{u}_{2}''\in p^{-\ell}\mathbb{Z}_p^{n-2},$ $\mathfrak{u}_{1}''\in p^{-r_1}\mathbb{Z}_p,$ $\mathfrak{c}\in p^{\ell-r_1}\mathbb{Z}_p^{n-2},$ $\alpha\in p^{-\ell}\mathbb{Z}_p,$ and 
$\mathfrak{u}'-\alpha p^\ell\begin{pmatrix}
\mathfrak{u}_{1}''\\
\mathfrak{u}_{2}''
\end{pmatrix}\in 
p^{\ell-e}\mathbb{Z}_p^{n-1}$.  
\end{enumerate}
		
\item[(B)] If $e<e'$, one of the following holds:
\begin{enumerate}
\item[(B.1)] We have $2e'-e=\ell$,
$\mathfrak{u}'\in p^{-e'}\mathbb{Z}_p^{n-1}$, and
$\alpha\in p^{-e'}\mathbb{Z}_p^{\times}$. Moreover,
\begin{align*}
wy=\begin{pmatrix}
p^{\ell+e'-e}\\
&p^{e'-e}I_{n-2}&\\
&&p^{-e}
\end{pmatrix}\begin{pmatrix}
1&\mathfrak{c}&\mathfrak{u}_{1}''\\
&I_{n-2}&\mathfrak{u}_2''\\
&&1
\end{pmatrix}k,\quad  \ k\in K_p',
\end{align*}
with $\mathfrak{c}\in p^{-\ell}\mathbb{Z}_p^{n-2},$  $\begin{pmatrix}
\mathfrak{u}_1''\\
\mathfrak{u}_2''
\end{pmatrix}-
\begin{pmatrix}
p^{-\ell}\\
&I_{n-2}
\end{pmatrix}p^{-e'}\alpha^{-1}\mathfrak{u}'\in \mathbb{Z}_p^{n-1}$.

\item[(B.2)] We have $2e'-e=2\ell$,
$\mathfrak{u}'\in p^{-e'}\mathbb{Z}_p^{n-1}$, and
$\alpha\in p^{-e'}\mathbb{Z}_p^{\times}$. Moreover,
\begin{align*}
wy=\begin{pmatrix}
p^{e'-e-\ell}\\
&p^{e'-e-\ell}I_{n-2}&\\
&&p^{-e}
\end{pmatrix}\begin{pmatrix}
1&&\mathfrak{u}_1''\\
&I_{n-2}&\mathfrak{u}_2''\\
&&1
\end{pmatrix}k,
\end{align*}
with $k\in K_p',$ and $\begin{pmatrix}
\mathfrak{u}_1''\\
\mathfrak{u}_2''
\end{pmatrix}-p^{\ell-e'}\alpha^{-1}\mathfrak{u}'\in \mathbb{Z}_p^{n-1}$.
			
\item[(B.3)] We have $\ell<2e'-e<2\ell$. Set $j=e'-\ell$ and
$r_1=3\ell-3e'+2e$. Then
$1\leq j-e+r_1<\ell$. Write $\mathfrak{u}'=\begin{pmatrix}
\mathfrak{u}_1'\\
\mathfrak{u}_2'\end{pmatrix}$; then $\alpha\in p^{-e'}\mathbb{Z}_p^{\times}$ and
\begin{align*}
wy=\begin{pmatrix}
p^{3\ell-3e'+e}\\
&p^{\ell-e'}I_{n-2}&\\
&&p^{-e}
\end{pmatrix}\begin{pmatrix}
1&\mathfrak{c}&\mathfrak{u}_1''\\
&I_{n-2}&\mathfrak{u}_2''\\
&&1
\end{pmatrix}k,
\end{align*}
$\mathfrak{c}\in p^{2e'-e-2\ell}\mathbb{Z}_p^{n-2},$ $k\in K_p',$
$\mathfrak{u}_2''\in p^{e'-e-\ell}\mathbb{Z}_p^{n-2},$
$\mathfrak{u}_1''\in p^{3e'-2e-3\ell}\mathbb{Z}_p,$ and
\begin{align*}
\begin{pmatrix}
\mathfrak{u}_1''\\
\mathfrak{u}_2''
\end{pmatrix}-
\begin{pmatrix}
p^{3e'-2e-3\ell}\\
&p^{e'-e-\ell}I_{n-2}
\end{pmatrix}\alpha^{-1}\mathfrak{u}'\in \begin{pmatrix}
p^{2e'-e-2\ell}\mathbb{Z}_p\\
\mathbb{Z}_p^{n-2}
\end{pmatrix}.
\end{align*}

\end{enumerate}
\end{enumerate}
\end{lemma}
\begin{proof}
Retain the notation from the proof of Lemma \ref{lem7}; in particular,
write $w\mathfrak{u}=\transp{(\alpha_1,\ldots,\alpha_n)}$,
$t'=\diag(p^{r_1},\ldots,p^{r_{n-1}})$, and
$\mathfrak{u}'=\transp{(\alpha_1,\ldots,\alpha_{n-1})}$.

The support condition is equivalent to
\begin{equation}\label{2.}
p^j\begin{pmatrix}
tu&w\mathfrak{u}\\
&1
\end{pmatrix}\in K_p\begin{pmatrix}
			p^{\ell}\\
			&I_{n-1}\\
			&&p^{-\ell}
		\end{pmatrix}K_p.
\end{equation}
for some $j\in\mathbb Z$. Necessarily $-\ell\leq j\leq \ell$.
\begin{enumerate}
\item[(A)] If $e\geq e'$, then $j-e=-\ell$, and row and column operations
reduce \eqref{2.} to
\begin{align*}
p^j\begin{pmatrix}
p^{-e}t'u'&&\mathfrak{u}'-\alpha_nt'\mathfrak{u}''\\
&p^{-e}&\alpha_n\\
&&1
\end{pmatrix}\in K_p\begin{pmatrix}
p^{\ell}\\
&I_{n-1}\\
&&p^{-\ell}
\end{pmatrix}K_p.
\end{align*}
Thus $\alpha_n\in p^{-e}\mathbb{Z}_p$ and
\begin{equation}\label{5}
p^j\begin{pmatrix}
p^{-e}t'u'&\mathfrak{u}'-\alpha_nt'\mathfrak{u}''\\
				&1
\end{pmatrix}\in K_p'\begin{pmatrix}
				p^{\ell}\\
				&I_{n-1}
\end{pmatrix}K_p'.
\end{equation} 
		
\begin{itemize}
\item[(A.1)] If $j-e+r_1=\ell$, use the identity
			\begin{equation}\label{3}
				\begin{pmatrix}
					a&aC&b\\
					&D&F\\
					&&1
				\end{pmatrix}^{-1}=\begin{pmatrix}
					a^{-1}&-CD^{-1}&-a^{-1}b+CD^{-1}F\\
					&D^{-1}&-D^{-1}F\\
					&&1
				\end{pmatrix},
			\end{equation}
to invert the matrix in \eqref{5}. Its $(n,n)$-entry forces $j=0$, and the
remaining entries give case (A.1). Moreover,
$\delta_{B'}^{-1}(t)\leq p^{2(n-1)\ell}$.

\item[(A.2)] If $j-e+r_1<\ell$, inversion of \eqref{5} using \eqref{3} gives
\begin{align*}
wy=\begin{pmatrix}
p^{-j-\ell+r_1}\\
&p^{-j}I_{n-2}&\\
&&p^{-j-\ell}
\end{pmatrix}\begin{pmatrix}
1&\mathfrak{c}&\mathfrak{u}_1''\\
&I_{n-2}&\mathfrak{u}_2''\\
&&1
\end{pmatrix}k,
\end{align*}
$k\in K_p',$	$-j+e\leq r_1<\ell-j+e,$ $r_2=\cdots=r_{n-1}=\ell,$ $\mathfrak{u}_2''\in p^{-\ell}\mathbb{Z}_p^{n-2},$ $\mathfrak{u}_1''\in p^{-r_1}\mathbb{Z}_p,$ $\mathfrak{c}\in p^{-j+e-r_1}\mathbb{Z}_p^{n-2},$  and $w\mathfrak{u}=\begin{pmatrix}
\mathfrak{u}'\\
\alpha_n\end{pmatrix}$ with $\mathfrak{u}'-\alpha_np^\ell\mathfrak{u}''\in 
p^{\ell-e}\mathbb{Z}_p^{n-1}$ and $\alpha_n\in p^{-\ell}\mathbb{Z}_p$. 
			
Comparison of determinants gives $3\ell=e+r_1$. Since $r_1\geq r_2=\ell$,
$j-e=-\ell$, and $-\ell\leq j\leq \ell$, we obtain $\ell\leq r_1<2\ell$ and
$e=j+\ell\leq2\ell$. This proves case (A.2), as well as
$\delta_{B'}^{-1}(t)<p^{2(n-1)\ell}$.
\end{itemize}

\item[(B)] If $e<e'$, then $j-e'=-\ell$. The pivot argument from the proof of
Lemma \ref{lem7} gives $e_p(\alpha_n)=-e'$, and
row and column operations reduce \eqref{2.} to
\begin{align*}
p^j\begin{pmatrix}
p^{-e}t'u'&p^{-e}t'\mathfrak{u}''-p^{-e}\alpha_n^{-1}\mathfrak{u}'&\\
				&p^{-e}&\alpha_n\\
				&&1
\end{pmatrix}\in K_p\begin{pmatrix}
				p^\ell\\
				&I_{n-1}\\
				&&p^{-\ell}
			\end{pmatrix}K_p.
\end{align*}
Consequently,
\begin{equation}\label{10}
p^j\begin{pmatrix}
p^{-e}t'u'&p^{-e}t'\mathfrak{u}''-p^{-e}\alpha_n^{-1}\mathfrak{u}'\\
				&p^{e'-e}
			\end{pmatrix}\in K_p'\begin{pmatrix}
				p^{\ell}\\
				&I_{n-1}
			\end{pmatrix}K_p'.
\end{equation}
		
\begin{itemize}
\item[(B.1)] If $j-e+r_1=\ell$, inversion of \eqref{10} gives
$j+e'-e=0$ and $r_2=\cdots=r_{n-1}=e'$. Hence \eqref{3} yields
\begin{align*}
wy=\begin{pmatrix}
p^{\ell+e'-e}\\
&p^{e'-e}I_{n-2}&\\
&&p^{-e}
\end{pmatrix}\begin{pmatrix}
1&\mathfrak{c}&\mathfrak{u}_1''\\
&I_{n-2}&\mathfrak{u}_2''\\
&&1
\end{pmatrix}k,
\end{align*}
$\mathfrak{c}\in p^{-\ell}\mathbb{Z}_p^{n-2},$ $k\in K_p',$ and $w\mathfrak{u}=\begin{pmatrix}
\mathfrak{u}'\\
\alpha_n\end{pmatrix},$ $\mathfrak{u}'\in p^{-e'}\mathbb{Z}_p^{n-1},$ $\alpha_n\in p^{-e'}\mathbb{Z}_p^{\times},$ $\mathfrak{u}''-t'^{-1}\alpha_n^{-1}\mathfrak{u}'\in \mathbb{Z}_p^{n-1}$. This is case (B.1).
			
Here $\ell>e'>e\geq0$, so $\delta_{B'}^{-1}(t)\leq p^{2(n-1)\ell}$.
			
\item[(B.2)] If $j+e'-e=\ell$, then $2e'-e=2\ell$. Inverting
\eqref{10} by \eqref{3} gives $0\leq r_1=\cdots=r_{n-1}=e'-\ell$ and
\begin{align*}
wy=\begin{pmatrix}
p^{e'-e-\ell}\\
&p^{e'-e-\ell}I_{n-2}&\\
&&p^{-e}
\end{pmatrix}\begin{pmatrix}
1&&\mathfrak{u}_1''\\
&I_{n-2}&\mathfrak{u}_2''\\
&&1
\end{pmatrix}k,\ k\in K_p',
\end{align*}
and $w\mathfrak{u}=\begin{pmatrix}
				\mathfrak{u}'\\
\alpha_n\end{pmatrix},$ $\mathfrak{u}'\in p^{-e'}\mathbb{Z}_p^{n-1},$ $\alpha_n\in p^{-e'}\mathbb{Z}_p^{\times},$ $\mathfrak{u}''-p^{\ell-e'}\alpha_n^{-1}\mathfrak{u}'\in \mathbb{Z}_p^{n-1}$. This is case (B.2).
			
Since $2e'-e=2\ell$ and $e<e'$, we have $e'<2\ell$. Therefore
$\delta_{B'}^{-1}(t)=p^{(n-1)(e'-\ell)}<p^{2(n-1)\ell}$.
			
\item[(B.3)] If $j-e+r_1<\ell$ and $j+e'-e<\ell$, then, by the argument in
case (B) of Lemma \ref{lem7}, the $(1,n)$-entry of the inverse of the
left-hand side of \eqref{10} has minimum valuation $-\ell$. Moreover,
\eqref{10} is equivalent to
\begin{align*}
p^{j+e'-e}\begin{pmatrix}
p^{-e'}t'u'&p^{-e'}t'\mathfrak{u}''-p^{-e'}\alpha_n^{-1}\mathfrak{u}'\\
&1
\end{pmatrix}\in K_p'\begin{pmatrix}
p^{\ell}\\
&I_{n-1}
\end{pmatrix}K_p',
\end{align*}
which has the form \eqref{5}. If $n=2$, comparison of determinants in
\eqref{10} gives
\begin{align*}
2j-2e+e'+r_1=\ell.
\end{align*}
If $n\geq3$, the argument of case (A.2) gives the same identity, together
with $j-e+r_m=0$ for some $1<m<n$. Since $j-e'=-\ell$, it follows in either
case that $3e'-2e+r_1=3\ell$ and
$-e+r_1=3\ell-3e'+e$. Write $w\mathfrak{u}=\begin{pmatrix}
\mathfrak{u}'\\
\alpha_n\end{pmatrix},$ where $\mathfrak{u}'=\begin{pmatrix}
\mathfrak{u}_1'\\
\mathfrak{u}_2'\end{pmatrix},$ $\alpha_n\in p^{-e'}\mathbb{Z}_p^{\times}$. Then  
\begin{align*}
wy=\begin{pmatrix}
p^{-e+r_1}\\
&p^{-j}I_{n-2}&\\
&&p^{-e}
\end{pmatrix}\begin{pmatrix}
1&\mathfrak{c}&\mathfrak{u}_1''\\
&I_{n-2}&\mathfrak{u}_2''\\
&&1
\end{pmatrix}k,\ k\in K_p',
\end{align*}
$\mathfrak{c}\in p^{-j+e-r_1}\mathbb{Z}_p^{n-2},$
$\mathfrak{u}_2''\in p^{e'-e-\ell}\mathbb{Z}_p^{n-2},$
$\mathfrak{u}_1''\in p^{3e'-2e-3\ell}\mathbb{Z}_p,$ and
\begin{align*}
\mathfrak{u}''-\begin{pmatrix}
p^{3e'-2e-3\ell}\\
&p^{e'-e-\ell}I_{n-2}
\end{pmatrix}\alpha_n^{-1}\mathfrak{u}'\in \begin{pmatrix}
p^{2e'-e-2\ell}\mathbb{Z}_p\\
\mathbb{Z}_p^{n-2}
\end{pmatrix}.
\end{align*}
This is case (B.3). Finally, $j-e>j-e'=-\ell$ and $j-e+r_1<\ell$, so
\begin{align*}
\delta_{B'}^{-1}(t)=p^{(n-1)r_1}<p^{2(n-1)\ell}.
\end{align*}
\end{itemize}	
\end{enumerate}
\end{proof}
\begin{remark}
Lemmas \ref{lem7} and \ref{lem7.} provide the local support and volume
constraints used to estimate every amplified orbital integral; see
\textsection\ref{sec5.4} and \textsection\ref{8.5.2}--\textsection\ref{sec10}.
Exploiting these constraints collectively yields the saving
$\delta\asymp n^{-3}$ in Corollary \ref{cor1.2}, compared with
$\delta\asymp n^{-5}$ in \cite{Nel21}.
\end{remark}

\section{Geometric Side: The Small Cell Orbital Integral}\label{8.5.1}
Let $\mathbf{s}=(s,0)\in\mathbb{C}^2$. By
\textsection\ref{2.2.1}, the small-cell orbital integral
$J^{\Reg}_{\Geo,\sm}(f,\mathbf{s})$ converges absolutely for $\Re(s)>0$ and
extends meromorphically to $\mathbb{C}$, with possible simple poles at
$s=0,-1$. Define
\begin{equation}\label{7.1'}
J_{\Geo,\Main}^{\sm}(f,\mathbf{s}):=J^{\Reg}_{\Geo,\sm}(f,\mathbf{s})-s^{-1}\underset{s=0}{\Res}\ J^{\Reg}_{\Geo,\sm}(f,\mathbf{s}).
\end{equation}

Then $J_{\Geo,\Main}^{\sm}(f,\mathbf{s})$ is holomorphic for $\Re(s)>-1$.
Our main estimate is the following.

\begin{prop}\label{prop54}
Let $f\in\big\{f_{p_0,i}^{\diag},f_{p_1,p_2}^{\mathrm{off}}\big\}$ and
$\varepsilon>0$, and set $\vartheta_f:=\max_{p\mid\nu(f)}\vartheta_p$, where $\vartheta_p$ is the admissible Ramanujan exponent for $\pi_p'$ fixed
in \textsection\ref{sec3.1.3}. Then
\begin{align*}	
J_{\Geo,\Main}^{\sm}(f,\mathbf{0})
\ll  T^{\frac{n}{2}+\varepsilon}M'^{2n}M^{n+\varepsilon} \mathcal{N}_f^{-1+2\vartheta_f+\varepsilon}\langle\phi',\phi'\rangle\prod_{p\mid M'}p^{ne_p(M)},
\end{align*}
where $\mathcal N_f$ is defined by \eqref{61}. The implied constant may
depend on $\varepsilon$, $c_\infty$, $C_\infty$, and the fixed data attached
to $\pi'$ and $p_\circ$.
\end{prop}

The proof is assembled in \textsection\ref{sec5.4} from the local estimates
below.

\subsection{Finite-Place Integrals}\label{6.2..}
For a finite prime $p$ and $\Re(s)>-1+2\vartheta_p$, define
\begin{equation}\label{eq7.2}
\mathcal{I}_p(f,s)=\int_{N'(\mathbb{Q}_p)\backslash G'(\mathbb{Q}_p)}\int_{G'(\mathbb{Q}_p)}\kappa(x,y)W_p'(x)\overline{W_p'(xy)}|\det x|_p^{1+s}dxdy,
\end{equation}
where
\begin{align*}
\kappa(x,y):=\int_{M_{n,1}(\mathbb{Q}_p)} f_p\left(\begin{pmatrix}
		y&\mathfrak{u}\\
		&1
	\end{pmatrix}\right)\psi_p(\eta x\mathfrak{u})d\mathfrak{u}.
\end{align*}
Here $W_p'$ is the local Whittaker vector fixed in
\textsection\ref{sec4.}; it is spherical and normalized by $W_p'(I_n)=1$
when $p\nmid M'$.

Write $y=w^{-1}tuk$ as in \eqref{dec}. At $p\mid\nu(f),$ $W_p'$ is spherical. Then 
\begin{equation}\label{6.7}
	\int_{K_p'}\overline{W_p'(bk'y)}dk'=\gamma_{\pi_p'}(t)\overline{W_p'(b)},
\end{equation}
where $\gamma_{\pi_p'}(t)$ is the normalized spherical function (see \cite{Mac71}). 

\begin{lemma}\label{lem7.2}
Suppose that $p\mid\nu(f)$, so that $p\nmid M'$. Then
\begin{multline*}
\mathcal{I}_p(f,0)=\frac{L_p(1,\pi_p'\times\widetilde\pi_p')}{\zeta_p(n)}
\int_{G'(\mathbb{Q}_p)}\int_{M_{n,1}(\mathbb{Q}_p)}\sum_{m\geq 0}p^{-mn}f_p\left(\begin{pmatrix}
		y&\mathfrak{u}\\
		&1
\end{pmatrix}\right)\\
\int_{\mathbb{Z}_p^{\times}}\psi_p(\eta p^m\beta \mathfrak{u})
d^{\times}\beta\,d\mathfrak{u}\,\gamma_{\pi_p'}(t)dy.
\end{multline*}	
\end{lemma}
\begin{proof}
Let $H=\diag(\mathrm{GL}(n-1),1)$ and $N_H=N'\cap H$. Write
$x=p^m\beta bk'$ with $m\in\mathbb Z$, $\beta\in\mathbb Z_p^\times$,
$b\in N_H(\mathbb Q_p)\backslash H(\mathbb Q_p)$, and $k'\in K_p'$.
Since $|\det(p^mI_n)|_p=p^{-mn}$, unfolding gives
\begin{multline*}
\mathcal{I}_p(f,0)=\int_{K_p'}\int_{G'(\mathbb{Q}_p)}\int_{M_{n,1}(\mathbb{Q}_p)}\sum_{m\in\mathbb{Z}}p^{-mn}f_p\left(\begin{pmatrix}
		y&\mathfrak{u}\\
		&1
\end{pmatrix}\right)\int_{\mathbb{Z}_p^{\times}}\psi_p(\eta p^m\beta k'\mathfrak{u})d^{\times}\beta d\mathfrak{u}\\
\int_{N_H(\mathbb{Q}_p)\backslash H(\mathbb{Q}_p)}W_p'(bk')\overline{W_p'(bk'y)}dbdydk'.
\end{multline*}	
At $s=0$, the Kirillov inner product is invariant under right translation by
$K_p'$. Moreover, left invariance under $K_p(M)$ and orthogonality under
$\mathfrak u\mapsto\mathfrak u+\mathbb Z_p^n$ show that only $m\geq0$
contributes. Changing variables
$\mathfrak u\mapsto k'^{-1}\mathfrak u$ and $y\mapsto k'^{-1}y$ gives
\begin{multline*}
\mathcal{I}_p(f,0)=\int_{G'(\mathbb{Q}_p)}\int_{M_{n,1}(\mathbb{Q}_p)}\sum_{m\geq 0}p^{-mn}f_p\left(\begin{pmatrix}
		y&\mathfrak{u}\\
		&1
\end{pmatrix}\right)\int_{\mathbb{Z}_p^{\times}}\psi_p(\eta p^m\beta \mathfrak{u})d^{\times}\beta d\mathfrak{u} \\
\int_{N_H(\mathbb{Q}_p)\backslash H(\mathbb{Q}_p)}W_p'(b)\Big[\int_{K_p'}\overline{W_p'(bk'y)}dk'\Big]dbdy.
\end{multline*}	
Applying \eqref{6.7} and the standard Rankin--Selberg norm identity
\begin{align*}
\int_{N_H(\mathbb Q_p)\backslash H(\mathbb Q_p)}|W_p'(b)|^2db
=\zeta_p(n)^{-1}L_p(1,\pi_p'\times\widetilde\pi_p')
\end{align*}
proves the lemma.
\end{proof}
\begin{remark}
The argument is special to $s=0$: for $s\neq0$, the weighted Kirillov pairing
\begin{align*}
\int_{N_H(\mathbb{Q}_p)\backslash H(\mathbb{Q}_p)}
W_p'(bk')\overline{W_p'(bk'y)}|\det b|_p^sdb
\end{align*} 
is no longer invariant under right translation by $k'$.
\end{remark}

\begin{lemma}\label{lem6.4.}
Let $p\nmid p_\circ M'\nu(f)$. Then
\begin{align*}
\mathcal{I}_p(f,s)=\Vol(\overline{K_p(M)})^{-1}\cdot L_p(1+s,\pi_p'\times\widetilde{\pi}_p'), \ \ \Re(s)>-1+2\vartheta_p.
\end{align*}
\end{lemma}
\begin{proof}

Unfolding as in the proof of Lemma \ref{lem7.2} gives
\begin{align*}
\sum_{m\geq 0}\frac{1}{p^{mn(1+s)}}\int_{K_p'}\int_{G'(\mathbb{Q}_p)}&\int_{\mathbb{Q}_p^n}f_p\left(\begin{pmatrix}
		y&\mathfrak{u}\\
		&1
	\end{pmatrix}\right)\int_{\mathbb{Z}_p^{\times}}\psi_p(\eta p^m\beta k'\mathfrak{u})d^{\times}\beta d\mathfrak{u} \Phi(s)dydk',
\end{align*}	
where   
\begin{align*}
\Phi(s):=\int_{N_H(\mathbb{Q}_p)\backslash H(\mathbb{Q}_p)}W_p'(b) \int_{K_p'}\overline{W_p'(bk'y)}dk' |\det b|_p^sdb.
\end{align*}
 
Since $p\nmid p_\circ\nu(f)$, the function $f_p\left(\begin{pmatrix}
		y&\mathfrak{u}\\
		&1
\end{pmatrix}\right)$ is supported on $y\in K_p'$ and
$\mathfrak u\in\mathbb Z_p^n$. Hence
\begin{align*}
\mathcal{I}_p(f,s)=\frac{1}{\Vol(\overline{K_p(M)})}
\sum_{m\geq 0}p^{-mn(1+s)}
\int_{N_H(\mathbb{Q}_p)\backslash H(\mathbb{Q}_p)}
|W_p'(b)|^2|\det b|_p^sdb.
\end{align*} 
The sum over $m$ is $\zeta_p(n(1+s))$, whereas the last integral is
$L_p(1+s,\pi_p'\times\widetilde\pi_p')/\zeta_p(n(1+s))$. Their product
gives the asserted identity.
\end{proof}

\begin{lemma}\label{lem7.aux.upper}
Let $p=p_\circ$. Then $\mathcal{I}_p(f,s)$ is holomorphic in a fixed
neighborhood of $s=0$, and
\begin{align*}
|\mathcal{I}_p(f,0)|+
\Big|\frac{d \mathcal{I}_p(f,s)}{ds}\Big|_{s=0}\Big|\ll 1.
\end{align*}
The implied constant may depend on $n$, $p_\circ$, and $\pi_{p_\circ}'$.
\end{lemma}
\begin{proof}
Since $p_\circ\nmid M'$, the representation $\pi_{p_\circ}'$ is unramified.
Unfolding as in Lemma \ref{lem7.2}, the standard Whittaker bound and
$2\vartheta_{p_\circ}<1$ give an absolutely convergent majorant throughout
a fixed strip $|\Re(s)|<\eta$, with $\eta>0$. The same majorant remains
summable after inserting the linear valuation weight produced by
differentiation. Thus $\mathcal I_p(f,s)$ is holomorphic in that strip and
its value and first derivative at $s=0$ are bounded. Since both the prime
and the local data are fixed, the bound is $O(1)$. This is the same
fixed-place argument as in \cite[Lemma 4.5]{Yan23b}.
\end{proof}

\begin{lemma}\label{lem11}
	Let notation be as before. Let $p\mid M'$. Then for $\Re(s)>2\vartheta_p-1,$ 
\begin{equation}\label{equ7.5}
|\mathcal{I}_p(f,s)|\ll p^{2ne_p(M)+2ne_p(M')-e_p(M'')}L_p(1+\Re(s),\pi_p'\times\widetilde{\pi}_p'),
\end{equation}
where the implied constant is absolute. In particular, 
\begin{equation}\label{eq7.5.}
|\mathcal{I}_p(f,0)|+\Big|\frac{d\mathcal{I}_p(f,s)}{ds}\Big|_{s=0}\Big|\ll p^{2ne_p(M)+2ne_p(M')-e_p(M'')}L_p(1,\pi_p'\times\widetilde{\pi}_p'),
\end{equation}
where the implied constant depends only on $n$.
\end{lemma}
\begin{proof}
By Cauchy--Schwarz,
\begin{align*}
|\mathcal{I}_p(f,s)|^2\leq  \mathcal{I}_p^{(1)}(f,s)\cdot \mathcal{I}_p^{(2)}(f,s)
\end{align*}
where 
\begin{align*}
\mathcal{I}_p^{(1)}(f,s):=&\int_{N'(\mathbb{Q}_p)\backslash G'(\mathbb{Q}_p)}\int_{G'(\mathbb{Q}_p)}|\kappa(x,y)||W_p'(x)|^2|\det x|_p^{1+\Re(s)}dxdy,\\
\mathcal{I}_p^{(2)}(f,s):=&\int_{N'(\mathbb{Q}_p)\backslash G'(\mathbb{Q}_p)}\int_{G'(\mathbb{Q}_p)}|\kappa(x,y)||W_p'(xy)|^2|\det x|_p^{1+\Re(s)}dxdy.
\end{align*}

Write $x=p^mak'\in N'(\mathbb{Q}_p)\backslash G'(\mathbb{Q}_p)$, where $a\in A'(\mathbb{Q}_p),$ and $k'\in K_p'$.
Then 
\begin{align*}
\kappa(x,y)=\int_{M_{n,1}(\mathbb{Q}_p)}f_p\left(\begin{pmatrix}
			y&\mathfrak{u}\\
			&1
		\end{pmatrix}\right)\psi_p(\eta p^mk'\mathfrak{u})d\mathfrak{u}.
\end{align*}
By the change of variable $\mathfrak{u}\mapsto y\mathfrak{u},$ we have
\begin{equation}\label{eq7.5}
\frac{\kappa(xy^{-1},y)}{|\det y|_p^{1+\Re(s)}}=|\det y|_p^{-\Re(s)}\int_{M_{n,1}(\mathbb{Q}_p)} f_p\left(\begin{pmatrix}
		y&y\mathfrak{u}\\
		&1
	\end{pmatrix}\right)\psi_p(\eta x\mathfrak{u})d\mathfrak{u}.
\end{equation}

Recall the explicit formula for $f_p$ in \eqref{3.13}.

Note that $f_p$ is bi-$\begin{pmatrix}
	I_n& M_{n,1}(\mathbb{Z}_p)\\
	&1
\end{pmatrix}$-invariant. By orthogonality, $\kappa(x,y)=0$ or $\kappa(xy^{-1},y)=0$ unless $m\geq 0,$ which amounts to $\eta x\in M_{1,n}(\mathbb{Z}_p)$. So
\begin{align*}
\mathcal{I}_p^{(1)}(f,s)=&\int\int_{G'(\mathbb{Q}_p)}|\kappa(x,y)|dy|W_p'(x)|^2|\det x|_p^{1+\Re(s)}\textbf{1}_{M_{1,n}(\mathbb{Z}_p)}(\eta x)dx,\\
\mathcal{I}_p^{(2)}(f,s)=&\int\int_{G'(\mathbb{Q}_p)}\frac{|\kappa(xy^{-1},y)|}{|\det y|_p^{1+\Re(s)}}dy|W_p'(x)|^2|\det x|_p^{1+\Re(s)}\textbf{1}_{M_{1,n}(\mathbb{Z}_p)}(\eta x)dx,
\end{align*}
where $x$ ranges through $N'(\mathbb{Q}_p)\backslash G'(\mathbb{Q}_p)$. 

Investigating the support of $f_p,$ we have $\kappa(xy^{-1},y)=0$ or $\kappa(x,y)=0$ unless $y\in K_p'$ and $\mathfrak{u}\in p^{-m'}M_{n,1}(\mathbb{Z}_p)$.  Hence, by triangle inequality,
\begin{align*}
\int_{G'(\mathbb{Q}_p)}|\kappa(x,y)|dy
&\leq \int_{K_p'}\int_{M_{n,1}(\mathbb{Q}_p)}\left|f_p\left(\begin{pmatrix}
y&\mathfrak{u}\\
			&1
		\end{pmatrix}\right)\right|d\mathfrak{u}dy\\
&\leq v_p^2\Vol(\overline{\supp f_p})\leq v_p^2p^{2nm'},
\end{align*}
which is $\ll p^{2ne_p(M)+2nm'-m''}$. Here the implied constant is absolute. Therefore,
\begin{multline*}
\mathcal{I}_p^{(1)}(f,s)
\ll p^{2ne_p(M)+2nm'-m''}
\int_{N'(\mathbb{Q}_p)\backslash G'(\mathbb{Q}_p)}|W_p'(x)|^2\\
\times|\det x|_p^{1+\Re(s)}
\textbf{1}_{M_{1,n}(\mathbb{Z}_p)}(\eta x)dx,
\end{multline*}
By \cite[Theorem 4.1]{Miy14}, this gives
\begin{align*}
\mathcal{I}_p^{(1)}(f,s)\ll p^{2ne_p(M)+2nm'-m''}L_p(1+\Re(s),\pi_p'\times\widetilde{\pi}_p'). 
\end{align*}

Likewise, by \eqref{eq7.5} and the change of variable $\mathfrak{u}\mapsto y^{-1}\mathfrak{u},$ we have 
\begin{align*}
\int_{G'(\mathbb{Q}_p)}\frac{|\kappa(xy^{-1},y)|}{|\det y|_p^{1+\Re(s)}}dy\leq \int_{K_p'}\int_{M_{n,1}(\mathbb{Q}_p)}\left|f_p\left(\begin{pmatrix}
			y&\mathfrak{u}\\
			&1
		\end{pmatrix}\right)\right|d\mathfrak{u}dy,
\end{align*}
which is $\ll p^{2ne_p(M)+2nm'-m''}$. As a consequence, 
\begin{align*}
\mathcal{I}_p^{(2)}(f,s)\ll p^{2ne_p(M)+2nm'-m''}L_p(1+\Re(s),\pi_p'\times\widetilde{\pi}_p'). 
\end{align*}

This proves \eqref{equ7.5}. Cauchy's formula on the circle
$|s|=n^{-100}$ gives
\begin{align*}
\Big|\frac{d\mathcal{I}_p(f,s)}{ds}\Big|_{s=0}\Big|\leq \frac{1}{2\pi}\int_{|s|=n^{-100}}\frac{|\mathcal{I}_p(f,s)|}{|s|^2}|ds|\ll p^{2ne_p(M)+2nm'-m''},
\end{align*}
where the implied constant depends only on $n$. The local Rankin--Selberg
factor is uniformly bounded on this circle by the admissible Ramanujan bound.
Together with $L_p(1,\pi_p'\times\widetilde\pi_p')\gg1$, this proves
\eqref{eq7.5.}.
\end{proof}
\begin{remark}
These bounds are more than sufficient because $\pi'$ is fixed.
\end{remark}

\begin{comment} 
and the spherical Whittaker function $W_p'$ supports in 
$$
a_p\in\mathcal{A}_p:=\big\{a_p=(p^{\alpha_1},\cdots, p^{\alpha_{n-1}},1),\ \alpha_1,\cdots, \alpha_{n-1}\in\mathbb{Z},\ 
\alpha_1\geq\cdots\geq \alpha_{n-1}\geq 0\big\}.
$$
Recall the Macdonald formula 
\begin{equation}\label{6}
	W'_p(p^{\boldsymbol{\ell}}n_p'k_p')=\delta_{B'}^{\frac{1}{2}}(p^{\boldsymbol{\ell}})\mathbf{S}_{\pi'_p}(p^{\boldsymbol{\ell}})W'_p(I_n),
\end{equation}
where $p^{\boldsymbol{\ell}}\in p^{\mathbb{Z}^{n}}$, $n_p'\in N'(\mathbb{Q}_p),$ $K_p'\in K_p',$ and $\mathbf{S}_{\pi'_p}(p^{\boldsymbol{\ell}})$ is the Schur polynomial and $\mathbf{S}_{\pi'_p}(p^{\boldsymbol{\ell}})=0$ unless $p^{\boldsymbol{\ell}}\in p^{\mathbb{Z}^{n}}$ is dominant. %Writing Schur polynomials in terms of semistandard Young tableaux we see that the coefficient of any monomial appearing in a Schur polynomial is non-negative. 
\end{comment}

\subsection{Amplification Primes}
We now estimate $\mathcal{I}_p(f,0)$ for $p\mid\nu(f)$.
  
\subsubsection{Notation}
Let $W_H\subset W_{G'}$ be generated by $w_1,\ldots,w_{n-2}$, so that
$W_{G'}=W_H\sqcup W_Hw_{n-1}W_H$.

Fix $f\in\{f_{p_0,i}^{\diag},f_{p_1,p_2}^{\mathrm{off}}\}$ and a prime
$p\mid\nu(f)$. Let $\mathcal T_p$ be the set of classes
$t\in Z'(\mathbb Q_p)\backslash T'(\mathbb Q_p)$ for which
\begin{align*}
f_p\left(\begin{pmatrix}
	w^{-1}tu&\mathfrak{u}\\
	&1
\end{pmatrix}\right)\neq 0
\end{align*}
for some $w\in W_{G'}$, $u\in N'(\mathbb Q_p)$, and
$\mathfrak u\in M_{n,1}(\mathbb Q_p)$. Lemmas \ref{lem7} and \ref{lem7.}
give $\#\mathcal T_p\ll_n1$.

For fixed $t,u,w$, let $\mathfrak U_p(t,u,w)$ be the set of
$\mathfrak u\in M_{n,1}(\mathbb Q_p)$ such that
\begin{align*}
f_p\left(\begin{pmatrix}
	w^{-1}tu& \mathfrak{u}\\
	&1
\end{pmatrix}\right)\neq 0.
\end{align*}
Set
\begin{align*}
\mathcal{J}_p^{\dagger}(t,u,w):=\sum_{m\geq 0}p^{-mn}\int_{\mathbb{Z}_p^{\times}}\int_{\mathfrak U_p(t,u,w)}\psi_p(\eta p^m\beta \mathfrak{u})d\mathfrak{u}d^{\times}\beta.
\end{align*}

With the measure normalization of \textsection\ref{notation}, the local
Ramanujan sum is
\begin{equation}\label{8}
\int_{\mathbb{Z}_p^{\times}}\psi_p(p^r\beta)d^{\times}\beta=\begin{cases}
1,& r\geq0,\\
-(p-1)^{-1},& r=-1,\\
0,& r\leq-2.
	\end{cases}
\end{equation}

\subsubsection{Off-diagonal terms}
Write $y=w^{-1}tuk$ as in \eqref{dec}, retaining the notation of Lemma
\ref{lem7}. Let $\mathcal U_p(t)$ denote the corresponding range of $u$, where $u=\begin{pmatrix}
I_{n-1}&\mathfrak{u}''\\
&1
\end{pmatrix}$.

\begin{lemma}[Off-diagonal integrals]\label{lem8}
Let $f=f_{p_1,p_2}^{\mathrm{off}}$, let $j\in\{1,2\}$, and set
$p=p_j$ and $\ell=\ell_p$. Then
	\begin{equation}\label{43}
		\mathcal{I}_p(f,0) \ll p^{(-1/2+\vartheta_p)\ell},
	\end{equation}
where the implied constant depends only on $n$.
\end{lemma}
\begin{proof}
The changes of variables $x\mapsto xy^{-1}$ and
$\mathfrak u\mapsto y\mathfrak u$ reduce the case $p=p_1$ to $p=p_2$.
We therefore assume $p=p_2$.

If $w\in W_H$, then \eqref{8} gives
$\mathcal J_p^\dagger(t,u,w)\ll1$. Hence Lemmas \ref{lem7} and
\ref{lem7.2} give
\begin{align*}
\sum_{w\in W_H}\sum_{t}\int_{\mathcal U_p(t)}|\gamma_{\pi_p'}(t)|\|f_p\|_{\infty}du
\ll p^{\vartheta_p \ell-\frac{n\ell}{2}}
\sum_{t}\delta_{B'}^{-\frac{1}{2}}(t)
\ll p^{(-\frac{1}{2}+\vartheta_p)\ell}.
\end{align*}	
Here we used $\|f_p\|_\infty\ll p^{-n\ell/2}$,
$\#\mathcal T_p\ll_n1$, Lemma \ref{lem7}, and Macdonald's bound
\begin{equation}\label{6.10}
|\gamma_{\pi_p'}(t)|\ll p^{\ell\vartheta_p}\cdot \delta_{B'}^{1/2}(t)
\end{equation}
along with $\Vol(K_p'tK_p')\ll\delta_{B'}^{-1}(t)$.
	
Now suppose $w\notin W_H$. Represent $t\in\mathcal T_p$ by
$\begin{pmatrix}
		I_{n-1}\\
		&p^{-a}
\end{pmatrix}$ modulo $Z'(\mathbb Q_p)$, where $0\leq a\leq\ell$.
Equation \eqref{8} and Lemma \ref{lem7} imply
\begin{equation}\label{6.10.}
\int_{\mathcal U_p(t)}\mathcal{J}_p^{\dagger}(t,u,w)du\ll p^{(n-1)a}.
\end{equation} 
	
Combining \eqref{6.10}, \eqref{6.10.}, and
$\|f_p\|_\infty\ll p^{-n\ell/2}$, we obtain
	\begin{align*}
\sum_{w\notin W_H}\sum_t|\gamma_{\pi_p'}(t)|\|f_p\|_\infty
\left|\int_{\mathcal U_p(t)}\mathcal J_p^\dagger(t,u,w)du\right|\ll p^{\vartheta_p\ell-\frac{n\ell}{2}}
p^{-\frac{(n-1)\ell}{2}}p^{(n-1)\ell}\ll p^{(-\frac12+\vartheta_p)\ell}.
	\end{align*}
This proves \eqref{43}.
\end{proof}

\subsubsection{Diagonal terms}
Retain the notation of Lemma \ref{lem7.}, write $y=w^{-1}tuk$, and let
$\mathcal U_p(t)$ denote the corresponding range of
$u=\begin{pmatrix}
1&\mathfrak{c}&\mathfrak{u}_1''\\
&I_{n-2}&\mathfrak{u}_2''\\
&&1
\end{pmatrix}$.

\begin{lemma}[Diagonal integrals]\label{lem9}
Let $f=f_{p_0,i}^{\diag}$, set $p=p_0$, and write $\ell=i$. Then
	\begin{equation}\label{13}
		\mathcal{I}_p(f,0)\ll p^{(-1+2\vartheta_p)\ell},
	\end{equation}
where the implied constant depends only on $n$.
\end{lemma}
\begin{proof}
If $w\in W_H$, then \eqref{8} implies that
$\mathcal J_p^\dagger(t,u,w)\ll1$. Combining Lemmas \ref{lem7.} and
\ref{lem7.2}, we obtain
\begin{align*}
\sum_{w\in W_H}\sum_t\int_{\mathcal U_p(t)}
|\gamma_{\pi_p'}(t)|\|f_p\|_\infty\,du
\ll
p^{2\vartheta_p\ell-n\ell}
\sum_t\delta_{B'}^{-1/2}(t)
\ll
p^{(-1+2\vartheta_p)\ell}.
\end{align*}	
Here we used Macdonald's bound
\begin{equation}\label{7.7}
|\gamma_{\pi_p'}(t)|\ll p^{2\ell\vartheta_p}\cdot \delta_{B'}^{1/2}(t)
\end{equation}
and $\Vol(K_p'tK_p')\ll\delta_{B'}^{-1}(t)$.

Suppose now that $w\notin W_H$. By Lemma \ref{lem7.}, every
$t\in\mathcal T_p$ has a representative
$\diag(p^b,I_{n-2},p^{-a})$ with $0\leq a,b\leq\ell$. We estimate the five
cases of that lemma separately.

Let $\mathcal{I}_p(*)$ denote the contribution of case
$*\in\{A.1,A.2,B.1,B.2,B.3\}$. By Lemma \ref{lem7.2},
\begin{equation}\label{6.13}
\mathcal{I}_p(*) \ll\sum_{w}\sum_{t}|\gamma_{\pi_p'}(t)|\|f_p\|_{\infty}\Big|\int_{\mathcal U_p(t)}\mathcal{J}_p^{\dagger}(t,u,w)\textbf{1}_{*}(y)du\Big|,
\end{equation}
where $\textbf{1}_{*}$ is the indicator function of case $(*)$.
	
Recall that $w\mathfrak{u}=\begin{pmatrix}
\mathfrak{u}'\\
	\alpha
\end{pmatrix}$. Write $\mathfrak{u}=\transp{(\alpha_1',\cdots,\alpha_n')}$. Let
$r_w$ be the position of $\alpha_n'$ in $\mathfrak{u}'$, and write the
corresponding entry of $\mathfrak{u}'':=\begin{pmatrix}
		\mathfrak{u}_1''\\
	\mathfrak{u}_2''
		\end{pmatrix}$ as $p^{\mu_{r_w}}\beta'$, where
$\beta'\in\mathbb{Z}_p^{\times}$.
	
\begin{enumerate}
\item[(A.1):] In the case (A.1) of Lemma \ref{lem7.}, $\mathcal{T}_p\ni t\in \diag(p^\ell,I_{n-2},p^{-\ell})Z'(\mathbb{Q}_p)$. By orthogonality, we derive that   $\mathcal{J}_p^{\dagger}(t,u,w)\ll p^{\ell},$ where the implied constant is absolute. In conjunction with \eqref{6.13} we obtain that 
\begin{align*}
\mathcal{I}_p(A.1)\ll \sum_{w}\sum_{t}|\gamma_{\pi_p'}(t)|\|f_p\|_{\infty}\int_{\mathcal U_p(t)}du\cdot p^{\ell}\ll p^{2\vartheta_p \ell-\ell}.
\end{align*}	
Here we use \eqref{7.7} and the facts that $\|f_p\|_{\infty}\ll p^{-n\ell}$ and 
$\Vol(\mathcal U_p(t))\ll p^{(2n-3)\ell},$ which follows from Lemma \ref{lem7.}.
		
\item[(A.2)] In the case (A.2) of Lemma \ref{lem7.}, $\mathcal{T}_p\ni t\in \diag(p^{2\ell-e},I_{n-2},p^{-\ell})Z'(\mathbb{Q}_p),$ where $\ell<e\leq 2\ell$. Note that $\alpha\in p^{-\ell}\mathbb{Z}_p,$ and 
		$\mathfrak{u}'-\alpha p^\ell\mathfrak{u}''\in 
		p^{\ell-e}\mathbb{Z}_p^{n-1}$. By \eqref{8} and the parametrization of $\mathcal U_p(t)$ we have 
\begin{equation}\label{6.14}
\int_{\mathcal U_p(t)}\mathcal{J}_p^{\dagger}(t,u,w)du\ll  p^{\ell-e}\cdot p^\ell\cdot p^{(n-1)(2\ell-e)}\cdot p^{(n-2)\ell},
\end{equation}
Here the four factors arise, respectively, from $\alpha$, the integration in
$\mathfrak u'$ together with the $r_w$-th oscillatory integral, the range of
$\mathfrak u''$, and the range of $\mathfrak c$ in Lemma \ref{lem7.}.
		
Therefore, by \eqref{6.13} and the fact that $\ell<e\leq 2\ell,$ we derive  
		\begin{align*}
			\mathcal{I}_p(A.2)
\ll& p^{2\vartheta_p \ell}p^{-n\ell}\sum_w\sum_{t}\delta_{B'}^{1/2}(t)\cdot p^{\ell-e}\cdot p^\ell\cdot p^{(n-1)(2\ell-e)}\cdot p^{(n-2)\ell},
\end{align*}
Since $\delta_{B'}^{-1}(t)=p^{(n-1)(3\ell-e)}$ and $e>\ell$, the
right-hand side is $\ll p^{(-1+2\vartheta_p)\ell}$.

\item[(B.1)] In the case (B.1) of Lemma \ref{lem7.}, $\mathcal{T}_p\ni t\in \diag(p^{\ell},I_{n-2},p^{-e'})Z'(\mathbb{Q}_p),$ where $e'<\ell$. Note that $\alpha\in p^{-e'}\mathbb{Z}_p^{\times},$ and 
		$\mathfrak{u}''-\begin{pmatrix}
			p^{-\ell}\\
			&I_{n-2}
		\end{pmatrix}p^{-e'}\alpha^{-1}\mathfrak{u}'\in \mathbb{Z}_p^{n-1}$. We distinguish whether $r_w=1$.
		
\begin{itemize}
\item Suppose $r_w=1$. Then the $r_w$-th entry of $\mathfrak{u}'$ lies in $p^{\mu_{r_w}}p^{e'+\ell}\alpha\beta'+p^{e'+\ell}\mathbb{Z}_p$. By orthogonality,
			\begin{align*}
\mathcal{J}_p^{\dagger}(t,u,w)\ll p^{e'}\cdot p^{-\ell}\cdot\Big|\sum_{m\geq 0}p^{-mn}\int_{\mathbb{Z}_p^{\times}}\psi_p( p^{m+\ell+\mu_{r_w}}\beta)d^{\times}\beta\Big|.
\end{align*}
where the factor $p^{e'}$ comes from the contribution from $\alpha,$ and $p^{-\ell}$ is the contribution from $\alpha_1',\cdots,\alpha_{n-1}'$.
			
\item Suppose $r_w>1$. Then the $r_w$-th entry of $\mathfrak{u}'$ lies in $p^{\mu_{r_w}}p^{e'}\alpha\beta'+p^{e'}\mathbb{Z}_p$. By orthogonality,
\begin{align*}
\mathcal{J}_p^{\dagger}(t,u,w)\ll p^{e'}\cdot \Big|\sum_{m\geq 0}p^{-mn}\int_{\mathbb{Z}_p^{\times}}\psi_p(p^{m+\mu_{r_w}}\beta)d^{\times}\beta\Big|.
\end{align*}
where the factor $p^{e'}$ comes from the contribution from $\alpha$. 	
\end{itemize}
		
Using \eqref{8} to handle the above integrals relative to $\beta$, we then obtain 
\begin{align*}
	\int_{\mathcal U_p(t)}\mathcal{J}_p^{\dagger}(t,u,w)du\ll p^{e'}\cdot p^{(n-1)(\ell+e')-e'}. 
\end{align*} 
Together with \eqref{6.13} and $e'<\ell$, this gives
\begin{align*}
\mathcal{I}_p(B.1)\ll& p^{2\vartheta_p \ell}p^{-n\ell}\sum_{w}\sum_{t}\delta_{B'}^{1/2}(t)\cdot p^{e'}\cdot p^{(n-1)(\ell+e')-e'}\ll p^{(-1+2\vartheta_p)\ell}.
\end{align*}
Here $\delta_{B'}^{-1}(t)=p^{(n-1)(\ell+e')}$ and $e'<\ell$.
		
\item[(B.2)] In the case (B.2) of Lemma \ref{lem7.}, $\mathcal{T}_p\ni t\in \diag(I_{n-1},p^{\ell-e'})Z'(\mathbb{Q}_p),$ where $\ell\leq e'<2\ell$. Note that $\alpha\in p^{-e'}\mathbb{Z}_p^{\times},$ and 
		$\mathfrak{u}''-p^{\ell-e'}\alpha^{-1}\mathfrak{u}'\in \mathbb{Z}_p^{n-1}$. So
\begin{align*}
\mathcal{J}_p^{\dagger}(t,u,w)\ll p^{e'}\cdot p^{-\ell}\cdot\Big|\sum_{m\geq \ell}p^{-mn}\int_{\mathbb{Z}_p^{\times}}\psi_p(p^{m-\ell+\mu_{r_w}}\beta)d^{\times}\beta\Big|,
\end{align*}
as in the case $r_w=1$ above.
		
Applying \eqref{8} to the $\beta$-integral and using \eqref{6.13},
\begin{align*}
\mathcal{I}_p(B.2)\ll& p^{2\vartheta_p \ell}p^{-n\ell}\sum_w\sum_{t}\delta_{B'}^{1/2}(t)\cdot p^{e'-\ell}\cdot p^{(n-1)(e'-\ell)-(e'-\ell)}\ll p^{(-1+2\vartheta_p)\ell}.
\end{align*}
Here $\delta_{B'}^{-1}(t)=p^{(n-1)(e'-\ell)}$ and $e'<2\ell$.
		
\item[(B.3)] In the case (B.3) of Lemma \ref{lem7.},
\begin{align*}
t\in \diag(p^{2\ell-2e'+e},I_{n-2},p^{e'-e-\ell})Z'(\mathbb{Q}_p),
\end{align*}
where $\ell<2e'-e<2\ell$ and $e'>e$. Moreover,
$\alpha\in p^{-e'}\mathbb{Z}_p^{\times}$ and
\begin{equation}\label{30}
\mathfrak{u}''-\begin{pmatrix}
p^{3e'-2e-3\ell}\\
&p^{e'-e-\ell}I_{n-2}
\end{pmatrix}\alpha^{-1}\mathfrak{u}'\in \begin{pmatrix}
p^{2e'-e-2\ell}\mathbb{Z}_p\\
\mathbb{Z}_p^{n-2}
\end{pmatrix}.
\end{equation}

We again distinguish whether $r_w=1$.
		
\begin{itemize}
\item Suppose $r_w=1$. According to \eqref{30}, the $r_w$-th entry of $\mathfrak{u}'$ ranges through $p^{\mu_{r_w}}p^{2e+3\ell-3e'}\alpha\beta'+p^{e+\ell-2e'}\mathbb{Z}_p$. By orthogonality,
\begin{align*}
\mathcal{J}_p^{\dagger}(t,u,w)\ll p^{e'}\cdot p^{e+\ell-2e'}\cdot\Big|\sum_{m\geq 2e'-e-\ell}p^{-mn}\int_{\mathbb{Z}_p^{\times}}\psi_p(p^{m+2e+3\ell-4e'+\mu_{r_w}}\beta)d^{\times}\beta\Big|,
\end{align*}
where the factor $p^{e'}$ comes from the contribution from $\alpha,$ and $p^{e+\ell-2e'}$ is the contribution from $\alpha_1',\cdots,\alpha_{n-1}'$. 
			
\item Suppose $r_w>1$. According to \eqref{30}, the $r_w$-th entry of $\mathfrak{u}'$ ranges through $p^{\mu_{r_w}}p^{e+\ell-2e'}\beta'+p^{e+\ell-2e'}\mathbb{Z}_p$. By orthogonality,
			\begin{align*}
\mathcal{J}_p^{\dagger}(t,u,w)\ll p^{e'}\cdot p^{e+\ell-2e'}\cdot\Big|\sum_{m\geq 2e'-e-\ell}p^{-mn}\int_{\mathbb{Z}_p^{\times}}\psi_p(p^{m+e+\ell-2e'+\mu_{r_w}}\beta)d^{\times}\beta\Big|.
\end{align*}
\end{itemize}
		
Using \eqref{8} to handle the above integrals relative to $\beta$, we then obtain 
\begin{align*}
\int_{\mathcal U_p(t)}\mathcal{J}_p^{\dagger}(t,u,w)du\ll p^{e'}\cdot p^{e+\ell-2e'}\cdot p^{(n-1)(2\ell-2e'+e)+e'-e-\ell}.
\end{align*} 		
		
Set $A:=2\ell-2e'+e$. By case (B.3) of Lemma \ref{lem7.},
$0<A<\ell$. Since $\delta_{B'}^{-1}(t)=p^{(n-1)A}$, \eqref{6.13}
gives
\begin{align*}
\mathcal{I}_p(B.3)
\ll p^{2\vartheta_p\ell-n\ell+\frac{n-1}{2}A}
\ll p^{(-1+2\vartheta_p)\ell}.
\end{align*}
\end{enumerate}

These estimates prove \eqref{13}.
\end{proof}

\subsubsection{Derivatives at amplification primes}
We next record the corresponding estimate for the derivative at $s=0$.

\begin{lemma}\label{lem7.7.}
Let $p\mid\nu(f)$. Set $\ell=\ell_p$ if
$f=f_{p_1,p_2}^{\mathrm{off}}$, and $\ell=i$ if
$f=f_{p_0,i}^{\diag}$. Then $\mathcal I_p(f,s)$ is holomorphic in a fixed
neighborhood of $s=0$, and
\begin{equation}\label{7.11}
\frac{d \mathcal{I}_p(f,s)}{ds}\Big|_{s=0} \ll
\begin{cases}
p^{(-1/2+\vartheta_p)\ell+\varepsilon},&
\text{if $f=f_{p_1,p_2}^{\mathrm{off}}$},\\
p^{(-1+2\vartheta_p)\ell+\varepsilon},&
\text{if $f=f_{p_0,i}^{\diag}$},
\end{cases} 
\end{equation}
where the implied constant depends at most on $n$ and $\varepsilon$.
\end{lemma}
\begin{proof}
Unfold as in Lemma \ref{lem7.2}, but retain the factor
$|\det x|_p^s$. The Cartan expansion is absolutely convergent for
$|s|<\eta_n$, where $\eta_n>0$ is fixed, because the admissible Ramanujan
bound gives a uniform gap $1-2\vartheta_p\gg_n1$. It may therefore be
differentiated termwise at $s=0$.

Writing the central exponent as $m$ and the dominant Cartan exponent as
$\lambda$, differentiation inserts a factor bounded by
\begin{align*}
(1+m+\|\lambda\|)\log p.
\end{align*}
The Casselman--Shalika and Macdonald bounds used in Lemmas \ref{lem8} and
\ref{lem9} majorize the remaining sums by geometric series whose ratios are
at most $p^{-c_n}$ for some $c_n>0$. Their first moments therefore satisfy
the same bounds, up to
$(1+\ell)^{O_n(1)}\log p$. The support restrictions and the Ramanujan-sum
estimates are unchanged. Consequently, the arguments of Lemmas \ref{lem8}
and \ref{lem9} give their respective bounds multiplied by
\begin{align*}
(1+\ell)^{O_n(1)}\log p\ll_{n,\varepsilon}p^\varepsilon,
\end{align*}
which proves \eqref{7.11}.
\end{proof}

\subsection{The Archimedean Integral}\label{sec5.3}

We now estimate the Archimedean factor using the construction of $f_\infty$;
see \cite[\S\S 14.5--14.6]{Nel23}.

\begin{lemma}\label{lem10}
Let $\mathcal I_\infty(f,s)$ denote the Archimedean analogue of
\eqref{eq7.2}. It is holomorphic in a fixed neighborhood of $s=0$, and
\begin{align*}
|\mathcal I_\infty(f,0)|+|\mathcal I_\infty'(f,0)|
\ll_\varepsilon\mathcal{A}_\infty
\ll_\varepsilon T^{\frac n2+\varepsilon}\|W_\infty'\|_2^2,
\end{align*}
where
\begin{multline*}
\mathcal{A}_\infty:=
\sum_{\delta\in\{-\varepsilon,\varepsilon\}}
\int_{K_{\infty}'}\int_{G'(\mathbb{R})}\int_{Z'(\mathbb{R})}
\Bigg|\int_{\mathbb{R}^n}f_{\infty}\left(\begin{pmatrix}
			y&u\\
			&1
\end{pmatrix}\right)\psi_{\infty}(\eta zku)du\Bigg||\det z|^{1+\varepsilon}_{\infty}d^{\times}z\\
\times\int_{N'(\mathbb{R})\backslash P'_0(\mathbb{R})}\big|W_{\infty}'(pk)\overline{W_{\infty}'(pky)}\big||\det p|_{\infty}^{1+\delta}d^{\times}pdydk
\end{multline*}
Here $\|W_{\infty}'\|_2$ is defined by \eqref{W_inf}.
\end{lemma}
\begin{proof}
For $k\in K_{\infty}'$, write $\eta k=(k_1,\ldots,k_n)$. By the construction
of $f_\infty$, the summand of $\mathcal{A}_\infty$ corresponding to
$\delta\in\{-\varepsilon,\varepsilon\}$ is, up to $O(T^{-\infty})$, bounded by
\begin{align*}
\int_{K_{\infty}'}\int_{G'(\mathbb{R})}\int_{\mathbb{R}^{\times}}\Bigg|\int_{\mathbb{R}^n}\tilde{f}_{\infty}^{\sharp}\left(\begin{pmatrix}
			y&u\\
			&1
\end{pmatrix}\right)\psi_{\infty}(\eta tku)du\Bigg| |t|^{n+n\varepsilon}d^{\times}t\,\mathcal{K}_\delta(k,y)\,dy\,dk,
\end{align*}
where $\tilde{f}_{\infty}^{\sharp}$ is defined in
\textsection\ref{3.2.1}; it is the untruncated operator
$\widetilde{Op}_h(a')$ of \cite[\S 14.6]{Nel23}, and
\begin{align*}
\mathcal{K}_\delta(k,y):=\int_{N'(\mathbb{R})\backslash P'_0(\mathbb{R})}
\big|W_{\infty}'(pk)\overline{W_{\infty}'(pky)}\big|
|\det p|_{\infty}^{1+\delta}\,d^{\times}p.
\end{align*}

Up to an error $O(T^{-\infty})$, we may restrict to the essential support
of $\tilde f_\infty^\sharp$. In particular, $y$ ranges over a fixed compact
set. For every fixed $A\in\mathbb R$, the
Kirillov--Plancherel argument used in the proof of Lemma \ref{lem4.6},
together with \cite[Lemma 21.27]{Nel21}, gives
\begin{align*}
\sup_{g\in C}
\int_{N'(\mathbb R)\backslash P_0'(\mathbb R)}
|W_\infty'(pg)|^2|\det p|_\infty^{1+A}\,d^\times p
\ll_{A,C}T^{o(1)}\|W_\infty'\|_2^2
\end{align*}
for every fixed compact set $C\subset G'(\mathbb R)$. Indeed, multiplication
by $|\det|^{A/2}$ is the determinant weight in the Kirillov model, and
$\det(T^{\rho^\vee})=1$, so the rescaling used to construct $W_\infty'$
does not alter this variable.

Applying this with $A=\delta$ to a compact set containing both $k$ and $ky$,
and then using Cauchy--Schwarz, gives, uniformly in $k$, $y$, and
$\delta\in\{-\varepsilon,\varepsilon\}$,
\begin{align*}
\mathcal{K}_\delta(k,y)\ll T^{o(1)}\|W_\infty'\|_2^2
\ll T^\varepsilon\|W_\infty'\|_2^2.
\end{align*}
Consequently,
\begin{multline*}
\mathcal{A}_\infty\ll T^\varepsilon\|W_\infty'\|_2^2
\int_{K_{\infty}'}\int_{G'(\mathbb{R})}\int_{\mathbb{R}^{\times}}\Bigg|\int_{\mathbb{R}^n}\tilde{f}_{\infty}^{\sharp}\left(\begin{pmatrix}
		y&u\\
		&1
\end{pmatrix}\right)\psi_{\infty}(\eta tku)du\Bigg|\\
\times|t|^{n+n\varepsilon}d^{\times}t\,dy\,dk+O(T^{-\infty}).
\end{multline*}

Set
\begin{align*}
\mathcal F(y,t,k):=
\int_{\mathbb R^n}\tilde f_\infty^\sharp
\left(\begin{pmatrix}y&u\\&1\end{pmatrix}\right)
\psi_\infty(\eta tku)\,du.
\end{align*}
The support of the partial Fourier transform in $u$ implies that
$\mathcal F(y,t,k)=O(T^{-\infty})$ unless
$|tk_j-T\tau_j|\ll T^{\frac{1}{2}+\varepsilon}$ for $1\leq j\leq n$,
where $(\tau_1,\ldots,\tau_n)$ is determined by $\tau$. The essential-support
estimate for $\tilde f_\infty^\sharp$ also gives
\begin{align*}
\mathcal F(y,t,k)=\int_{\mathcal{U}}\tilde{f}_{\infty}^{\sharp}\left(\begin{pmatrix}
		y&u\\
		&1
	\end{pmatrix}\right)\psi_{\infty}(\eta tku)du+O(T^{-\infty}),
\end{align*}
where
\begin{align*}
\mathcal U:=\big\{\transp{(u_1,\ldots,u_n)}\in\mathbb R^n:
|u_j|\ll T^{-1/2+\varepsilon},\ 1\leq j\leq n\big\}.
\end{align*}
	
Absorbing the sign of $t$ into $k$, we may take $t>0$. The pushforward of
Haar measure on $K_\infty'$ under $k\mapsto\eta k$ is the invariant measure
on $\mathbb S^{n-1}$. After absorbing factors of size $T^{O(\varepsilon)}$,
we obtain
\begin{multline*}
\mathcal{A}_{\infty}\ll_{\varepsilon}T^{\varepsilon}\|W_\infty'\|_2^2
\int_{G'(\mathbb{R})}\int_{\mathcal{U}}\bigg|\tilde{f}_{\infty}^{\sharp}\left(\begin{pmatrix}
			y&u\\
			&1
		\end{pmatrix}\right)\bigg|du\,dy\\
\times\int_0^\infty
\int_{\mathbb{S}_{T,\tau}^{n-1}(t)}d\mathbf{k}\,t^{n-1}dt
+O(T^{-\infty}),
\end{multline*}
where $\mathbb S_{T,\tau}^{n-1}(t)$ is defined as
\begin{equation}\label{98}
\big\{(k_1,\ldots,k_n)\in\mathbb R^n:\
k_1^2+\cdots+k_n^2=1,\
|tk_j-T\tau_j|\ll T^{\frac12+\varepsilon},\ 1\leq j\leq n\big\}.
\end{equation}

Put $r=T^{1/2+\varepsilon}$ and
$\tau_{\mathrm{row}}=(\tau_1,\ldots,\tau_n)$. By stability,
$\|\tau_{\mathrm{row}}\|\asymp1$. If
$\mathbb S_{T,\tau}^{n-1}(t)$ is nonempty, then $\big|t-T\|\tau_{\mathrm{row}}\|\big|\ll r$. 
Thus the relevant $t$-range has length $O(r)$ and satisfies $t\asymp T$.
For each such $t$, the set $\mathbb S_{T,\tau}^{n-1}(t)$ is contained in a
spherical cap of angular radius $O(r/t)$. Therefore
\begin{equation}\label{7.23}
\int_0^\infty
\Vol\big(\mathbb S_{T,\tau}^{n-1}(t)\big)t^{n-1}\,dt
\ll r^n\ll T^{\frac n2+n\varepsilon}.
\end{equation}

The support and size estimates in the construction of
$\tilde f_\infty^\sharp$ give
\begin{align*}
\Vol(\mathcal U)\ll T^{-\frac{n}{2}+\varepsilon},\qquad 
\int_{G'(\mathbb R)}
\sup_{u\in\mathcal U}
\left|\tilde f_\infty^\sharp
\left(\begin{pmatrix}y&u\\&1\end{pmatrix}\right)\right|dy \ll T^{\frac n2+\varepsilon}.
\end{align*}
The second estimate follows from the microlocal support condition
\eqref{245}, stability of $\tau$, and the Fourier-inversion estimate underlying
\eqref{250}; it is the corresponding uniform form of the argument in
\cite[\S 14.9]{Nel23}. Consequently,
\begin{align*}
\int_{G'(\mathbb R)}\int_{\mathcal U}
\left|\tilde f_\infty^\sharp
\left(\begin{pmatrix}y&u\\&1\end{pmatrix}\right)\right|du\,dy
\ll T^\varepsilon.
\end{align*}
Combining this estimate with \eqref{7.23} gives
\begin{align*}
\mathcal{A}_\infty
\ll_\varepsilon
T^{\frac n2+\varepsilon}\|W_\infty'\|_2^2,
\end{align*}
as required. Finally, differentiation at $s=0$ inserts
$\log|\det x|_\infty$. On the effective support, $|t|\asymp T$, while
\begin{align*}
1+|\log v|\ll_\varepsilon v^\varepsilon+v^{-\varepsilon}
\qquad (v>0).
\end{align*}
Thus $\mathcal{A}_\infty$ dominates both the value and the first derivative
of the local factor. The same majorants, with slightly smaller
$\varepsilon$, also prove holomorphy near $s=0$.
\end{proof}
\begin{remark}
The function $f_\infty$ is not $\iota(K_\infty')$-invariant, so the
spherical argument used at the non-Archimedean places is unavailable here.
\end{remark}

\subsection{Proof of Proposition \ref{prop54}}\label{sec5.4}
Let
$S=\{\infty,p_\circ\}\bigsqcup\,\{p:p\mid\nu(f)M'\}$ and
$\mathbf{s}=(s,0)$. By \textsection\ref{2.2.1} and Lemma
\ref{lem6.4.},
\begin{equation}\label{7.1}
J^{\Reg}_{\Geo,\sm}(f,\mathbf{s})=\prod_{\substack{p\mid M,\ p\nmid M'}}\Vol(\overline{K_p(M)})^{-1}\cdot L^{(S)}(1+s,\pi'\times\widetilde{\pi}')\cdot \mathcal{I}_S(f,s), 
\end{equation}
where $L^{(S)}$ denotes the partial $L$-function with the factors in $S$
removed and $\mathcal{I}_S(f,s)=\prod_{p\in S}\mathcal{I}_p(f,s)$.
For $p=\infty$, the definition of $\mathcal{I}_p$ is the Archimedean analogue
of \eqref{eq7.2}, with $\mathbb Q_p=\mathbb R$.

Fix $0<\varepsilon<4/(n(n+1)+2)$. Equation \eqref{7.1} gives the
meromorphic continuation of $J^{\Reg}_{\Geo,\sm}(f,\mathbf{s})$ to
$\Re(s)>-\varepsilon$, while $\mathcal{I}_S(f,s)$ is holomorphic there.

Write the Laurent and Taylor expansions at $s=0$ as
\begin{align*}
L^{(S)}(1+s,\pi'\times\widetilde{\pi}')=\frac{a_{-1}}s+a_0+a_1s+\cdots,\quad 
\mathcal{I}_S(f,s)=b_0+b_1s+b_2s^2+\cdots.
\end{align*}
Then
\begin{align*}
J^{\Reg}_{\Geo,\sm}(f,\mathbf{s})=\prod_{p\mid M,\ p\nmid M'}\Vol(\overline{K_p(M)})^{-1}\cdot\Big[\frac{a_{-1}b_0}{s}+a_0b_0+a_{-1}b_1+O(s)\Big].
\end{align*} 
Consequently,
\begin{align*}
J_{\Geo,\Main}^{\sm}(f,\mathbf{0})
=(a_0b_0+a_{-1}b_1)\prod_{p\mid M,\ p\nmid M'}
\Vol(\overline{K_p(M)})^{-1}.
\end{align*}

Standard bounds at $s=1$ give
$|a_{-1}|+|a_0|\ll_\varepsilon C(\pi')^\varepsilon$. Moreover,
$b_0=\mathcal{I}_S(f,0)$ and
$b_1=\mathcal{I}_S'(f,0)$. Combining Lemmas \ref{lem11},
\ref{lem7.aux.upper}, \ref{lem8}, \ref{lem9}, \ref{lem7.7.}, and
\ref{lem10}, and absorbing the number of prime divisors into
$M^\varepsilon M'^\varepsilon$, gives
\begin{multline*}
\prod_{p\mid M,\ p\nmid M'}
\Vol(\overline{K_p(M)})^{-1}
(|a_{-1}|+|a_0|)(|b_0|+|b_1|)\\ 
\ll_{\varepsilon}
T^{\frac n2+\varepsilon}M'^{2n}M^{n+\varepsilon}
\mathcal N_f^{-1+2\vartheta_f+\varepsilon}
\langle\phi',\phi'\rangle
\prod_{p\mid M'}p^{ne_p(M)}.
\end{multline*}
This proves Proposition \ref{prop54}.

\section{Geometric Side: The Dual Orbital Integral}\label{8.5.2}
Let $\mathbf{s}=(s,0)\in\mathbb{C}^2$. Recall that
\begin{align*}
J^{\Reg}_{\Geo,\du}(f,\mathbf{s})
:=
\int_{G'(\mathbb A)}\int_{[\overline{G'}]}
\phi'(x)\overline{\phi'(xy)}E(s,x;f,y)\,dx\,dy,
\end{align*}
where $E(s,x;f,y)$ is defined in \textsection\ref{2.2.1}. This integral
converges absolutely for $\Re(s)>1$ and extends meromorphically to
$\mathbb C$, with possible simple poles at $s=0,1$. Define
\begin{equation}\label{8.1'}
J_{\Geo,\Main}^{\du}(f,\mathbf{s}):=J^{\Reg}_{\Geo,\du}(f,\mathbf{s})-s^{-1}\underset{s=0}{\Res}\ J^{\Reg}_{\Geo,\du}(f,\mathbf{s}),
\end{equation}
which is holomorphic for $\Re(s)<1$. We prove the following estimate.

\begin{prop}\label{prop8.1}
Let $f\in\big\{f_{p_0,i}^{\diag},f_{p_1,p_2}^{\mathrm{off}}\big\}$ and
$\varepsilon>0$, and set
$\vartheta_f:=\max_{p\mid\nu(f)}\vartheta_p$. Then
\begin{align}\label{8.1}
J_{\Geo,\Main}^{\du}(f,\mathbf{0})
\ll
T^{\frac n2+\varepsilon}M'^{2n}M^{n+\varepsilon}
\mathcal N_f^{-1+2\vartheta_f+\varepsilon}
\langle\phi',\phi'\rangle
\prod_{p\mid M'}p^{ne_p(M)}.
\end{align}
Here $\mathcal N_f$ is defined by \eqref{61}. The implied constant may
depend on $\varepsilon$, $c_\infty$, $C_\infty$, and the fixed data attached
to $\pi'$ and $p_\circ$.
\end{prop}
\begin{proof}
For $\Re(s)<0$, the functional equation of the Eisenstein series and
unfolding give
\begin{multline*}
J^{\Reg}_{\Geo,\du}(f,\mathbf{s})
=\int_{P_0'(\mathbb Q)\backslash G'(\mathbb A)}
\int_{G'(\mathbb A)}\int_{M_{1,n}(\mathbb A)}
f\left(
\begin{pmatrix}I_n&\\ \mathfrak c&1\end{pmatrix}
\begin{pmatrix}y&\\ &1\end{pmatrix}
\right)
\psi(\eta x\transp{\mathfrak c})\\
\times\overline{\phi'(x)}\phi'(x\transp{y}^{-1})
|\det x|^{1-s}\,d\mathfrak c\,dy\,dx.
\end{multline*}
It follows that
\begin{align*}
J^{\Reg}_{\Geo,\du}(f,\mathbf{s})
=\prod_{p\leq\infty}\mathcal J_p(f,s),\qquad \Re(s)<0,
\end{align*}
where, with $\mathbb Q_\infty=\mathbb R$,
\begin{multline*}
\mathcal J_p(f,s)
:=\int_{N'(\mathbb Q_p)\backslash G'(\mathbb Q_p)}
\int_{G'(\mathbb Q_p)}\int_{M_{1,n}(\mathbb Q_p)}
f_p\left(
\begin{pmatrix}I_n&\\ \mathfrak c&1\end{pmatrix}
\begin{pmatrix}y&\\ &1\end{pmatrix}
\right)
\psi_p(\eta x\transp{\mathfrak c})\\
\times\overline{W_p'(x)}W_p'(x\transp{y}^{-1})
|\det x|_p^{1-s}\,d\mathfrak c\,dy\,dx.
\end{multline*}

Set $S=\{\infty,p_\circ\}\bigsqcup\{p:p\mid\nu(f)M'\}$. The unramified calculation of Lemma \ref{lem6.4.} gives
\begin{align*}
\mathcal J_p(f,s)
=\Vol(\overline{K_p(M)})^{-1}
L_p(1-s,\pi_p'\times\widetilde\pi_p'),\qquad p\notin S.
\end{align*}
Therefore
\begin{align*}
J^{\Reg}_{\Geo,\du}(f,\mathbf{s})
=\prod_{\substack{p\mid M,\ p\nmid M'}}
\Vol(\overline{K_p(M)})^{-1}
L^{(S)}(1-s,\pi'\times\widetilde\pi')\mathcal J_S(f,s),
\end{align*}
where $\mathcal J_S(f,s):=\prod_{p\in S}\mathcal J_p(f,s)$.

To compare the remaining local factors with those in
\textsection\ref{6.2..}, define $f_p^\vee(g):=f_p(g^\vee)$. 
The changes of variables
$y\mapsto\transp{y}^{-1}$ and
$\mathfrak u=-\transp{\mathfrak c}$ transform the test-function term into $f_p^\vee\left(\begin{pmatrix}y&\mathfrak u\\ &1\end{pmatrix}\right)$. Thus the absolute-value estimates for $\mathcal J_p(f,s)$ are precisely the
estimates of \textsection\ref{6.2..} with $s$ replaced by $-s$, the additive
character replaced by its inverse, and $f_p$ replaced by $f_p^\vee$.
These changes do not affect the support or size estimates: the diagonal
amplifier is preserved up to $K_p$-double translation, while the two factors
of an off-diagonal amplifier are interchanged. The estimates at $p_\circ$,
at the ramified primes, and at infinity are likewise unchanged. In
particular, $\mathcal J_S(f,s)$ is holomorphic near $s=0$.

Write
\begin{align*}
L^{(S)}(1-s,\pi'\times\widetilde\pi')=-\frac{a_{-1}}s+a_0-a_1s+\cdots,\quad
\mathcal J_S(f,s)=c_0+c_1s+c_2s^2+\cdots.
\end{align*}
Taking the finite part at $s=0$ gives
\begin{align*}
J_{\Geo,\Main}^{\du}(f,\mathbf 0)
=(a_0c_0-a_{-1}c_1)
\prod_{\substack{p\mid M,\ p\nmid M'}}
\Vol(\overline{K_p(M)})^{-1}.
\end{align*}
The standard bounds at $s=1$ give
$|a_{-1}|+|a_0|\ll_\varepsilon C(\pi')^\varepsilon$. Applying
Lemmas \ref{lem11}, \ref{lem7.aux.upper}, \ref{lem8}, \ref{lem9},
\ref{lem7.7.}, and \ref{lem10} to the transformed local factors yields
\begin{multline*}
\prod_{\substack{p\mid M,\ p\nmid M'}}
\Vol(\overline{K_p(M)})^{-1}
(|a_{-1}|+|a_0|)(|c_0|+|c_1|)\\
\ll_\varepsilon
T^{\frac n2+\varepsilon}M'^{2n}M^{n+\varepsilon}
\mathcal N_f^{-1+2\vartheta_f+\varepsilon}
\langle\phi',\phi'\rangle
\prod_{p\mid M'}p^{ne_p(M)}.
\end{multline*}
This proves \eqref{8.1}.
\end{proof}

\section{Geometric Side: The Regular Orbital Contribution}\label{sec10}
By \cite[Theorem 6.6]{Yan22}, and after the change of variables
$\boldsymbol\xi\mapsto(1-t)\boldsymbol\xi$, which is valid since $t\neq1$,
the regular orbital contribution admits the absolutely convergent expansion
\begin{align*}
J^{\Reg,\RNum{2}}_{\Geo,\bi}(f,\mathbf{0})=\sum_{\substack{(\transp{\boldsymbol{\xi}},t)\in \mathbb{Q}^n\\
(\transp{\boldsymbol{\xi}},t)\neq\mathbf 0\\
t\neq 1}}
\int_{P_0'(\mathbb Q)\backslash G'(\mathbb A)}
\int_{G'(\mathbb A)}
f\left(\iota(x)^{-1}\gamma(\boldsymbol\xi,t)\iota(xy)\right)
\phi'(x)\overline{\phi'(xy)}\,dy\,dx,
\end{align*}
where $\gamma(\boldsymbol\xi,t):=
\begin{pmatrix}
I_{n-1}&&(1-t)\boldsymbol\xi\\
&1&t\\
&1&1
\end{pmatrix}$. This formula already combines the
regular orbits arising from the orbit-alteration procedure. We exploit this
collective structure to reduce the estimate to a simultaneous count of the
rational parameters $(\boldsymbol{\xi},t)$.

The main result of this section is the following.

\begin{thm}\label{thm9.1}
Let the notation be as before. Then
\begin{align*}
\frac{J^{\Reg,\RNum{2}}_{\Geo,\bi}(f,\mathbf{0})}
{\langle\phi',\phi'\rangle}
\ll
(TMM')^{\varepsilon}T^{\frac{n-1}{2}}M'^{2n}M^n
\mathcal N_f^{n-2}
\max_{1\leq\ell\leq n}
\left(\frac{M'^2\mathcal N_f^2}{M^\dagger}\right)^{\ell}\prod_{p\mid M'}p^{ne_p(M)},
\end{align*}
where $M^{\dagger}:=\prod_{p\mid M,\ p\nmid M'}p^{e_p(M)}$. The implied constant may depend on $\varepsilon$, $c_\infty$, $C_\infty$,
and the fixed data attached to $\pi'$ and $p_\circ$.
\end{thm}
\begin{remark}
When $M'$ is fixed and $M^\dagger\asymp M$, the maximum is attained at an
endpoint, and the bound contains the two terms
$M^{n-1}\mathcal N_f^{n}$ and $\mathcal N_f^{3n-2}$. The saving in the
first term is essential for the level-aspect nonvanishing argument in
\textsection\ref{sec12}. More precisely, writing
$X=M'^2\mathcal N_f^2/M^\dagger$, one has
$\max_{1\leq\ell\leq n}X^\ell=\max\{X,X^n\}$. The two endpoints correspond
to one and $n$ active rational coordinates, respectively. Thus this factor is
sharp for the lattice count based solely on the local support conditions;
any further improvement would require additional relations among the rational
parameters or cancellation between regular orbits.
\end{remark}

Theorem \ref{thm9.1} refines \cite[Theorem 6.6]{Yan22} by incorporating
amplification and tracking the dependence on the varying automorphic weight
$\phi'$.

\subsection{Preliminary Reduction}
By Cauchy--Schwarz,
\begin{equation}\label{9.1.}
\big|J^{\Reg,\RNum{2}}_{\Geo,\bi}(f,\mathbf{0})\big|\leq \sqrt{\mathcal{J}_{\Geo,\bi}^{\Reg,\RNum{2}}(f)\widetilde{\mathcal{J}}_{\Geo,\bi}^{\Reg,\RNum{2}}(f)},
\end{equation}
where
\begin{align*}
&\mathcal{J}_{\Geo,\bi}^{\Reg,\RNum{2}}(f)
:=\sum_{\substack{(\transp{\boldsymbol{\xi}},t)\in M_{1,n}(\mathbb Q)\\
(\transp{\boldsymbol{\xi}},t)\neq\mathbf 0,\ t\neq1}}
\int_{G'(\mathbb A)}\int_{P_0'(\mathbb Q)\backslash G'(\mathbb A)}
\left|f\left(\iota(x)^{-1}\gamma(\boldsymbol\xi,t)\iota(xy)\right)\right||\phi'(x)|^2\,dx\,dy,\\
&\widetilde{\mathcal{J}}_{\Geo,\bi}^{\Reg,\RNum{2}}(f)
:=\sum_{\substack{(\transp{\boldsymbol{\xi}},t)\in M_{1,n}(\mathbb Q)\\
(\transp{\boldsymbol{\xi}},t)\neq\mathbf 0,\ t\neq1}}
\int_{G'(\mathbb A)}\int_{P_0'(\mathbb Q)\backslash G'(\mathbb A)}
\left|f\left(\iota(x)^{-1}\gamma(\boldsymbol\xi,t)\iota(xy)\right)\right||\phi'(xy)|^2\,dx\,dy.
\end{align*}

\begin{lemma}\label{lem9.2.}
Define $f^{-1}(g):=f(g^{-1})$. Then
\begin{align*}
\widetilde{\mathcal J}_{\Geo,\bi}^{\Reg,\RNum{2}}(f)
=\mathcal J_{\Geo,\bi}^{\Reg,\RNum{2}}(f^{-1}).
\end{align*}
\end{lemma}
\begin{proof}
This is the inversion symmetry proved in \cite[Lemma 6.7]{Yan22}. Briefly,
the changes of variables $x\mapsto xy^{-1}$ and $y\mapsto y^{-1}$ replace
the orbital matrix by its inverse. Since $t\neq1$, rational left and right
translations in $G'$ transform this inverse back into the same family of
matrices, after the substitutions
$x\mapsto\diag(I_{n-1},-1)x$ and
$\boldsymbol{\xi}\mapsto-\boldsymbol{\xi}$. The remaining central scalar
has modulus one under the central character and disappears after taking
absolute values.
\end{proof}

By \eqref{9.1.} and Lemma \ref{lem9.2.}, Theorem \ref{thm9.1} follows from the following proposition.
\begin{prop}\label{prop9.3.}
Set $\mathcal{J}:=\max\big\{\mathcal{J}_{\Geo,\bi}^{\Reg,\RNum{2}}(f), \mathcal{J}_{\Geo,\bi}^{\Reg,\RNum{2}}(f^{-1})\big\}$. Then 
\begin{align*}
\mathcal{J}\ll (TMM')^{\varepsilon}T^{\frac{n-1}{2}}M'^{2n}M^n\mathcal{N}_f^{n-2}\langle\phi',\phi'\rangle\prod_{p\mid M'}p^{ne_p(M)} \cdot \max_{1\leq \ell_0\leq n}\Bigg[\frac{M'^{2}\mathcal{N}_f^{2}}{M^{\dagger}}\Bigg]^{\ell_0},
\end{align*}	
with the same dependence as in Theorem \ref{thm9.1}.
\end{prop}

The remainder of the section proves Proposition \ref{prop9.3.}.

\subsection{Notation and Local Factorization}

\subsubsection{Siegel sets and auxiliary constants}\label{Sie}
By the Siegel-set decomposition of \textsection\ref{sec4.4}, every class in
$[P_0']$ has a representative
$b=ak^*$ with $a\in A_H(\mathbb R)$ and $k^*\in\Omega^*$, where
$\Omega^*\subset P_0'(\mathbb A)$ is fixed and compact. Let $\mathcal P$ be
the finite set consisting of $p_\circ$ and the primes at which
$\Omega_p^*$ is not contained in $G'(\mathbb Z_p)$. At these fixed places,
we enlarge the constants below, if necessary, to absorb the support of the
fixed local test function. Enlarging $L$ if necessary, we may assume that
$\mathcal P\cap\mathcal L=\varnothing$.

Write $\Omega^*=\otimes_{p\leq\infty}\Omega_p^*$. For $p\in\mathcal P$,
choose a fixed integer $\kappa_p\leq0$ controlling the denominators in
$\Omega_p^*$ and in the fixed local test function; set $\kappa_p=0$
otherwise. Define
\begin{equation}\label{equ9.3}
e_{\min}(\Omega_p^*):=\min\big\{\kappa_p,\ e_p(E_{i,j}(x)),\
e_p(E_{i,j}(x^{-1})):
x\in\Omega_p^*,\ 1\leq i,j\leq n\big\}.
\end{equation}
, where $E_{i,j}(x)$ denotes the $(i,j)$-entry of $x$. Thus
$e_{\min}(\Omega_p^*)=O_{\Omega^*}(1)$ and is zero for $p\notin\mathcal P$.
Set
$C_{\Omega_p^*}:=p^{e_{\min}(\Omega_p^*)}\leq1$ and
$C_{\Omega^*}:=\prod_{p<\infty}C_{\Omega_p^*}$. This is a fixed positive
constant depending only on $\Omega^*$.

\subsubsection{Local components}\label{sec9.2.2} 
Let $f\in\big\{f_{p_0,i}^{\diag},f_{p_1,p_2}^{\mathrm{off}}\big\}$. We
analyze $f\bigl(\iota(x)^{-1}\gamma(\boldsymbol\xi,t)\iota(xy)\bigr)$
to determine the support of $x$ and $y$.

Write $x=zbk$ with $z\in\mathbb A^\times$, $b=ak^*$,
$a\in A_H(\mathbb R)$, $k^*\in\Omega^*$, and $k\in K'$. The preceding
expression factors as $\mathfrak F_\infty\mathfrak F_{\fin}$, where
\begin{equation}\label{55}
\mathfrak{F}_{\infty}:=f_{\infty}\left(\iota(z_{\infty}ak_{\infty}^*k_{\infty})^{-1}
\gamma(\boldsymbol\xi,t)\iota(z_{\infty}ak_{\infty}^*k_{\infty}y_{\infty})\right)
\end{equation}
and $\mathfrak{F}_{\fin}:=\prod_{p< \infty}\mathfrak{F}_p,$ with 
\begin{equation}\label{58}
\mathfrak{F}_p:=f_p\left(\iota(z_pk_p^*k_p)^{-1}
\gamma(\boldsymbol\xi,t)\iota(z_pk_p^*k_py_p)\right),\ \ p<\infty. 
\end{equation}

\subsection{Counting Rational Points}\label{sec9.2}
We now determine which rational parameters $(\boldsymbol{\xi},t)$ can
contribute.
\subsubsection{The non-Archimedean constraint}
Define $\mathfrak{X}(f)$ to be the set 
\begin{align*}
\bigg\{(\xi_1,\ldots,\xi_{n-1},t)\in\mathbb Q^n:
t\neq1,\ 
\xi_j\in \frac{C_{\Omega^*}^{2}M^{\dagger}}
{M'^{2}\mathcal{N}_f^{2}}\mathbb{Z},\
\frac{t}{t-1}\in \frac{C_{\Omega^*}^2M^{\dagger}}
{M'^{2}\mathcal{N}_f^{2}}\mathbb{Z},\
1\leq j<n\bigg\}.
\end{align*}

\begin{lemma}\label{lem9.2}
The product $\mathfrak F_{\fin}$ vanishes unless
$(\boldsymbol{\xi},t)\in\mathfrak X(f)$.\end{lemma}
\begin{proof}
Write $\boldsymbol{\xi}=\transp{(\xi_1,\cdots,\xi_{n-1})}\in \mathbb{Q}^{n-1}$. Let $p<\infty,$ write $z_p=p^rI_n$.

%By the construction of $f_{\infty},$ \eqref{55} vanishes unless $\xi_j\ll z_{\infty}a_j,$ $1\leq j\leq n-1,$ and $t\ll z_{\infty}\ll 1,$ and $y_{\infty}\in \mathcal{Y}_{\infty},$ which is a fixed compact neighborhood of $I_n\in G'(\mathbb{R})$. 

\begin{enumerate}
	\item[(A)] Suppose $r\geq 0$. By Cramer's rule there exists $k_{p}^{(1)}\in G(\mathbb{Z}_p)$ such that 
\begin{align*}
\iota(z_p)^{-1}\gamma(\boldsymbol\xi,t)\iota(z_p)=(1-t)k_{p}^{(1)}\begin{pmatrix}
\frac{1}{1-t}I_{n-1}&&p^{-r}\boldsymbol{\xi}\\
&\frac{1}{1-t}&\frac{p^{-r}t}{1-t}\\
&&1
\end{pmatrix}.
\end{align*}
	
\item[(B)] Suppose $r< 0$. By Cramer's rule there exists $k_{p}^{(2)}\in G(\mathbb{Z}_p)$ such that 
\begin{align*}
\iota(z_p)^{-1}\gamma(\boldsymbol\xi,t)\iota(z_p)=p^{-r}(1-t)k_{p}^{(2)}\begin{pmatrix}
\frac{p^r}{1-t}I_{n-1}&&\boldsymbol{\xi}\\
&\frac{p^{2r}}{1-t}&\frac{p^{r}}{1-t}\\
&&1
\end{pmatrix}.
\end{align*}
\end{enumerate}
Here $p^{-r}$ and $(1-t)$ are scalars identified with $p^{-r}I_{n+1}$ and $(1-t)I_{n+1},$ respectively. 

Let $t\in \mathbb{Q}-\{1\}$. Consider the various cases as follows.
%Since $k_{\fin}^*k_{\fin}$ ranges over a compact set, the support of $f_{\fin}$ and the above cases (A) and (B) forces that \eqref{58} vanishes unless $y_{\fin}$ lies in a compact subset of $G'(\mathbb{A}_{\fin})$. 
\begin{enumerate}
\item Let $p\nmid MM'\nu(f)$ and $p\notin\mathcal{P}$. 
\begin{itemize}
\item Suppose $e_p(t-1)\geq 0$. From case (B) and the support of $f_p$, it follows that $r\geq 0$ by analyzing the $(n,n+1)$-th entry in case (B). So 
\begin{align*}
\mathfrak{F}_p=f_p\left(
\begin{pmatrix}
\frac{1}{1-t}I_{n-1}&&p^{-r}\boldsymbol{\xi}\\
&\frac{1}{1-t}&\frac{p^{-r}t}{1-t}\\
			&&1
\end{pmatrix}\iota(k_p^*k_py_p)
	\right)\neq 0
\end{align*}
if and only if $r\geq 0$ and there exists some $\lambda_p\in \mathbb{Q}_p^{\times}$ such that 
\begin{align*}
\lambda_p\begin{pmatrix}
\frac{1}{1-t}I_{n-1}&&p^{-r}\boldsymbol{\xi}\\
&\frac{1}{1-t}&\frac{p^{-r}t}{1-t}\\
&&1
\end{pmatrix}\iota(k_p^*k_py_p)\in G(\mathbb{Z}_p),
\end{align*}
which forces that $\lambda_p\in \mathbb{Z}_p^{\times}$ and 
\begin{align*}
\begin{cases}
e_p(t)-e_p(t-1)\geq r\geq 0\\
e_p(\xi_j)\geq r\geq 0,\ 1\leq j\leq n-1.
\end{cases}
\end{align*}
					
\item Suppose $e_p(t-1)\leq -1$. Then $e_p(t)=e_p(t-1)$. From the case (B) we obtain $e_p(t-1)\leq r\leq -1;$ in the case (A) we have $0\leq r\leq e_p(t)-e_p(t-1),$ implying $r=0$. So 
\begin{equation}\label{9.4..}
\mathfrak{F}_p=\begin{cases}
			f_p\left(
\begin{pmatrix}
\frac{1}{1-t}I_{n-1}&&\boldsymbol{\xi}\\
&\frac{1}{1-t}&\frac{t}{1-t}\\
			&&1
\end{pmatrix}\iota(k_p^*k_py_p)
	\right),&\text{if $r=0$}\\
	f_p\left(
\begin{pmatrix}
\frac{p^r}{1-t}I_{n-1}&&\boldsymbol{\xi}\\
&\frac{p^{2r}}{1-t}&\frac{p^{r}}{1-t}\\
			&&1
\end{pmatrix}\iota(k_p^*k_py_p)
	\right),& \text{if $r<0$}.
\end{cases}
\end{equation}
		
As a consequence, we have
\begin{itemize}
\item If $r=0$, then $\mathfrak{F}_p$ can be nonzero only if
$\boldsymbol{\xi}\in \mathbb{Z}_p^{n-1}$ and
$e_p(t)-e_p(t-1)\geq 0$.
\item If $r<0$, then $\mathfrak{F}_p$ can be nonzero only if
$\boldsymbol{\xi}\in \mathbb{Z}_p^{n-1}$ and $r-e_p(t-1)\geq0$; hence
$e_p(t-1)\leq r<0$.
\end{itemize}
\end{itemize}
	
Thus $\mathfrak{F}_p$ can be nonzero only if
\begin{align*}
\begin{cases}
e_p(t)-e_p(t-1)\geq 0\\
e_p(\xi_j)\geq  0,\ 1\leq j\leq n-1.
\end{cases}
\end{align*} 
	
\item Let $p\mid M'$. By \eqref{3.13},
$\supp f_p=Z(\mathbb Q_p)D_p$, where
\begin{align*} 
D_p:=\bigcup_{\substack{\alpha_j,\ \beta_j\in \mathbb{Z}/p^{m'}\mathbb{Z}\\ 1\leq j<n}}\bigcup_{\substack{\alpha_n \in (\mathbb{Z}/p^{m''}\mathbb{Z})^{\times}\\ \beta_n \in (\mathbb{Z}/p^{m''}\mathbb{Z})^{\times}}}\begin{pmatrix}
I_n&\mathfrak{b}_{\boldsymbol{\alpha}}\\
&1
\end{pmatrix}K_p\begin{pmatrix}
I_n&\mathfrak{b}_{\boldsymbol{\beta}}\\
&1
\end{pmatrix}.
\end{align*} 
Here $\mathfrak{b}_{\boldsymbol{\alpha}}=\transp{(\alpha_1p^{-m'},\cdots, \alpha_{n-1}p^{-m'},\alpha_n p^{-m''})},$ and $\mathfrak{b}_{\boldsymbol{\beta}}$ is defined similarly.

By the case (A) and case (B), there exists some $\lambda_p\in \mathbb{Q}_p^{\times}$ such that 
\begin{equation}\label{eq9.4}
\lambda_p\begin{pmatrix}
\frac{1}{1-t}I_{n-1}&&p^{-r}\boldsymbol{\xi}\\
			&\frac{1}{1-t}&\frac{p^{-r}t}{1-t}\\
			&&1
\end{pmatrix}\iota(k_p^*k_py_p)\in \iota(k_p^*k_p)D_p
\end{equation}
if $r\geq 0;$ and if $r<0,$ we have 
\begin{equation}\label{eq9.5}
\lambda_p\begin{pmatrix}
\frac{p^r}{1-t}I_{n-1}&&\boldsymbol{\xi}\\
&\frac{p^{2r}}{1-t}&\frac{p^{r}}{1-t}\\
			&&1
\end{pmatrix}\iota(k_p^*k_py_p)\in \iota(k_p^*k_p)D_p.
\end{equation}
	
Recall that  $m'=e_p(M')$. The constraint \eqref{eq9.4}  implies that 
\begin{equation}\label{eq9.9}
\begin{cases}
-m'+e_{\min}(\Omega_p^*)\leq e_p(\lambda_p)\leq m'-e_{\min}(\Omega_p^*),\ \ r\geq 0,\\
e_p(t)-e_p(t-1)-r+e_p(\lambda_p)\geq -m'+e_{\min}(\Omega_p^*)\\
e_p(\xi_j)-r+e_p(\lambda_p)\geq  -m'+e_{\min}(\Omega_p^*),\ 1\leq j\leq n-1,
\end{cases}
\end{equation}
and \eqref{eq9.5} implies that 
\begin{equation}\label{eq9.11}
\begin{cases}
-m'+e_{\min}(\Omega_p^*)\leq e_p(\lambda_p)\leq m'-e_{\min}(\Omega_p^*),\ \ r<0,\\
r-e_p(t-1)+e_p(\lambda_p)\geq -m'+e_{\min}(\Omega_p^*)\\
e_p(\xi_j)+e_p(\lambda_p)\geq  -m'+e_{\min}(\Omega_p^*),\ 1\leq j\leq n-1,
\end{cases}
\end{equation}
where $e_{\min}(\Omega_p^*)$ was defined by \eqref{equ9.3}. The upper bound
for $e_p(\lambda_p)$ follows by inverting \eqref{eq9.4} and \eqref{eq9.5} and
examining the $(n+1,n+1)$-entry:
\begin{align*}
e_p(\lambda_p)\leq m'-e_{\min}(\Omega_p^*)
\end{align*}
Consequently, \eqref{eq9.9} and \eqref{eq9.11} imply
\begin{align*}
\begin{cases}
e_p(t)-e_p(t-1)\geq -2m'+2e_{\min}(\Omega_p^*)\\
e_p(\xi_j)\geq  -2m'+2e_{\min}(\Omega_p^*),\ 1\leq j\leq n-1.
\end{cases}
\end{align*} 

\item Let $p\in\mathcal P$ and $p\nmid M'$. Since both
$\Omega_p^*$ and $\supp f_p$ are fixed and compact modulo the center, the
same calculation gives fixed lower bounds for
$e_p(t)-e_p(t-1)$ and $e_p(\xi_j)$. These bounds are absorbed into
$C_{\Omega^*}$.

\item Let $p\mid\nu(f)$. Then $p\notin\mathcal{P},$ namely, $\Omega_p^*=K_p$. Recall definitions in \textsection\ref{3.6.2},  $\supp f_p=Z(\mathbb{Q}_p)D_p',$ where 
\begin{align*}
D_p':=\begin{cases}
K_p\diag(p^{\ell_p},I_{n})K_p,\ \ & \text{if $f=f_{p_1,p_2}^{\mathrm{off}}$ and $p=p_1,$}\\
K_p\diag(I_{n},p^{-\ell_p})K_p,\ \ & \text{if $f=f_{p_1,p_2}^{\mathrm{off}}$ and $p=p_2,$}\\
K_p\diag(p^{i},I_{n-1},p^{-i})K_p,\ \ & \text{if $f=f_{p_0,i}^{\diag}$ and $p=p_0$.}
\end{cases}
\end{align*}

By the case (A) and case (B), there exists some $\lambda_p\in \mathbb{Q}_p^{\times}$ such that 
\begin{align*}
\lambda_p\begin{pmatrix}
\frac{1}{1-t}I_{n-1}&&p^{-r}\boldsymbol{\xi}\\
			&\frac{1}{1-t}&\frac{p^{-r}t}{1-t}\\
			&&1
\end{pmatrix}\iota(k_p^*k_py_p)\in \iota(k_p^*k_p)D_p'=D_p'
\end{align*}
if $r\geq 0;$ and if $r<0,$ we have 
\begin{align*}
\lambda_p\begin{pmatrix}
\frac{p^r}{1-t}I_{n-1}&&\boldsymbol{\xi}\\
&\frac{p^{2r}}{1-t}&\frac{p^{r}}{1-t}\\
			&&1
\end{pmatrix}\iota(k_p^*k_py_p)\in \iota(k_p^*k_p)D_p'=D_p'.
\end{align*}

The classifications in Lemmas \ref{lem7} and \ref{lem7.} show that
$\mathfrak{F}_p\neq0$ only if
\begin{equation}\label{9.6.3}
\begin{cases}
r\geq 0,\
e_p(t)-e_p(t-1)-r\geq -\delta_p\\
e_p(\xi_j)-r\geq  -\delta_p,\ 1\leq j\leq n-1,
\end{cases}
\end{equation} 
or 
\begin{equation}\label{9.6.4}
\begin{cases}
r< 0,\
r-e_p(t-1)\geq -\delta_p\\
e_p(\xi_j)\geq  -\delta_p,\ 1\leq j\leq n-1,
\end{cases}
\end{equation} 
where $\delta_p=\ell_p$ if $f=f_{p_1,p_2}^{\mathrm{off}}$ and $p\mid p_1p_2,$ and $\delta_p=2i$ if $f=f_{p_0,i}^{\diag}$ and $p=p_0$. Both \eqref{9.6.3} and \eqref{9.6.4} yield that 
\begin{align*}
\begin{cases}
e_p(t)-e_p(t-1)\geq -\delta_p\\
e_p(\xi_j)\geq  -\delta_p,\ 1\leq j\leq n-1.
\end{cases}
\end{align*} 
\item Finally we consider the case that $p\mid M,$ $p\nmid M'$ and $p\notin \mathcal{P}$. By Iwasawa decomposition, we can write $k_p^*k_py_p=p^{r'}a_p'u_p'k_p'$ for some $r'\in\mathbb{Z},$ $a_p'\in A'(\mathbb{Q}_p),$ $u_p'\in N'(\mathbb{Q}_p)$ and $k_p'\in K_p'$. Then by definition in \textsection\ref{11.1.3},  
\begin{align*}
\mathfrak{F}_p=f_p\left(\iota(z_pk_p^*k_p)^{-1}
\gamma(\boldsymbol\xi,t)\iota(z_pk_p^*k_py_p)\right)\neq 0
\end{align*}   
unless 
\begin{equation}\label{eq9.10}
\lambda_p\iota(z_pk_p^*k_p)^{-1}\gamma(\boldsymbol\xi,t)
\iota(z_pk_p^*k_py_p)\in K_p(M)
\end{equation}
for some $\lambda_p\in \mathbb{Q}_p^{\times}$. Observe that \eqref{eq9.10} amounts to 
\begin{align*}
\lambda_p\begin{pmatrix}
*&*&p^{-r}(1-t)\boldsymbol{\xi}\\
&p^{r'}&p^{-r}t\\
&p^{r+r'}&1
\end{pmatrix}\in K_p(M),
\end{align*}
which leads to that 
\begin{equation}\label{cons}
\begin{cases}
e_p(\lambda_p)+r'=0,\ e_p(\lambda_p)\geq 0,\\
e_p(\lambda_p)+r+r'\geq e_p(M),\\
e_p(\lambda_p)-r+e_p(t)\geq 0,\\
e_p(\lambda_p)-r+e_p(1-t)+e_p(\xi_j)\geq 0,\ 1\leq j\leq n-1,\\
2e_p(\lambda_p)+r'+e_p(1-t)=0.
\end{cases}
\end{equation}
Notice that the last constraint in \eqref{cons} comes from the determinant of the lower right $2\times 2$-corner. So $r\geq e_p(M),$ $e_p(1-t)=-e_p(\lambda_p)\leq 0,$ and $e_p(\xi_j)\geq r\geq e_p(M),$ $1\leq j\leq n-1$. 
		\end{enumerate}

Combining these cases proves the lemma.
\end{proof}

\subsubsection{The Archimedean constraint}

\begin{lemma}\label{lem9.1}
We have $\mathfrak F_\infty=0$ unless
\begin{align*}
|\xi_j|\ll |z_\infty|a_j,\quad 1\leq j<n,
\ \ \text{and}\ \  |t|\ll |z_\infty|.
\end{align*}
The implied constants depend only on $c_\infty$, $C_\infty$, and
$\Omega_\infty^*$.
\end{lemma}
\begin{proof}
Recall the support of $f_{\infty}$ (see \eqref{245} in \textsection\ref{3.2.1}). We have $\mathfrak{F}_{\infty}=0$ unless	
\begin{equation}\label{eq9.13}
\iota(z_{\infty}ak_{\infty}^*k_{\infty})^{-1}\gamma(\boldsymbol\xi,t)
\iota(z_{\infty}ak_{\infty}^*k_{\infty}y_{\infty})
\in I_{n+1}+O(T^{-\varepsilon}).
\end{equation}

The upper-left $n\times n$ block gives
$y_\infty=I_n+O(T^{-\varepsilon})$, while determinants give
$|1-t|=1+O(T^{-\varepsilon})$.

Consider the $(j,n+1)$-th entry of \eqref{eq9.13}. We obtain 
\begin{align*}
\gamma(\boldsymbol\xi,t)\in \iota(z_{\infty}ak_{\infty}^*k_{\infty})(I_{n+1}+O(T^{-\varepsilon})),
\end{align*}
which implies
\begin{align*}
\begin{cases}
|(1-t)\xi_j|\ll |z_{\infty}|a_j,\ \ 1\leq j\leq n-1,\\
|t|\ll |z_{\infty}|.
\end{cases}
\end{align*}
Since $|1-t|\asymp1$, the result follows.
\end{proof}

\subsubsection{The global constraint}
Let $z_\infty\in\mathbb R^\times$ and $a\in A_H(\mathbb R)$. Define 
\begin{align*}
\mathfrak X(f;a,z_\infty):=
\big\{(\boldsymbol{\xi},t)\in\mathfrak X(f):
|\xi_j|\ll|z_\infty|a_j,\ 1\leq j<n,\
|t|\ll|z_\infty|\ll1\big\}.
\end{align*}

\begin{cor}\label{cor9.4}
Uniformly for
$z_{\fin}\in\mathbb A_{\fin}^{\times}$, $k^*\in\Omega^*$, $k\in K'$, and
$y\in G'(\mathbb A)$, we have 
\begin{align}\label{9.6}
f\left(\iota(z_{\infty}z_{\fin}ak^*k)^{-1}\gamma(\boldsymbol\xi,t)
\iota(z_{\infty}z_{\fin}ak^*ky)\right)=0
\end{align}
unless $(\boldsymbol{\xi},t)\in \mathfrak{X}(f;a,z_{\infty})$. The implied constants are those of Lemma \ref{lem9.1}.
\end{cor}

\begin{defn}
Define
\begin{align*}
\mathfrak{X}^*(f;a,z_{\infty}):=\Big\{(\boldsymbol{\xi},t)\in \mathfrak{X}(f;a,z_{\infty}):\ (\boldsymbol{\xi},t)\neq \mathbf{0}\Big\}. 
\end{align*}
\end{defn}

Counting lattice points in each coordinate gives
\begin{equation}\label{9.7}
\big|\mathfrak{X}^*(f;a,z_{\infty})\big|\ll \max\Bigg\{1,\frac{M'^2\mathcal{N}_f^2|z_{\infty}|}{C_{\Omega^*}^2M^{\dagger}}\Bigg\}\cdot \prod_{j=1}^{n-1}\max\Bigg\{1,\frac{M'^{2}\mathcal{N}_f^2|z_{\infty}|a_j}{C_{\Omega^*}^2M^{\dagger}}\Bigg\}, 
\end{equation}
where the implied constant depends only on $c_\infty$, $C_\infty$, and
$\Omega_\infty^*$.

For $(\boldsymbol\xi,t)\in\mathfrak X^*(f;a,z_\infty)$, either $t\neq0$
or some $\xi_j\neq0$. Since $a_1\geq\cdots\geq a_{n-1}$, the defining
lattice constraints imply
\begin{equation}\label{10.17}
C_{\Omega^*}^2M^{\dagger}M'^{-2}\mathcal N_f^{-2}
\left(\textbf{1}_{t\neq0}+a_1^{-1}\textbf{1}_{t=0}\right)
\ll |z_\infty|\ll1.
\end{equation}
When $n=1$, the term involving $a_1$ is omitted.

\subsection{Reduction to the Essential Range}\label{sec9.3}
With the above parametrization, 
\begin{align*}
\mathcal{J}_{\Geo,\bi}^{\Reg,\RNum{2}}(f)=\int_{P_0'(\mathbb Q)\backslash G'(\mathbb A)}\int_{{G'}(\mathbb{A})}
\sum_{(\boldsymbol{\xi},t)\in \mathfrak{X}^*(f;a,z_{\infty})}
\Big|f\bigl(\iota(x)^{-1}\gamma(\boldsymbol\xi,t)\iota(xy)\bigr)\Big|
\big|\phi'(x)\big|^2dydx.
\end{align*}
The sum is finite by \eqref{9.7}.

\begin{defn}\label{defn9.9}
Apply Lemma \ref{lem4.9} with $c_0=n+11$, and let $c_1$ be the resulting
exponent. Set $c:=2c_1+10n^2$ and define
\begin{align*}
\mathfrak A_{c,\varepsilon}:=
\big\{a\in A_H(\mathbb R):
\|a\|\leq T^c,\ 
T^{-\varepsilon}\ll|\det a|_\infty\ll T^\varepsilon\big\},
\end{align*}
where
$\|a\|:=\prod_{j=1}^{n-1}\max\{a_j,a_j^{-1}\}$.
\end{defn}

\begin{lemma}\label{lem9.7.}
For $a\in A_H(\mathbb R)\setminus\mathfrak A_{c,\varepsilon}$ and
$k^*\in\Omega^*$, $k\in K'$, one has
\begin{equation}\label{9.10.}
|\phi'(ak^*k)|\cdot \|\phi'\|_{\infty}\ll T^{-10n^2}\|a\|^{-n-10},
\end{equation}	
where the implied constant depends only on the fixed data.
\end{lemma}
\begin{proof}
Lemma \ref{lem4.9}, applied at the identity, gives
$\|\phi'\|_\infty\ll T^{c_1}$. We distinguish three cases.
\begin{enumerate}
\item If $\|a\|>T^c$, use Lemma \ref{lem4.9} with $c_2=0$:
\begin{align*}
|\phi'(ak^*k)|\|\phi'\|_\infty
\ll \|a\|^{-n-11}T^{2c_1}
\ll T^{-10n^2}\|a\|^{-n-10}.
\end{align*}
\item If $\|a\|\leq T^c$ and $|\det a|_\infty\gg T^\varepsilon$, use $c_2=-(10n^2+2c_1+(n+10)c)\varepsilon^{-1}$. 
Then Lemma \ref{lem4.9} gives
\begin{align*}
|\phi'(ak^*k)|\|\phi'\|_\infty
\ll T^{-10n^2-(n+10)c}
\ll T^{-10n^2}\|a\|^{-n-10}.
\end{align*}
\item The case $\|a\|\leq T^c$ and
$|\det a|_\infty\ll T^{-\varepsilon}$ is identical, with the opposite
choice of sign for $c_2$.
\end{enumerate}
\end{proof}

Write $x=zbk$ and $b=ak^*=ak_{\infty}^*k_{\fin}^*$ as before. By \eqref{87} and Corollary \ref{cor9.4}, 
\begin{align*}
\mathcal{J}_{\Geo,\bi}^{\Reg,\RNum{2}}(f)\ll \mathcal{J}_{\Geo,\bi}^{\Reg,\RNum{2},\dagger}(f)+\mathcal{J}_{\Geo,\bi}^{\Reg,\RNum{2},\heartsuit}(f),
\end{align*}
 where  
\begin{multline*}
\mathcal{J}_{\Geo,\bi}^{\Reg,\RNum{2},\dagger}(f):= 
\int_{Z'(\mathbb{A})}\int_{K'}\int_{G'(\mathbb{A})}
\int_{\Omega^*}\int_{A_H(\mathbb{R})-\mathfrak{A}_{c,\varepsilon}}\big|\phi'(ak^*k)\big|^2\\
\sum_{\substack{(\boldsymbol{\xi},t)\in \mathfrak{X}^*(f;a,z_{\infty})}}|\mathfrak F_\infty\mathfrak{F}_{\fin}|\cdot \delta_{B'}^{-1}(a)d^{\times}a\,dk^*\,dy\,dk\,d^{\times}z,
\end{multline*}
and $\mathcal{J}_{\Geo,\bi}^{\Reg,\RNum{2},\heartsuit}(f)$ is defined similarly except that $a$ ranges through $\mathfrak{A}_{c,\varepsilon}$.

We show that the complementary-range contribution
$\mathcal{J}_{\Geo,\bi}^{\Reg,\RNum{2},\dagger}(f)$ is negligible, reducing
the problem to the essential range.

\subsubsection{Estimate of non-Archimedean local integrals}\label{9.2.1}
Define  
\begin{align*}
\mathcal{J}_p^{\heartsuit}(\boldsymbol{\xi},t;k_{p}^*k_p):=
\int_{Z'(\mathbb{Q}_p)}\int_{G'(\mathbb{Q}_{p})}
|\mathfrak{F}_p|dy_{p}d^{\times}z_p,\quad p<\infty.
\end{align*}

\begin{lemma}\label{lem9.7}
Let notation be as before. Then 
\begin{align*}
	\prod_{p<\infty}\mathcal{J}_p^{\heartsuit}(\boldsymbol{\xi},t;k_{p}^*k_p)\ll M'^{2n+\varepsilon}M^n \mathcal{N}_f^{n-2}\prod_{p\mid M'}p^{ne_p(M)}\cdot \prod_{p<\infty}E_p(\boldsymbol{\xi},t),
\end{align*}
where $E_p(\boldsymbol{\xi},t):=\big|\min\{e_p(t)-e_p(t-1), e_p(\xi_j):\ 1\leq j\leq n-1\}\big|+\big|e_p(t-1)\big|+2e_p(M')-2e_{\min}(\Omega_p^*)+1,$ and the implied constant depends only on $\varepsilon$.
\end{lemma}
\begin{proof}
Consider various finite places as follows. We will make use of case (A) and case (B) in the proof of Lemma \ref{lem9.2}.
\begin{enumerate}
\item Let $p\nmid MM'\nu(f)$ and $p\notin\mathcal{P}$. Then $f_p$ is bi-invariant under $\Omega^*_pK_p'=K_p'$.
\begin{itemize}
\item Suppose $e_p(t-1)\geq 0$. Then from the case (B) and the support of $f_p$ one must have $r\geq 0$. Change variable $y_p\mapsto (k_p^*k_p)^{-1}y_p$ to obtain  
\begin{align*}
\mathcal{J}_p^{\heartsuit}(\boldsymbol{\xi},t;k_{p}^*k_p)=\sum_{r\geq 0}\int_{G'(\mathbb{Q}_{p})}f_p\left(\begin{pmatrix}
I_{n-1}&&p^{-r}(1-t)\boldsymbol{\xi}\\
&1&p^{-r}t\\
&&1-t
\end{pmatrix}\iota(y_p)\right)dy_{p}.
\end{align*}
By definition, $f_p\left(\begin{pmatrix}
I_{n-1}&&p^{-r}(1-t)\boldsymbol{\xi}\\
			&1&p^{-r}t\\
			&&1-t
\end{pmatrix}\iota(y_p)\right)\neq 0$ if and only if 
\begin{align*}
 y_p\in (1-t)K_p'\ \ \text{and}\ \ 0\leq r\leq \min\{e_p(t)-e_p(t-1),\ e_p(\xi_j):\  1\leq j\leq n-1\}.
\end{align*}
		
Therefore, we have 
\begin{equation}\label{56}
\mathcal{J}_p^{\heartsuit}(\boldsymbol{\xi},t;k_{p}^*k_p)=\textbf{1}_{e_p(\xi_j)\geq e_p(t-1)\geq 0,\ 1\leq j\leq n-1}\sum_{\substack{0\leq r\leq e_p(t)-e_p(t-1)\\
r\leq e_p(\xi_j),\ 1\leq j\leq n-1}}1.
\end{equation}
		
\item Suppose $e_p(t-1)\leq -1$. Then $e_p(t)=e_p(t-1)$. From the case (B) we obtain $e_p(t-1)\leq r\leq -1;$ in the case (A) we have $0\leq r\leq e_p(t)-e_p(t-1),$ implying $r=0$. Changing variable $y_p\mapsto (k_p^*k_p)^{-1}y_p,$ 
\begin{align*}
\mathcal{J}_p^{\heartsuit}(\boldsymbol{\xi},t;k_{p}^*k_p)=&\int_{G'(\mathbb{Q}_{p})}f_p\left(\begin{pmatrix}
I_{n-1}&&(1-t)\boldsymbol{\xi}\\
&1&t\\
&&1-t
\end{pmatrix}\iota(y_p)\right)dy_{p}\\
&+\sum_{e_p(t-1)\leq r\leq -1}\int_{G'(\mathbb{Q}_{p})}f_p\left(\begin{pmatrix}
I_{n-1}&&p^{-r}(1-t)\boldsymbol{\xi}\\
&p^{r}&1\\
&&p^{-r}(1-t)
\end{pmatrix}\iota(y_p)\right)dy_{p}.
\end{align*}
		
The first integrand can be nonzero only when $y_p\in(1-t)K_p'$, and the
second only when
$y_p\in(1-t)\diag(p^{-r}I_{n-1},p^{-2r})K_p'$. Hence
\begin{equation}\label{57}
\mathcal{J}_p^{\heartsuit}(\boldsymbol{\xi},t;k_{p}^*k_p)=(1-e_p(t-1))\textbf{1}_{e_p(\xi_j)\geq 0,\ 1\leq j\leq n-1}\textbf{1}_{e_p(t-1)\leq -1}.
\end{equation}
\end{itemize}
	
From the above discussions we see that $\mathcal{J}_p^{\heartsuit}(\boldsymbol{\xi},t;k_{p}^*k_p)=1$ if $e_p(t)=e_p(t-1)=0$ and $e_p(\xi_j)\geq 0,$ $1\leq j\leq n-1;$ and $y_p$ ranges over a compact set depending on $t$. As a consequence, for fixed $\boldsymbol{\xi}$ and $t,$  $\mathcal{J}_p^{\heartsuit}(\boldsymbol{\xi},t;k_{p}^*k_p)=1$ for all but finitely many places $p$.
		
\item Let $p\mid M'$ and put $m'=e_p(M')$. Let
$D_p'\subset G'(\mathbb Q_p)$ be defined by
$\iota(D_p')=\iota(G'(\mathbb Q_p))\cap D_p$. The support conditions
\eqref{eq9.4}--\eqref{eq9.11} show that
$|e_p(\lambda_p)|\leq m'-e_{\min}(\Omega_p^*)$. They also restrict $r$ to
\begin{align*}
0\leq r\leq
\min\{e_p(t)-e_p(t-1),e_p(\xi_j):1\leq j\leq n-1\}
+2m'-2e_{\min}(\Omega_p^*)
\end{align*}
when $r\geq0$, and to
\begin{align*}
e_p(t-1)-2m'+2e_{\min}(\Omega_p^*)\leq r<0
\end{align*}
otherwise. In the two cases, $y_p$ ranges over a left translate of $D_p'$ and a left
translate of $(k_p^*k_p)^{-1}\diag(I_{n-1},p^{-r})(k_p^*k_p)D_p'$, respectively. By the invariance of Haar measure and $\|f_p\|_{\infty}\ll p^{2ne_p(M)-m''}$, we obtain
\begin{align*}
\mathcal{J}_p^{\heartsuit}(\boldsymbol{\xi},t;k_{p}^*k_p)\ll (m'-e_{\min}(\Omega_p^*)+1)E_p(\boldsymbol{\xi},t)\Vol(D_p')\ll  p^{2nm'+2ne_p(M)+\varepsilon}E_p(\boldsymbol{\xi},t).	
\end{align*}

\item Let $p\in\mathcal P$ and $p\nmid M'$. Since $p$ belongs to a fixed
finite set and both $\Omega_p^*$ and $f_p$ are fixed, compactness modulo the
center gives
\begin{align*}
\mathcal J_p^{\heartsuit}(\boldsymbol\xi,t;k_p^*k_p)
\ll E_p(\boldsymbol\xi,t),
\end{align*}
with a constant depending only on the fixed local data.
	
\item Let $p\mid M^{\dagger}:=\prod_{p\mid M,\ p\nmid M'}p^{e_p(M)}$ and $p\notin \mathcal{P}$. By a similar argument,
\begin{align*}
\mathcal{J}_p^{\heartsuit}(\boldsymbol{\xi},t;k_{p}^*k_p)\ll p^{ne_p(M^{\dagger})}\cdot E_p(\boldsymbol{\xi},t).
\end{align*}

\item Let $p\mid\nu(f)$. Assume $\mathfrak{F}_p\neq 0$. Let 
\begin{align*}
M_{n,1}(\mathbb{Q}_p)\ni  \mathfrak{u}:=\begin{cases}
\transp{(\xi_1,\cdots, \xi_{n-1},\frac{p^{-r}t}{1-t})},\ &\text{if $r\geq 0$,}\\
\transp{(\xi_1,\cdots, \xi_{n-1},\frac{p^{r}}{1-t})},\ &\text{if $r<0$.}
\end{cases}
\end{align*}

Let $\mathcal{D}_p(\mathfrak{u})$ be the set of $y_p$ such that 	$f_p\left(\begin{pmatrix}
		y_p&\mathfrak{u}\\
		&1
	\end{pmatrix}\right)\neq 0$. The structure of $\mathcal{D}_p(\mathfrak{u})$ has been classified by Lemmas \ref{lem7} and \ref{lem7.}. Set 
\begin{align*}
\delta_p:=\begin{cases}
\ell_p, &\text{if $f=f_{p_1,p_2}^{\mathrm{off}}$ and $p\mid p_1p_2,$}\\ 
2i, &\text{if $f=f_{p_0,i}^{\diag}$ and $p=p_0$.}
\end{cases}
\end{align*}	
	
By \eqref{9.6.3} and \eqref{9.6.4}, $\mathfrak{F}_p\neq 0$ implies that   	 
\begin{align*}
\begin{cases}
r\geq \min\{e_p(t-1),0\}-\delta_p,\\
r\leq \min\{e_p(t)-e_p(t-1), e_p(\xi_j):\ 1\leq j\leq n-1\}+\delta_p,\\
y_p\in (1-t)\mathcal{D}_p(\mathfrak{u}),\ \text{if $r\geq 0,$}\\
y_p\in (1-t)\diag(p^{-r}I_{n-1},p^{-2r})\mathcal{D}_p(\mathfrak{u}),\ \text{if $r<0$.}
\end{cases}
\end{align*}

Therefore,
\begin{align*}
\frac{\mathcal{J}_p^{\heartsuit}(\boldsymbol{\xi},t;k_{p}^*k_p)}{\|f_p\|_\infty}\ll \sum_{0\leq r\leq d_3}\int_{(1-t)\mathcal{D}_p(\mathfrak{u})}dy_p+\sum_{d_4\leq r< 0}\int_{(1-t)\diag(p^{-r}I_{n-1},p^{-2r})\mathcal{D}_p(\mathfrak{u})}dy_p,
\end{align*}
where $d_3:=\min\{e_p(t)-e_p(t-1), e_p(\xi_j):\ 1\leq j\leq n-1\}+\delta_p$ and $d_4:=\min\{e_p(t-1),0\}-\delta_p$.
Since $\delta_p\leq 2(n+1)$, it follows that
\begin{align*}
\mathcal{J}_p^{\heartsuit}(\boldsymbol{\xi},t;k_{p}^*k_p)\ll \|f_p\|_\infty E_p(\boldsymbol{\xi},t)\cdot \max_{\mathfrak{u}}\Vol(\mathcal{D}_p(\mathfrak{u})).
\end{align*}
	
Lemmas \ref{lem7} and \ref{lem7.} give the required volume bounds for
$\mathcal{D}_p(\mathfrak{u})$.
	
\begin{itemize}
\item If $\nu(f)=p_1p_2$ and $p\mid\nu(f)$, cases (A) and (B) of
Lemma \ref{lem7} give
\begin{align*}
\mathcal J_p^{\heartsuit}(\boldsymbol\xi,t;k_p^*k_p)
\ll p^{-n\ell_p/2}p^{(n-1)\ell_p}E_p(\boldsymbol\xi,t)
=p^{(n/2-1)\ell_p}E_p(\boldsymbol\xi,t),
\end{align*}
where $\|f_p\|_\infty\ll p^{-n\ell_p/2}$.
\item If $\nu(f)=p_0$ and $p\mid\nu(f)$, the cases in Lemma \ref{lem7.}
give
\begin{align*}
\mathcal J_p^{\heartsuit}(\boldsymbol\xi,t;k_p^*k_p)
\ll p^{-ni}p^{2(n-1)i}E_p(\boldsymbol\xi,t)
=p^{(n-2)i}E_p(\boldsymbol\xi,t),
\end{align*}
where $\|f_p\|_\infty\ll p^{-ni}$. In case (A.2), for example,
\begin{align*}
\max_{\mathfrak{u}}\Vol(\mathcal{D}_p(\mathfrak{u}))\ll p^{(n-2)(r_1-i)+(n-1)(e-i)}\ll p^{2(n-1)i}.
\end{align*}
Here $r_1$ and $e$ are defined in Lemma \ref{lem7.}, case (A.2). 
\end{itemize}
These estimates use the relation between $\mathfrak{u}$ and the unipotent
coordinate $\mathfrak{u}''$ of $y_p$ in Lemmas \ref{lem7} and \ref{lem7.}.
In either the diagonal or the off-diagonal case, their product over
$p\mid\nu(f)$ is $\ll\mathcal N_f^{n-2}\prod_{p\mid\nu(f)}E_p(\boldsymbol\xi,t)$.
	
\end{enumerate}

 Combining the local estimates proves the lemma.
 \end{proof}

\subsubsection{Separating the automorphic weight}
By Lemma \ref{lem4.6}, one has
$\langle\phi',\phi'\rangle\gg_{\pi'}1$. In conjunction with Lemma
\ref{lem9.7.}, this gives, for $a\notin \mathfrak{A}_{c,\varepsilon}$,
\begin{equation}\label{9.13}
\big|\phi'(ak^*k)\big|^2\leq |\phi'(ak^*k)|\cdot \|\phi'\|_{\infty}\ll T^{-10n^2}\|a\|^{-n-10}\ll \frac{\langle\phi',\phi'\rangle}{T^{9n^2}\|a\|^{n+10}}.
\end{equation}

This estimate separates the initially nonfactorizable automorphic weight from
the remaining local integrals.

\subsubsection{Bounding $\mathcal{J}_{\Geo,\bi}^{\Reg,\RNum{2},\dagger}(f)$}
\begin{prop}\label{prop9.11.}
 Let notation be as before. Then
 \begin{align*}
 \frac{\mathcal{J}_{\Geo,\bi}^{\Reg,\RNum{2},\dagger}(f)}{\langle\phi',\phi'\rangle}\ll (TMM')^{\varepsilon}T^{-8n^2}M'^{2n}M^n\mathcal{N}_f^{n-2}\prod_{p\mid M'}p^{ne_p(M)}\max_{1\leq \ell_0\leq n}\Bigg[\frac{M'^{2}\mathcal{N}_f^{2}}{M^{\dagger}}\Bigg]^{\ell_0},
 \end{align*}
with the same dependence as in Theorem \ref{thm9.1}.
 \end{prop}
\begin{proof}
Substituting \eqref{9.13} into the definition of
$\mathcal{J}_{\Geo,\bi}^{\Reg,\RNum{2},\dagger}(f)$ gives
\begin{multline*}
\mathcal{J}_{\Geo,\bi}^{\Reg,\RNum{2},\dagger}(f)\ll T^{-9n^2}\langle\phi',\phi'\rangle\int_{K'}
	\int_{\Omega^*}\int_{A_H(\mathbb{R})-\mathfrak{A}_{c,\varepsilon}}\|a\|^{-n-10}
	\int_{Z'(\mathbb{R})}\int_{G'(\mathbb{R})}\\
\sum_{\substack{(\boldsymbol{\xi},t)\in \mathfrak{X}^*(f;a,z_{\infty})}}\prod_{p<\infty}\mathcal{J}_p^{\heartsuit}(\boldsymbol{\xi},t;k_{p}^*k_p)|\mathfrak F_\infty|\,dy_{\infty}d^{\times}z_{\infty}
\frac{d^{\times}a}{\delta_{B'}(a)}dk^*dk.
\end{multline*}

By Lemma \ref{lem9.7} and the fact that $\delta_{B'}(a)\geq \|a\|^{-n},$ the above integral is further bounded by 
\begin{multline}\label{e9.23}
C^{\dag}\int_{K'_{\infty}}
	\int_{\Omega^*_{\infty}}\int_{A_H(\mathbb{R})-\mathfrak{A}_{c,\varepsilon}}\|a\|^{-10}
	\int_{Z'(\mathbb{R})}\sum_{\substack{(\boldsymbol{\xi},t)\in \mathfrak{X}^*(f;a,z_{\infty})}}\prod_{p<\infty}E_p(\boldsymbol{\xi},t)\\
\int_{G'(\mathbb{R})} |\mathfrak F_\infty|\,dy_{\infty}d^{\times}z_{\infty}
d^{\times}adk_{\infty}^*dk_{\infty},
\end{multline}
where $C^{\dag}:=T^{-9n^2}M'^{2n+\varepsilon}M^n\cdot \mathcal{N}_f^{n-2}\cdot \prod_{p\mid M'}p^{ne_p(M)}$.

For $(\boldsymbol{\xi},t)\in\mathfrak X^*(f;a,z_\infty)$, the
non-Archimedean constraints also give
\begin{equation}\label{eq9.24.}
\prod_{p<\infty}E_p(\boldsymbol{\xi},t)
\ll (TMM')^{\varepsilon}\|a\|^{\varepsilon}.
\end{equation}
Indeed, the support conditions bound the heights of the rational parameters
by $(TMM'\mathcal N_f\|a\|)^{O(1)}$, so \eqref{eq9.24.} follows from the
standard divisor bound and $\log\mathcal N_f\ll\log T$. In the complementary
range, the factor $\|a\|^\varepsilon$ is absorbed by rapid decay; in the
essential range, it is absorbed into $T^\varepsilon$ by Definition
\ref{defn9.9}. We suppress this factor below.

\begin{itemize}
\item Suppose  $M'^2\mathcal{N}_f^2|z_{\infty}|\gg C_{\Omega^*}^2M^{\dagger}$. Then 
\begin{equation}\label{eq9.23}
\big|\mathfrak{X}^*(f;a,z_{\infty})\big|\ll \frac{M'^{2n}\mathcal{N}_f^{2n}|z_{\infty}|}{(C_{\Omega^*}^2M^{\dagger})^n}\prod_{j=1}^{n-1}\max\{a_j,a_j^{-1}\}=\frac{M'^{2n}\mathcal{N}_f^{2n}|z_{\infty}|}{(C_{\Omega^*}^2M^{\dagger})^n}\|a\|.
\end{equation}

\item Suppose  $M'^2\mathcal{N}_f^2|z_{\infty}|=o(C_{\Omega^*}^2M^{\dagger})$. Then by \eqref{10.17} we have $t=0,$ and thus 
\begin{align*}
M'^2\mathcal{N}_f^2|z_{\infty}|a_1\geq c' C_{\Omega^*}^2M^{\dagger}
\end{align*}
for some constant $c'>0$. Here $c'$ depends only on $c_{\infty}$ and $C_{\infty}$ defined in \textsection\ref{sec3.14}. Recall that $a_1\geq a_2\geq  \cdots \geq a_{n-1}>0$. Let $1\leq \ell_0\leq n-1$ be the largest integer such that $M'^2\mathcal{N}_f^2|z_{\infty}|a_{\ell_0}\geq 10^{-1}c' C_{\Omega^*}^2M^{\dagger}$. Then 
\begin{align*}
\big|\mathfrak{X}^*(f;a,z_{\infty})\big|\ll \prod_{j=1}^{n-1}\max\Bigg\{1,
\frac{M'^{2}\mathcal{N}_f^2|z_{\infty}|a_j}{C_{\Omega^*}^2M^{\dagger}}\Bigg\}\ll \Bigg[\frac{M'^{2}\mathcal{N}_f^{2}|z_{\infty}|}
{C_{\Omega^*}^2M^{\dagger}}\Bigg]^{\ell_0}
\prod_{j=1}^{\ell_0}\max\{a_j,a_j^{-1}\}.
\end{align*}

Since $|z_{\infty}|\ll 1$ and $\max\{a_j,a_j^{-1}\}\geq 1$ for all $1\leq j\leq n-1,$ then 
\begin{equation}\label{eq9.24}
\big|\mathfrak{X}^*(f;a,z_{\infty})\big|\ll \Bigg[\frac{M'^{2}\mathcal{N}_f^{2}}{C_{\Omega^*}^2M^{\dagger}}\Bigg]^{\ell_0}\cdot |z_{\infty}|\cdot \|a\|.
\end{equation}
\end{itemize}

Substituting \eqref{eq9.24.},  \eqref{eq9.23} and \eqref{eq9.24} into \eqref{e9.23} we derive that 
\begin{align*}
\frac{\mathcal{J}_{\Geo,\bi}^{\Reg,\RNum{2},\dagger}(f)}{\langle\phi',\phi'\rangle}\ll C^{\dag}T^{\frac{n(n+1)}{2}}(TMM')^{\varepsilon}\max_{1\leq \ell_0\leq n}\Bigg[\frac{M'^{2}\mathcal{N}_f^{2}}{C_{\Omega^*}^2M^{\dagger}}\Bigg]^{\ell_0}\int \frac{1}{\|a\|^{9}}\int |z_{\infty}|d^{\times}z_{\infty}d^{\times}a,
\end{align*}
where $a\in A_H(\mathbb{R})-\mathfrak{A}_{c,\varepsilon},$ and  $C_{\Omega^*}^2M^{\dagger}M'^{-2}\mathcal{N}_f^{-2}a_1^{-1}\ll |z_{\infty}|\ll 1$.

Since $|z_{\infty}|d^{\times}z_{\infty}$ is additive Haar measure up to normalization,
\begin{align*}
\int |z_{\infty}|d^{\times}z_{\infty}\leq \int_0^{O(1)}dz_{\infty}\ll 1.
\end{align*}

In conjunction with $\int \|a\|^{-9}d^{\times}a\ll 1$ we obtain 
\begin{align*}
\frac{\mathcal{J}_{\Geo,\bi}^{\Reg,\RNum{2},\dagger}(f)}{\langle\phi',\phi'\rangle}\ll (TMM')^{\varepsilon}T^{-8n^2}M'^{2n}M^n\mathcal{N}_f^{n-2}\prod_{p\mid M'}p^{ne_p(M)}\max_{1\leq \ell_0\leq n}\Bigg[\frac{M'^{2}\mathcal{N}_f^{2}}{C_{\Omega^*}^2M^{\dagger}}\Bigg]^{\ell_0}.
\end{align*}

Therefore, Proposition \ref{prop9.11.} follows.
\end{proof}

\subsection{Estimates in the Essential Range}\label{sec9.4}
By definition,
\begin{multline*}
\mathcal{J}_{\Geo,\bi}^{\Reg,\RNum{2},\heartsuit}(f):=
	\int_{\mathfrak{A}_{c,\varepsilon}}
	\int_{K'}
	\int_{\Omega^*}\big|\phi'(ak^*k)\big|^2\int_{G'(\mathbb{A})}
	\int_{Z'(\mathbb{A})}\\
\sum_{\substack{(\boldsymbol{\xi},t)\in \mathfrak{X}^*(f;a,z_{\infty})}}\big|\mathfrak{F}_{\infty}\cdot\mathfrak{F}_{\fin}\big|d^{\times}zdydk^*dk\frac{d^{\times}a}{\delta_{B'}(a)}.
\end{multline*}

%We now estimate the contribution from the essential range.

\begin{prop}\label{prop9.12} 
 Let notation be as before. Then
 \begin{align*}
 \frac{\mathcal{J}_{\Geo,\bi}^{\Reg,\RNum{2},\heartsuit}(f)}{\langle\phi',\phi'\rangle}\ll (TMM')^{\varepsilon}T^{\frac{n-1}{2}}M'^{2n}M^n\mathcal{N}_f^{n-2}\prod_{p\mid M'}p^{ne_p(M)} \cdot \max_{1\leq \ell_0\leq n}\Bigg[\frac{M'^{2}\mathcal{N}_f^{2}}{M^{\dagger}}\Bigg]^{\ell_0},
 \end{align*}
 with the same dependence as in Theorem \ref{thm9.1}.
 \end{prop}

\subsubsection{Separating the automorphic weight}\label{sec9.5.1}
By Lemma \ref{norm} we have 
\begin{equation}\label{9.14}
\int_{K'}\int_{\Omega^*}\frac{\big|\phi'(ak^*k)\big|^2}{\langle\phi',\phi'\rangle}dk^*dk\ll T^{o(1)}\delta_{B'}(a)h(a),
\end{equation}
where $h(a):=\min\big\{|a_1|_{\infty}^{-n}|\det a|_{\infty},|\det a|_{\infty}^{-1}\big\}$. 

\subsubsection{Further reductions}

Interchanging the order of integration gives
 \begin{align*}
\mathcal{J}_{\Geo,\bi}^{\Reg,\RNum{2},\heartsuit}(f)\ll 
	\int_{\mathfrak{A}_{c,\varepsilon}}I(a)\Bigg[\max_{\substack{k^*\in\Omega^*\\ k\in K'}}\int_{G'(\mathbb{A})}\int_{Z'(\mathbb{A})}\sum_{\substack{(\boldsymbol{\xi},t)}}\big|\mathfrak{F}_{\infty}\cdot \mathfrak{F}_{\fin}\big|d^{\times}zdy\Bigg]\delta_{B'}^{-1}(a)d^{\times}a,
\end{align*}
where $\mathfrak{F}_{\infty}$ and $\mathfrak{F}_{\fin}$ are defined by \eqref{55} and \eqref{58}, $(\boldsymbol{\xi},t)\in \mathfrak{X}^*(f;a,z_{\infty})$,
 and
\begin{align*}
 I(a):=\int_{K'}\int_{\Omega^*}\big|\phi'(ak^*k)\big|^2dk^*dk.
\end{align*}

Applying \eqref{9.14} to $I(a)$, we obtain
\begin{align*}
\frac{\mathcal{J}_{\Geo,\bi}^{\Reg,\RNum{2},\heartsuit}(f)}
{\langle\phi',\phi'\rangle}
\ll & T^{\varepsilon}\int_{\mathfrak{A}_{c,\varepsilon}}\bigg[
\max_{\substack{k^*\in\Omega^*\\ k\in K'}}
\int_{G'(\mathbb{R})}\int_{Z'(\mathbb{R})}
\sum_{(\boldsymbol{\xi},t)\in \mathfrak{X}^*(f;a,z_{\infty})}|\mathfrak{F}_{\infty}|\\
&\hspace{35mm}\times
\prod_{p<\infty}\mathcal{J}_p^{\heartsuit}
(\boldsymbol{\xi},t;k_{p}^*k_p)d^{\times}z_{\infty}dy_{\infty}
\bigg]h(a)d^{\times}a,
\end{align*}
where 
$\mathcal{J}_p^{\heartsuit}(\boldsymbol{\xi},t;k_{p}^*k_p)$ is defined in
\textsection\ref{9.2.1} and $h(a)$ is defined in \textsection\ref{sec9.5.1}.  

Lemma \ref{lem9.7} now gives
\begin{align*}
\mathcal{J}_{\Geo,\bi}^{\Reg,\RNum{2},\heartsuit}(f)\ll &T^{\varepsilon}M'^{2n+\varepsilon}M^n\mathcal{N}_f^{n-2}\langle\phi',\phi'\rangle\prod_{p\mid M'}p^{ne_p(M)} \int_{\mathfrak{A}_{c,\varepsilon}}\max_{\substack{k_{\infty}^*\in\Omega_{\infty}^*\\ k_{\infty}\in K_{\infty}'}}\int_{Z'(\mathbb{R})}
	\\
&\sum_{\substack{(\boldsymbol{\xi},t)\in \mathfrak{X}^*(f;a,z_{\infty})}}\int_{G'(\mathbb{R})}\big|\mathfrak{F}_{\infty}\big|dy_{\infty}h(a)\prod_{p<\infty}E_p(\boldsymbol{\xi},t)d^{\times}z_{\infty}d^{\times}a.
\end{align*}

\subsubsection{The Archimedean inner integral}
Denote by 
\begin{equation}\label{9.19.}
g=\iota(z_{\infty}ak_{\infty}^*k_{\infty})^{-1}\gamma(\boldsymbol\xi,t)
\iota(z_{\infty}ak_{\infty}^*k_{\infty})\in G(\mathbb{R}).
\end{equation}

\begin{lemma}\label{lem9.13.}
Let notation be as before. Then 
\begin{align*}
	\int_{G'(\mathbb{R})}\big|\mathfrak{F}_{\infty}\big|dy_{\infty}\ll T^{\frac{n}{2}+\varepsilon}\cdot\min\Bigg\{1,\frac{T^{-1/2}}{d_{G'}(g)}\Bigg\},
\end{align*}
where $d_{G'}$ is defined by \eqref{dg} and $g$ by \eqref{9.19.}. The
implied constant has the same dependence as in Theorem \ref{thm9.1}.
\end{lemma}
\begin{proof}
Put $r=T^{-1/2+\varepsilon}$. By \eqref{245}, the integrand vanishes unless
$g\iota(y_{\infty})$ lies in a fixed neighborhood of the identity and
\begin{equation}\label{9.18.}
\dist\bigl(g\iota(y_{\infty})\tau,\tau\bigr)\ll r.
\end{equation}
In local coordinates on $G'(\mathbb R)$, condition \eqref{9.18.} restricts
$y_{\infty}$ to a ball of radius $O(r)$. Since $\dim G'=n^2$, its volume is
$O(r^{n^2})$. Averaging in the fixed central neighborhood $\mathcal Z$ from
Proposition \ref{prop3.1} and applying that proposition gives the sharper
volume estimate
\begin{equation}\label{9.18..}
\Vol\bigl\{y_{\infty}:\eqref{9.18.}\text{ holds}\bigr\}
\ll r^{n^2}\min\left\{1,\frac{r}{d_{G'}(g)}\right\}.
\end{equation}
Here right multiplication by elements in a fixed neighborhood of $G'$ changes
$d_{G'}(g)$ only by a bounded factor. Finally, \eqref{250} and
\eqref{9.18..} yield
\begin{align*}
\int_{G'(\mathbb{R})}|\mathfrak{F}_{\infty}|dy_{\infty}
\ll T^{\frac{n(n+1)}{2}+\varepsilon}r^{n^2}
\min\left\{1,\frac{r}{d_{G'}(g)}\right\}
\ll T^{\frac n2+\varepsilon}
\min\left\{1,\frac{T^{-1/2}}{d_{G'}(g)}\right\}.
\end{align*}
\end{proof}

\subsubsection{The outer integrals}
We use the preceding estimates \eqref{9.7}, \eqref{10.17}, and
\eqref{eq9.24.}.

The main results in this section are the following two lemmas.
\begin{lemma}\label{lem9.14}
Uniformly for $a\in\mathfrak{A}_{c,\varepsilon}$,
$k_\infty^*\in\Omega_\infty^*$, and $k_\infty\in K_\infty'$, we have
\begin{equation}\label{lem9.14.}
\int_{Z'(\mathbb{R})}\frac{T^{-1/2}}{d_{G'}(g)}
\sum_{\substack{(\boldsymbol{\xi},t)\in
\mathfrak{X}^*(f;a,z_{\infty})\\
\boldsymbol{\xi}=\mathbf{0}}}
h(a)\prod_{p<\infty}E_p(\boldsymbol{\xi},t)
d^{\times}z_{\infty}
\ll
\frac{(TMM')^{\varepsilon}M'^2\mathcal{N}_f^2}
{C_{\Omega^*}^2M^{\dagger}T^{1/2}},
\end{equation}
where $g$ is defined by \eqref{9.19.}. The implied constant has the same
dependence as in Theorem \ref{thm9.1}.
\end{lemma}
\begin{proof}
Since $\boldsymbol\xi=\mathbf0$ forces $t\neq0$, the defining constraints of
$\mathfrak X^*(f;a,z_\infty)$ give
\begin{equation}\label{9.22e}
\sum_{\substack{(\boldsymbol{\xi},t)\in
\mathfrak{X}^*(f;a,z_{\infty})\\
\boldsymbol{\xi}=\mathbf{0}}}1
\ll
\frac{M'^2\mathcal{N}_f^2|z_{\infty}|}
{C_{\Omega^*}^2M^{\dagger}}.
\end{equation}
	
Note that 
$h(a)\leq |\det a|_{\infty}^{-1}\ll T^{\varepsilon}$ for all $a\in \mathfrak{A}_{c,\varepsilon}$.  Combining \eqref{eq9.24.} with \eqref{9.22e} we have 
\begin{align*}
\sum_{\substack{(\boldsymbol{\xi},t)\in \mathfrak{X}^*(f;a,z_{\infty})\\ \boldsymbol{\xi}=\mathbf{0}}}h(a)\prod_{p<\infty}E_p(\boldsymbol{\xi},t) \ll \frac{M'^2\mathcal{N}_f^2|z_{\infty}|}{C_{\Omega^*}^2M^{\dagger}}\cdot (TMM')^{\varepsilon}.
\end{align*}

From \eqref{dg} and the explicit matrix in \eqref{9.19.}, uniformly for
$a\in\mathfrak A_{c,\varepsilon}$ and the compact variables
$k_{\infty}^*,k_{\infty}$,
\begin{align*}
	d_{G'}(g)\gg\min\big\{1,|z_{\infty}|\big\}\gg |z_{\infty}|. 
\end{align*}

Consequently, the left-hand side of \eqref{lem9.14.} is
\begin{align*}
\ll
\frac{(TMM')^{\varepsilon}M'^2\mathcal{N}_f^2}
{C_{\Omega^*}^2M^{\dagger}T^{1/2}}
\int\frac{d|z_\infty|}{|z_\infty|}
\ll
\frac{(TMM')^{\varepsilon}M'^2\mathcal{N}_f^2}
{C_{\Omega^*}^2M^{\dagger}T^{1/2}},
\end{align*}
since \eqref{10.17} and $\log\mathcal N_f\ll\log T$ give
$\int d|z_\infty|/|z_\infty|\ll(TMM')^\varepsilon$.
\end{proof}

\begin{lemma}\label{lem9.15}
Uniformly for $a\in\mathfrak{A}_{c,\varepsilon}$,
$k_\infty^*\in\Omega_\infty^*$, and $k_\infty\in K_\infty'$, we have
\begin{equation}\label{lem9.15.}
\int \frac{T^{-\frac{1}{2}}}{d_{G'}(g)}\sum_{\substack{(\boldsymbol{\xi},t)\in
\mathfrak X^*(f;a,z_\infty)\\ \boldsymbol{\xi}\neq \mathbf{0}}}h(a)\prod_{p<\infty}E_p(\boldsymbol{\xi},t)d^{\times}z_{\infty}\ll \frac{(TMM')^{\varepsilon}(C_{\Omega^*}^2M^{\dagger}+M'^{2}\mathcal{N}_f^2)}{T^{\frac{1}{2}}(C_{\Omega^*}^2M^{\dagger})^{n}(M'^{2}\mathcal{N}_f^2)^{1-n}},
\end{equation}
where $z_{\infty}\in Z'(\mathbb{R})$,
$g$ is defined by \eqref{9.19.}, and the implied constant has the same dependence as in
Theorem \ref{thm9.1}.
\end{lemma}

\begin{proof}
For $n=1$ the sum is empty, so assume $n\geq2$. If
$\boldsymbol{\xi}\neq \mathbf{0}$, then
$M'^{2}\mathcal{N}_f^2|z_{\infty}|a_1\gg C_{\Omega^*}^2M^{\dagger}$.
For $a\in\mathfrak A_{c,\varepsilon}$, Definition \ref{defn9.9} and
$a_1\geq\cdots\geq a_{n-1}>0$ give
\begin{equation}\label{eq9.35}
T^{-c}\leq a_1^{-1}\ll T^{\varepsilon},
\end{equation}
where $c=2c_1+10n^2$ and $c_1$ is the constant in Lemma \ref{lem4.9}. Note that $c_1$ (and thus $c$) depends only on $n$.

Using \eqref{9.7}, we deduce that
\begin{align*}
\sum_{\substack{(\boldsymbol{\xi},t)\in \mathfrak{X}^*(f;a,z_{\infty})\\ \boldsymbol{\xi}\neq \mathbf{0}}}1\ll \max\Bigg\{1,\frac{M'^2\mathcal{N}_f^2|z_{\infty}|}{C_{\Omega^*}^2M^{\dagger}}\Bigg\}\cdot\Bigg[\frac{M'^{2}\mathcal{N}_f^2|z_{\infty}|}{C_{\Omega^*}^2M^{\dagger}}\Bigg]^{n-1}\cdot a_1^{n-1}.
\end{align*}

Using $h(a)\leq |\det a|_\infty/a_1^n$ and \eqref{eq9.35}, and since
$|z_\infty|\ll1$ and $n\geq2$, it follows that
\begin{equation}\label{9.36}
\sum_{\substack{(\boldsymbol{\xi},t)\in \mathfrak{X}^*(f;a,z_{\infty})\\ \boldsymbol{\xi}\neq \mathbf{0}}}h(a)\prod_{p<\infty}E_p(\boldsymbol{\xi},t)
\ll CT^{\varepsilon}|z_\infty|^{n-1}\ll CT^{\varepsilon}|z_\infty|,
\end{equation}
where $C:=(TMM')^{\varepsilon}\cdot \max\big\{1,M'^2\mathcal{N}_f^2(C_{\Omega^*}^2M^{\dagger})^{-1}\big\}\cdot\big[M'^{2}\mathcal{N}_f^2(C_{\Omega^*}^2M^{\dagger})^{-1}\big]^{n-1}$. 

Moreover,
\begin{equation}\label{eq9.35z}
C_{\Omega^*}^2M^{\dagger}(M'^{2}\mathcal{N}_f^2)^{-1}a_1^{-1}\ll |z_{\infty}|\ll 1.
\end{equation}

As in the proof of Lemma \ref{lem9.14}, the explicit form of $g$ gives
$d_{G'}(g)\gg |z_{\infty}|$. Therefore, by \eqref{eq9.35}, \eqref{9.36},
and \eqref{eq9.35z}, the left-hand side of \eqref{lem9.15.} is
\begin{align*}
	\ll CT^{\varepsilon}\int\frac{T^{-\frac{1}{2}}}{|z_{\infty}|}d|z_{\infty}|\ll (TMM')^{\varepsilon}T^{-\frac{1}{2}}\max\Bigg\{1,\frac{M'^2\mathcal{N}_f^2}{C_{\Omega^*}^2M^{\dagger}}\Bigg\}\cdot\Bigg[\frac{M'^{2}\mathcal{N}_f^2}{C_{\Omega^*}^2M^{\dagger}}\Bigg]^{n-1},
\end{align*}
where the range of $z_\infty$ follows from \eqref{eq9.35} and
\eqref{eq9.35z}.
\end{proof}

\subsubsection{Completion of the estimate}
Set $X_0:=(C_{\Omega^*}^2M^\dagger)^{-1}M'^2\mathcal N_f^2$. For $n\geq2$, Lemmas \ref{lem9.14} and \ref{lem9.15} contribute
$X_0$ and $X_0^{n-1}+X_0^n$, respectively; for $n=1$, the latter
contribution is absent. Since $C_{\Omega^*}$ is fixed, in either case their
sum is
\begin{align*}
\ll\max_{1\leq\ell_0\leq n}
\left((M^\dagger)^{-1}M'^2\mathcal N_f^2\right)^{\ell_0}.
\end{align*}
Together with Lemma \ref{lem9.13.} and
$\int_{\mathfrak A_{c,\varepsilon}}d^\times a\ll T^\varepsilon$, this gives
\begin{align*}
\frac{\mathcal{J}_{\Geo,\bi}^{\Reg,\RNum{2},\heartsuit}(f)}{\langle\phi',\phi'\rangle}\ll (TMM')^{\varepsilon}T^{\frac{n-1}{2}}M'^{2n}M^n\mathcal{N}_f^{n-2}\prod_{p\mid M'}p^{ne_p(M)}\cdot \max_{1\leq \ell_0\leq n}\Bigg[\frac{M'^{2}\mathcal{N}_f^{2}}{M^{\dagger}}\Bigg]^{\ell_0},
\end{align*}
where the implied constant depends on $\varepsilon,$ parameters $c_{\infty}$ and $C_{\infty}$ defined in \textsection\ref{sec3.14}, and the  conductor of $\pi_{\infty}'$. Hence 
Proposition \ref{prop9.12} follows.

\subsection{Proof of Proposition \ref{prop9.3.}}\label{sec9.5}

Combining
\begin{align*}
\mathcal{J}_{\Geo,\bi}^{\Reg,\RNum{2}}(f)\ll \mathcal{J}_{\Geo,\bi}^{\Reg,\RNum{2},\dagger}(f)+\mathcal{J}_{\Geo,\bi}^{\Reg,\RNum{2},\heartsuit}(f)
\end{align*}
with Propositions \ref{prop9.11.} and \ref{prop9.12}, we obtain
\begin{align*}
	\frac{\mathcal{J}_{\Geo,\bi}^{\Reg,\RNum{2}}(f)}{\langle\phi',\phi'\rangle}\ll (TMM')^{\varepsilon}T^{\frac{n-1}{2}}M'^{2n}M^n\mathcal{N}_f^{n-2}\prod_{p\mid M'}p^{ne_p(M)} \max_{1\leq \ell_0\leq n}\Bigg[\frac{M'^{2}\mathcal{N}_f^{2}}{M^{\dagger}}\Bigg]^{\ell_0}.
\end{align*}	
The same argument applies to $f^{-1}$: inversion preserves all relevant local
support and sup-norm bounds, interchanges the two amplification primes in the
off-diagonal case, and leaves $\mathcal N_f$ unchanged. This proves
Proposition \ref{prop9.3.}.

\section{Applications to the Subconvexity Problem}\label{proof}
\subsection{The Spectral Side}
The spectral lower bound established in \textsection\ref{sec5} is recalled
below.
\thmf* 

\subsection{The Geometric Side}\label{sec11.2}
By Lemma \ref{lem3.2},
$J_{\Geo,\bi}^{\Reg,\RNum{1}}(f,\mathbf{s})\equiv0$, and hence
\begin{equation}\label{eq10.1}
J_{\Geo}^{\Reg,\heartsuit}(f,\mathbf{s})
=
J_{\Geo,\Main}^{\sm}(f,\mathbf{s})
+
J_{\Geo,\Main}^{\du}(f,\mathbf{s})
+
J_{\Geo,\bi}^{\Reg,\RNum{2}}(f,\mathbf{s}).
\end{equation}
Combining \eqref{eq10.1} with Propositions \ref{prop54} and \ref{prop8.1}
and Theorem \ref{thm9.1}, we obtain the following estimate.
\begin{cor}\label{cor11.3}
Let
$f\in\big\{f_{p_0,i}^{\diag},f_{p_1,p_2}^{\mathrm{off}}\big\}$ and set
$\vartheta_f:=\max_{p\mid\nu(f)}\vartheta_p$. Then
\begin{align*}
\frac{J_{\Geo}^{\Reg,\heartsuit}(f,\mathbf{0})}{\langle\phi',\phi'\rangle}&\ll T^{\frac{n}{2}+\varepsilon}M'^{2n}M^{n+\varepsilon} \mathcal{N}_f^{-1+2\vartheta_f+\varepsilon}\prod_{p\mid M'}p^{ne_p(M)}\\
	&+(TMM')^{\varepsilon}T^{\frac{n-1}{2}}M'^{2n}M^n\mathcal{N}_f^{n-2}\prod_{p\mid M'}p^{ne_p(M)} \max_{1\leq \ell_0\leq n}\Bigg[\frac{M'^{2}\mathcal{N}_f^{2}}{M^{\dagger}}\Bigg]^{\ell_0},
\end{align*}
where the implied constant depends on $\varepsilon$, $c_{\infty}$,
$C_{\infty}$, and the fixed data attached to $\pi'$ and $p_\circ$.
\end{cor}

The preceding estimate suffices for the diagonal terms. For the
off-diagonal main terms, however, summing the local bounds separately loses a
power of the amplifier. We instead retain their Hecke structure. Set $\rho:=\widetilde\pi'\boxplus(\omega\omega')$, a unitary isobaric representation of $G(\mathbb A)$ with central character
$\omega$. Its local components are unramified at every $p\in\mathcal L$.

\begin{lemma}[The Local Hecke Transform]\label{lem10.hecke}
Let $p\in\mathcal L$ and $\ell\geq0$. For $s$ near $0$, set $\mathcal H_p(h,s):=L_p(1+s,\pi_p'\times\widetilde\pi_p')^{-1}\mathcal I_p(h,s)$. Then
\begin{align*}
\mathcal H_p(T_{p^\ell},0)=p^{-\ell/2}\Lambda_{\rho,p}(\ell),\quad
\mathcal H_p(T_{p^\ell}^*,0)=p^{-\ell/2}\overline{\Lambda_{\rho,p}(\ell)},
\end{align*}
where $\Lambda_{\rho,p}(\ell)$ is the normalized Satake transform of $T_{p^\ell}$
evaluated at the Satake parameters of $\rho_p$. Moreover, for
$0\leq\ell\leq n+1$, we have $|\Lambda_{\rho,p}(\ell)|^2
\ll_n a_{\rho\times\widetilde\rho}(p^\ell)$ and 
\begin{align*}
\left|\mathcal H_p'(T_{p^\ell},0)\right|
+\left|\mathcal H_p'(T_{p^\ell}^*,0)\right|
\ll_n p^{-\ell/2}(1+\ell)
a_{\rho\times\widetilde\rho}(p^\ell)^{1/2}\log p,
\end{align*}
where $a_{\rho\times\widetilde\rho}$ denotes the Rankin--Selberg
coefficients.
\end{lemma}
\begin{proof}
Let $\alpha_{1,p},\ldots,\alpha_{n,p}$ be the Satake parameters of $\pi_p'$
and put
\begin{align*}
\boldsymbol\beta_p
:=(\alpha_{1,p}^{-1},\ldots,\alpha_{n,p}^{-1},\Omega_p),\qquad
\Omega_p:=\omega_p(p)\omega_p'(p).
\end{align*}
Applying the Cartan decomposition and Macdonald's formula to the unfolded
expression in Lemma \ref{lem7.2} gives
\begin{align*}
p^{\ell/2}\mathcal H_p(T_{p^\ell},0)
=\Lambda_{\rho,p}(\ell),
\end{align*}
where $\Lambda_{\rho,p}(\ell)$ is the Hall--Littlewood polynomial associated
with $K_p\diag(p^\ell,I_n)K_p$, evaluated at $\boldsymbol\beta_p$. Its
Kostka expansion (see \cite[\textsection 4]{BM15}) is of the form
\begin{align*}
\Lambda_{\rho,p}(\ell)
=\sum_{\lambda\vdash\ell}c_{\ell,\lambda}(p^{-1})
s_\lambda(\boldsymbol\beta_p),\quad
c_{\ell,\lambda}(p^{-1})\ll_n1.
\end{align*}
Here $\lambda$ ranges over partitions of $\ell$, and $s_\lambda$ denotes the
corresponding Schur polynomial. 
The local Cauchy identity gives
\begin{align*}
a_{\rho\times\widetilde\rho}(p^\ell)
=\sum_{\lambda\vdash\ell}|s_\lambda(\boldsymbol\beta_p)|^2.
\end{align*}
The first estimate follows by Cauchy--Schwarz. Before setting $s=0$, the same
Cartan calculation weights the finitely many terms above by powers $p^{-rs}$
with $|r|\ll_n\ell$. Differentiation gives the second estimate, while the
adjoint identity follows by taking contragredients.
\end{proof}

\begin{remark}
This parallels \cite[\textsection 5.8]{Nel21}, where the main term retains
the $H$-Hecke eigenvalue of the wave packet. Here the genuine cuspidal
$H$-weight and the central-character projection produce instead the
Satake polynomials attached to $\widetilde\pi'\boxplus(\omega\omega')$.
When $\ell=1$, one has $\Lambda_{\rho,p}(1)=\lambda_\rho(p)$.
\end{remark}

\begin{prop}[Collective Estimate for the Main Terms]\label{prop10.main}
Suppose that $|\alpha_p|\leq1$ for every $p\in\mathcal L$, and let
\begin{align*}
\mathcal M_{\mathrm{off}}:=\sum_{\substack{p_1,p_2\in\mathcal L\\p_1\neq p_2}}
\alpha_{p_1}\overline{\alpha_{p_2}}\cdot
\left(
J_{\Geo,\Main}^{\sm}(f_{p_1,p_2}^{\mathrm{off}},\mathbf0)
+J_{\Geo,\Main}^{\du}(f_{p_1,p_2}^{\mathrm{off}},\mathbf0)
\right).
\end{align*}
Then
\begin{align*}
|\mathcal M_{\mathrm{off}}|
\ll
(TMM')^\varepsilon T^{n/2}M'^{2n}M^n
\langle\phi',\phi'\rangle
\prod_{p\mid M'}p^{ne_p(M)}\cdot L.
\end{align*}
\end{prop}
\begin{proof}
After dividing the amplification-prime factors by
$L_p(1+s,\pi_p'\times\widetilde\pi_p')$, the global Rankin--Selberg factor is
independent of $p_1$ and $p_2$. Lemma \ref{lem10.hecke} therefore reduces the
off-diagonal value terms to products of
$b_p:=\alpha_p p^{-\ell_p/2}\Lambda_{\rho,p}(\ell_p)$ and their conjugates.
Their sum over $p_1\neq p_2$ is bounded by
\begin{align*}
\big|\sum_{p\in\mathcal L}b_p\big|^2+\sum_{p\in\mathcal L}|b_p|^2.
\end{align*}
Set $s_L:=1+(\log L)^{-1}$. Since $p\asymp L$ and
$1\leq\ell_p\leq n+1$, we have
$p^{-\ell_p}\ll_n p^{-\ell_p s_L}$. Lemma \ref{lem10.hecke} therefore gives
\begin{align*}
\sum_{p\in\mathcal L}p^{-\ell_p}|\Lambda_{\rho,p}(\ell_p)|^2\ll_n
\sum_{p\in\mathcal L}
\frac{a_{\rho\times\widetilde\rho}(p^{\ell_p})}
{p^{\ell_p s_L}}\leq L(s_L,\rho\times\widetilde\rho)
\ll_\varepsilon (MM'L)^\varepsilon.
\end{align*}
Cauchy--Schwarz and $\#\mathcal L\ll L/\log L$ therefore give
\begin{align*}
\big|\sum_{p\in\mathcal L}b_p\big|^2\ll_\varepsilon (MM'L)^\varepsilon L.
\end{align*}
The derivative bound in Lemma \ref{lem10.hecke} is treated in the same way,
using $\rho\times\widetilde\rho$; its logarithmic factors are absorbed into
$L^\varepsilon$. The dual term gives the contragredient coefficients and
satisfies the same estimate. The remaining local and Archimedean factors are
bounded as in Propositions \ref{prop54} and \ref{prop8.1}, completing the
proof.
\end{proof}

\begin{comment}
\begin{prop}\label{prop11.3}
Let notation be as before. Assume that $\mathcal{N}_f\leq T$. Then
\begin{align*}
		\mathcal{J}_{\Geo}^{\heartsuit}(\boldsymbol{\alpha},\boldsymbol{\ell})\ll & T^{\frac{n}{2}+\varepsilon} M^{n} \Bigg[\sum_{p\in \mathcal{L}}|\alpha_p|^2+\bigg|\sum_{p\in \mathcal{L}}|\alpha_p|p^{-\frac{2\vartheta_p+1}{2}\ell_{p}}\bigg|^2+\frac{1}{T^{\frac{1}{2}}M^{n}}\bigg|\sum_{p\in \mathcal{L}}|\alpha_p|p^{2n\ell_{p}}\bigg|^2\Bigg],
\end{align*}
where the implied constant relies on $\varepsilon$ and $M'$.
\end{prop}
\begin{proof}
\begin{align*}
	J_{\Geo}^{\Reg,\heartsuit}(f,\mathbf{0})&\ll_{\varepsilon, \pi'}T^{\frac{n}{2}+\varepsilon} M^{n} (\mathcal{N}_f^{-1-2\vartheta_p}+T^{-\frac{1}{2}+\varepsilon}\mathcal{N}_f^{2n})+T^{\frac{n}{2}-\frac{1}{2}+\varepsilon} (M^{n-1} \mathcal{N}_f^{2n}+\mathcal{N}_f^{4n}).
\end{align*}

\end{proof}

\end{comment}

\subsection{A General Geometric-Side Bound}

\subsubsection{Choice of the set $\mathcal{L}$}\label{sec11.3.1}
Let $L\gg1$ satisfy $L\ll(TM)^{O(1)}$, and define
\begin{equation}\label{L}
\mathcal L:=\big\{\text{prime $p$}:\ L<p\leq2L,\ p\nmid p_\circ MM'\big\}.
\end{equation}
For the choices of $L$ made below, the number of prime divisors of
$p_\circ MM'$ in $(L,2L]$ is negligible. Hence the prime number theorem gives $\#\mathcal L\asymp L/\log L$. 

\subsubsection{Choice of $\boldsymbol{\ell}$ and $\boldsymbol{\alpha}$}\label{alpha}
Let $\mu\in i\mathfrak a_Q^*/i\mathfrak a_G^*$ be as in Theorem
\ref{prop11.2}. The Hecke relations \cite[Corollary 4.3]{BM15} give
\begin{align*}
\sum_{j=1}^{n+1}|\lambda_{\pi_{\mu}}(p^j)|\gg 1,\quad p\nmid M.
\end{align*}
Thus, for each $p\in\mathcal L$, there exists
$\ell_p\in\{1,\ldots,n+1\}$ such that
$|\lambda_{\pi_\mu}(p^{\ell_p})|\gg1$. Fix one such $\ell_p$ and set
$\alpha_p:=\overline{\lambda_{\pi_{\mu}}(p^{\ell_p})}/
|\lambda_{\pi_{\mu}}(p^{\ell_p})|$.

\subsubsection{The amplified bound}
Recall from \textsection\ref{sec3.7.2} that
\begin{align*}
\mathcal{J}_{\Geo}^{\heartsuit}(\boldsymbol{\alpha},\boldsymbol{\ell})=\sum_{p_1\neq p_2\in\mathcal{L}}\alpha_{p_1}\overline{\alpha_{p_2}}J_{\Geo}^{\Reg,\heartsuit}(f_{p_1,p_2}^{\mathrm{off}},\mathbf{0})+\sum_{p_0\in\mathcal{L}}\sum_{i=0}^{\ell_{p_0}}c_{p_0,i}|\alpha_{p_0}|^2J_{\Geo}^{\Reg,\heartsuit}(f_{p_0,i}^{\diag},\mathbf{0}).
\end{align*}

\begin{prop}\label{cor11.4}
Let $\mathcal L$ be as in \eqref{L}, choose
$\boldsymbol\ell$ and $\boldsymbol\alpha$ as in \textsection\ref{alpha}, and
let $\phi'$ be the vector of Definition \ref{def4.2}. Set
$R:=M/M^\dagger$. Then
\begin{multline}\label{10.4}
\frac{\mathcal J_{\Geo}^{\heartsuit}
(\boldsymbol\alpha,\boldsymbol\ell)}
{\langle\phi',\phi'\rangle T^{n/2}}
\ll (TMM')^\varepsilon
\big(M^nR^nL
+T^{-1/2}M^{n-1}R^{n+1}L^{n(n+1)+2}\big)\\
+(TMM')^\varepsilon
T^{-1/2}R^{2n}L^{(3n-2)(n+1)+2},
\end{multline}
where the implied constant depends on $\varepsilon$, $c_\infty$,
$C_\infty$, and the conductor of $\pi'$.
\end{prop}

\begin{proof}
The factor $M'^{2n}$ in Corollary \ref{cor11.3} is absorbed into the
implied constant, while $\prod_{p\mid M'}p^{ne_p(M)}=R^n$. Proposition \ref{prop10.main} bounds the off-diagonal small-cell and dual
terms by the first term on the right of \eqref{10.4}. For the diagonal main
terms, Corollary \ref{cor11.3} and the uniform gap
$1-2\vartheta_p\gg_n1$ give
\begin{align*}
\sum_{p\in\mathcal L}\sum_{i=0}^{\ell_p}
|c_{p,i}|p^{-i(1-2\vartheta_p)+\varepsilon}
\ll_\varepsilon L^{1+\varepsilon};
\end{align*}
here the term $i=0$ dominates.

It remains to sum the regular orbital contribution. By Theorem
\ref{thm9.1}, $|\alpha_p|=1$, $c_{p,i}\ll1$,
$\#\mathcal L\ll L/\log L$, and $\mathcal N_f\ll L^{n+1}$, its total is
bounded, after division by $\langle\phi',\phi'\rangle T^{n/2}$, by
\begin{align*}
(TMM')^\varepsilon T^{-1/2}
\big(M^{n-1}R^{n+1}L^{n(n+1)+2}
+R^{2n}L^{(3n-2)(n+1)+2}\big).
\end{align*}
Combining these estimates proves \eqref{10.4}.
\end{proof}

%Here we also appeal to Lemma \ref{lem4.6} to bound $\langle\phi',\phi'\rangle$ by $O(T^{\varepsilon})$.

\subsection{Completion of the General Case}\label{sec11.4}
We now combine the spectral and geometric estimates.
\begin{thmx}\label{E}
Retain the notation and assumptions of Theorem \ref{B}. Then
\begin{equation}\label{nb}
\frac{L(1/2,\pi\times\pi')}
{\sqrt{\langle\phi',\phi'\rangle}
\prod_{j=1}^m\sqrt{L(1,\pi_j,\Ad)}}
\ll
T^{\frac{n(n+1)}{4}+\varepsilon}
M^{\frac n2+\varepsilon}R^{\frac n2}\mathbf L^{-1},
\end{equation}
where the implied constant depends on $\varepsilon$, $c_\infty$,
$C_\infty$, and the conductor of $\pi'$.
\end{thmx}

\begin{proof}
The case $|L(1/2,\pi\times\pi')|<2$ follows from the standard lower bounds
at $s=1$ for the adjoint $L$-values and Lemma \ref{lem4.6}. We may therefore
assume that $|L(1/2,\pi\times\pi')|\geq2$.
By the choice of $\boldsymbol\ell$ and $\boldsymbol\alpha$,
\begin{align*}
	\sum_{p\in \mathcal{L}}\alpha_p\lambda_{\pi_{\mu}}(p^{\ell_p})=\sum_{p\in \mathcal{L}}\big|\lambda_{\pi_{\mu}}(p^{\ell_p})\big|\gg \#\mathcal{L}\gg \frac{L}{\log L}.
\end{align*}
Theorem \ref{prop11.2} and Proposition \ref{cor11.4} therefore imply that
\begin{multline*}
\frac{\big|L(1/2,\pi\times\pi')\big|^2}
{\langle\phi',\phi'\rangle T^{\frac{n(n+1)}{2}}}
\prod_{j=1}^m\frac{1}{L(1,\pi_j,\Ad)}
\ll (TML)^\varepsilon
\left(
\frac{M^nR^n}{L}
+\frac{R^{2n}L^{3n^2+n-2}}{T^{1/2}}
\right)\\
+(TML)^\varepsilon
\frac{M^{n-1}R^{n+1}L^{n(n+1)}}{T^{1/2}}.
\end{multline*}
If $(M^\dagger)^{n^2-n-1}\leq T^{n+1}$, balance the first term with
$T^{-1/2}R^{2n}L^{3n^2+n-2}$; otherwise balance it with
$T^{-1/2}M^{n-1}R^{n+1}L^{n(n+1)}$. In either case the optimal amplifier
length is $L=\mathbf L^2$, and \eqref{nb} follows.
\end{proof}

Since $\pi'$ is fixed, Lemma \ref{lem4.6} gives
$1\ll_{\pi'}\langle\phi',\phi'\rangle\ll_{\pi'}T^{o(1)}$; the latter factor
is absorbed into $T^{\varepsilon}$. Hence Theorem \ref{B} follows.

\begin{remark}
Corollary \ref{cor11.3} also shows that the implied constant in Theorem
\ref{B} depends at most polynomially on the conductor of $\pi'$.
\end{remark}

\subsection{\texorpdfstring{The Geometric-Side Bound in the $t$-Aspect}
{The Geometric-Side Bound in the t-Aspect}}\label{sec12.5}
\subsubsection{Choice of the amplifier}\label{12.4.1}
Take $\mathcal L$ as in \textsection\ref{sec11.3.1}, set $\ell_p=1$ for
every $p\in\mathcal L$, and let
$\boldsymbol\alpha=(\alpha_p)_{p\in\mathcal L}$ satisfy $|\alpha_p|\leq1$.

\subsubsection{The amplified bound}
The same argument, now with $\ell_p=1$, gives the following sharper analogue
of Proposition \ref{cor11.4}.
 
\begin{prop}\label{prop12.5}
Let notation be as before. Let $\mathcal{L}$, $\boldsymbol{\ell}$, and
$\boldsymbol{\alpha}$ be defined in \textsection\ref{12.4.1}, and let
$\phi'\in\pi'$ be the vector of Definition \ref{def4.2}. Set
$R:=M/M^\dagger$. Then
\begin{align*}
\frac{\mathcal{J}_{\Geo}^{\heartsuit}(\boldsymbol{\alpha},\boldsymbol{\ell})}
{\langle\phi',\phi'\rangle T^{n/2}}
\ll (TMM')^\varepsilon
\left(
M^nR^nL
+T^{-1/2}M^{n-1}R^{n+1}L^{n+2}
+T^{-1/2}R^{2n}L^{3n}
\right),
\end{align*}
where the implied constant depends on $\varepsilon$, $c_{\infty}$,
$C_{\infty}$, and the conductor of $\pi'$.
\end{prop}
\begin{proof}
Repeat the proof of Proposition \ref{cor11.4}, using
$\mathcal N_f\asymp L$ for the off-diagonal terms. The diagonal contribution
is again dominated by $i=0$.
\end{proof}

\subsubsection{Proof of Theorem \ref{A}: the tempered case}
The tempered estimate in Theorem \ref{A} follows from the next result.
\begin{thmx}\label{J}
Suppose that $\pi$ is tempered. Then
\begin{equation}\label{nb.}
L(1/2+it,\pi\times\pi')
\ll_{\varepsilon,\pi,\pi'}
(1+|t|)^{
\frac{n(n+1)}{4}
-\frac{1}{4(3n-1)}
+\varepsilon}.
\end{equation}
\end{thmx}

\begin{proof}
Set $\sigma:=\pi\otimes|\det|^{it}$ and $T:=1+|t|$, and let $\mu$ be the
deformation associated with $\sigma$ in Theorem \ref{prop11.2}. For $T$
sufficiently large, we may assume that $L>M$. Since $\pi$ is
tempered, $|\lambda_{\sigma_\mu}(p)|
=
|\lambda_{\pi_\mu}(p)|
\ll_n 1$ for all $p\nmid M$. Choose
\begin{align*}
\alpha_p:=
\begin{cases}
\overline{\lambda_{\sigma_\mu}(p)}/|\lambda_{\sigma_\mu}(p)|,
&\lambda_{\sigma_\mu}(p)\neq0,\\
0,&\lambda_{\sigma_\mu}(p)=0.
\end{cases}
\end{align*}
Then Rankin--Selberg theory gives
\begin{align*}
\big|\sum_{p\in\mathcal L}
\alpha_p\lambda_{\sigma_\mu}(p)\big|
=
\sum_{p\in\mathcal L}
|\lambda_{\sigma_\mu}(p)|
\gg \sum_{p\in\mathcal L}|\lambda_{\sigma_\mu}(p)|^2
\asymp_{\pi}
\frac{L}{\log L}.
\end{align*}

Theorem \ref{prop11.2} and Proposition \ref{prop12.5} therefore yield
\begin{align*}
\frac{|L(1/2,\sigma\times\pi')|^2}
{\langle\phi',\phi'\rangle
T^{\frac{n(n+1)}{2}}}
\prod_{j=1}^m\frac{1}{L(1,\pi_j,\Ad)}
\ll_{\varepsilon,\pi,\pi'}
(TL)^\varepsilon
\left(
L^{-1}+T^{-1/2}L^{3n-2}
\right).
\end{align*}
Here we have used $L(1,\pi_j\otimes|\det|^{it},\Ad)=L(1,\pi_j,\Ad)$. Taking $L=T^{1/(2(3n-1))}$ and using Lemma \ref{lem4.6} proves \eqref{nb.}.
\end{proof}

When $n=1$, the preceding argument remains valid with the empty products and
the $A_H$-integrals omitted, as stipulated in \textsection\ref{sec10}. The
general amplifier gives the exponent $5/12$ in \eqref{1}, while the tempered
choice $\ell_p=1$ gives $3/8$ in \eqref{1.4*}. Taking
$\pi=\textbf{1}\boxplus\textbf{1}$ then gives the exponent $3/16$ in
Corollary \ref{cor1.2}. This completes the proof of Theorem \ref{A} for
$n=1$.

\section{\texorpdfstring{Quantitative Nonvanishing for
$\mathrm{GL}(n+1)\times\mathrm{GL}(n)$}
{Quantitative Nonvanishing for GL(n+1) x GL(n)}}\label{sec12}
We combine Theorems \ref{C} and \ref{E} to prove Theorem \ref{thmC}.
\subsection{Choice of Local and Global Data}\label{sec12.1}
 
\subsubsection{Automorphic weights $\phi_j'\in\pi_j'$}
Let $M\geq1$, and fix a prime $p_*\nmid M$. For $j=1,2$, let
$\pi_j'=\otimes_{p\leq\infty}\pi_{j,p}'$ be a unitary cuspidal
representation of $G'(\mathbb A)$ with central character $\omega_j$ and
conductor $M_j'$. Assume that
\begin{align*}
M_1'M_2'>1,\quad (M_1',M_2')=(p_*M,M_1'M_2')=1,\quad
\pi_{1,\infty}'\simeq\pi_{2,\infty}'.
\end{align*}
For $p\mid M_1'M_2'$, the essential-vector identity \cite{Kim10} gives
\begin{equation}\label{12.3}
\int_{N'(\mathbb Q_p)\backslash\overline{G'}(\mathbb Q_p)}
W_{1,p}'(x)\overline{W_{2,p}'(x)}|\det x|_p\,dx
=L_p(1,\pi_{1,p}'\times\widetilde\pi_{2,p}'),
\end{equation}
where $W_{j,p}'$ is the normalized Whittaker new vector. Let $W_\infty'$
be the common Archimedean vector of Definition \ref{def4.2}. Let
$\phi_j'\in\pi_j'$ be the cusp form whose Whittaker function has local
components $W_\infty'$ and $W_{j,p}'$. The conductor assumptions imply
$\pi_1'\not\simeq\pi_2'$.

\subsubsection{The cuspidal family}\label{12.1.3}
Fix a unitary Hecke character $\omega$ of
$\mathbb Q^\times\backslash\mathbb A^\times$, trivial at infinity, and a
supercuspidal representation $\pi_{p_*}$ of $G(\mathbb Q_{p_*})$ with
central character $\omega_{p_*}$. Also fix
a unitary irreducible admissible representation $\pi_\infty$ of
$\mathrm{PGL}_{n+1}(\mathbb R)$ with uniform parameter growth of size
$(T;c_\infty,C_\infty)$.
Let $\mathcal U_\infty(\pi_\infty)$ denote the fixed-width spectral
neighborhood of $\pi_\infty$ selected by the Archimedean microlocal test
function; its rapidly decaying tails will be absorbed throughout.

\begin{defn}\label{defn11.1}
With the fixed character $\omega$ suppressed from the notation, let
$\mathcal{A}_0(T,M;\pi_{\infty},\pi_{p_*},\pi_1',\pi_2')$ be the set of
equivalence classes of cuspidal representations
$\sigma=\otimes_{p\leq\infty}\sigma_p\in\mathcal{A}_0([G];\omega)$ such that
\begin{itemize}
 \item $\sigma_\infty\in\mathcal U_\infty(\pi_\infty)$ and
 $\sigma_{p_*}\simeq\pi_{p_*}$;
\item $\sigma_p$ has a nonzero $K_p(MM_1')$-fixed vector at $p\nmid p_*M_2'$;
\item $\sigma_p$ has a nonzero $I_p(M_2')$-fixed vector at $p\mid M_2'$.
\end{itemize}
\end{defn}
 
The microlocal construction localizes the Archimedean spectrum to
$\mathcal U_\infty(\pi_\infty)$. The uniform parameter growth condition
identifies the Plancherel mass of this window, and the Weyl law gives
\begin{equation}\label{12.2.2}
|\mathcal{A}_0(T,M;\pi_{\infty},\pi_{p_*},\pi_1',\pi_2')|\asymp  T^{\frac{n(n+1)}{2}+o(1)}M^n,
\end{equation}
with constants depending on $M_1'$, $M_2'$, and $p_*$. 
 
Let $W_{p_*}$ be a normalized Whittaker new vector of $\pi_{p_*}$, and,
for $j=1,2$, let $W_{j,p_*}'$ be the normalized spherical Whittaker vector
of $\pi_{j,p_*}'$. The local Rankin--Selberg identity \cite{JPSS81} gives
\begin{equation}\label{12.1}
\mathcal P_{j,p_*}:=\int_{N'(\mathbb Q_{p_*})\backslash G'(\mathbb Q_{p_*})}
W_{p_*}(\iota(x))W_{j,p_*}'(x)\,dx=L_{p_*}(1/2,\pi_{p_*}\times\pi_{j,p_*}')\neq 0.
\end{equation}

\subsubsection{The test function}
Define $f=\otimes_{p\leq\infty}f_p$ as follows.
\begin{itemize}
\item $f_{\infty}$ is defined as in \textsection\ref{11.1.1};
\item  $p\mid M,$ $f_p$ is defined in \textsection\ref{sec3.3.}; 
\item $p=p_*,$ $f_p(g):=\langle \pi_{p_*}(g)W_{p_*},W_{p_*}\rangle,$ $g\in G(\mathbb{Q}_{p_*});$
\item  $p\nmid p_*MM_1'M_2'$ and $p<\infty,$ $f_p$ is defined in \textsection\ref{11.1.6};
\item $p\mid M_1'$. Define the function 
\begin{align*}
\widehat{\Phi}_p(b_1,\cdots, b_n)=\begin{cases}
\overline{\omega}_{1,p}(b_n)\omega_{2,p}(b_n),\ &\text{if $b_1,\cdots, b_{n-1}\in M_1'\mathbb{Z}_p,$ $b_n\in\mathbb{Z}_p^{\times},$}\\
0,\ & \text{otherwise}.
\end{cases}
\end{align*}
Let $\psi_p$ be the standard unramified additive character of $\mathbb{Q}_p$ (see \textsection\ref{notation}). Let $\Phi_p$ be the Fourier inversion of $\widehat{\Phi}_p$ relative to $\psi_p$. Define $f_p(g_p):=\int_{Z(\mathbb{Q}_p)}\tilde{f}_p(z_pg_p)\omega_p(z_p)d^{\times}z_p,$ $g_p\in G(\mathbb{Q}_p),$ where 
\begin{align*}
\tilde{f}_p(g_p)=\frac{\textbf{1}_{M_{n,n}(\mathbb{Z}_p)}(A)\Phi_p(\mathfrak{b})\textbf{1}_{M_{1,n}(\mathbb{Z}_p)}(\mathfrak{c})\textbf{1}_{\mathbb{Z}_p}(d)\textbf{1}_{G(\mathbb{Z}_p)}(g_p)}{\Vol(K_p'(M_1'))}
\end{align*}
for all $g_p=\begin{pmatrix}
A&\mathfrak{b}\\
\mathfrak{c}&d
\end{pmatrix}\in G(\mathbb{Q}_p)$. Here $K_p'(M_1')$ is the Hecke congruence subgroup of $K_p'$ of level $e_p(M_1')$. Recall \eqref{3.4.} for its definition. 
 		
In particular, $f_p$ is right $K_p(M_1')$-invariant, where $K_p(M')$ is defined in the manner of \eqref{equ3.4}. 
 
\item If $p\mid M_2'$, define $\widehat\Phi_p$ and $\Phi_p$ as above with
$M_2'$ in place of $M_1'$. Set
$f_p(g_p):=\int_{Z(\mathbb{Q}_p)}\tilde{f}_p(z_pg_p)\omega_p(z_p)d^{\times}z_p,$ where
\begin{align*}
\tilde{f}_p(g_p)=\frac{\textbf{1}_{K_p'(M_2')}(A)\Phi_p(\mathfrak{b})\textbf{1}_{M_{1,n}(\mathbb{Z}_p)}(\mathfrak{c})\textbf{1}_{\mathbb{Z}_p}(d)\textbf{1}_{G(\mathbb{Z}_p)}(g_p)}{\Vol(K_p'(M_2')}
\end{align*}
for all $g_p=\begin{pmatrix}
A&\mathfrak{b}\\
\mathfrak{c}&d
\end{pmatrix}\in G(\mathbb{Q}_p)$. Then $f_p$ is right $I_p(M_2')$-invariant, where
\begin{align*}
I_p(M_2'):=\big\{(g_{i,j})\in G(\mathbb Z_p):
g_{i,j}\in M_2'\mathbb Z_p\ \text{for }1\leq j<i\leq n+1\big\}.
\end{align*}
 \end{itemize}

\subsection{The Relative Trace Formula}
Applying Theorem \ref{C} to the data $(f,\phi_1',\phi_2')$ at $\mathbf s=\mathbf0$ gives
$J_{\Spec}^{\Reg,\heartsuit}(f,\mathbf0)=J_{\Geo}^{\Reg,\heartsuit}(f,\mathbf0)$.

\subsubsection{The spectral side}
\begin{lemma}\label{lemma12.1}
	Let notation be as before. Then 
	\begin{align*}
	\big|J_{\Spec}^{\Reg,\heartsuit}(f,\mathbf{0})\big|\ll \sum_{\sigma\in \mathcal{A}_0(T,M;\pi_{\infty},\pi_{p_*},\pi_1',\pi_2')}\frac{T^{-\frac{n^2}{2}+\varepsilon}\big|L(1/2,\sigma\times\pi_1')L(1/2,\sigma\times\pi_2')\big|}{L(1,\sigma,\Ad)},
\end{align*}
where the implied constant depends on $\pi_{\infty}$, $\pi_{p_*}$,
$\phi_1'$, and $\phi_2'$.
\end{lemma}
\begin{proof}
The matrix coefficient $f_{p_*}$ annihilates the noncuspidal spectrum.
The remaining local data restrict the cuspidal contribution to
$\mathcal A_0(T,M;\pi_\infty,\pi_{p_*},\pi_1',\pi_2')$, up to the rapidly
decaying Archimedean tails fixed above. Unfolding the two
periods gives a product, rather than an absolute square, because the
automorphic weights are $\phi_1'$ and $\phi_2'$. The local estimates used in
the proof of Proposition \ref{thm6}, together with Lemma \ref{prop4.11} and
\eqref{12.1} for $j=1,2$, give the stated bound.
\end{proof}

\subsubsection{The geometric side}
\begin{lemma}\label{lem12.2}
	Let notation be as before. Then $J^{\Reg}_{\Geo,\du}(f,\mathbf{s})\equiv 0$.
\end{lemma}
\begin{proof}
Take $\mathbf s=(s,0)$ with $\Re(s)>1$, where the dual term converges
absolutely. At any $p\mid M_1'M_2'$, exactly one of $W_{1,p}'$ and
$W_{2,p}'$ is spherical. The support of $f_p$ reduces the corresponding local
factor to the pairing of spherical data with a ramified new vector, which
vanishes by $K_p'$-orthogonality. Since $M_1'M_2'>1$, at least one such prime
exists; hence the dual term vanishes in this half-plane, and therefore
identically by meromorphic continuation.
\end{proof}

\begin{prop}\label{proposition12.2}
Let notation be as before. Then
\begin{equation}\label{12.44}
	\big|J_{\Geo}^{\Reg,\heartsuit}(f,\mathbf{0})\big|\gg_{\varepsilon,\pi_1',\pi_2',\pi_{p_*}} T^{\frac{n}{2}-\varepsilon}M^n\|W_{\infty}'\|_{L^2}^2,
\end{equation}
where $\|\cdot\|_{L^2}$ denotes the $L^2$-norm.
\end{prop}
\begin{proof}
The supercuspidal matrix coefficient at $p_*$ annihilates
$J^{\Reg,\RNum{1}}_{\Geo,\bi}$, while Lemma \ref{lem12.2} annihilates the
dual term. Moreover, Cauchy--Schwarz and the proof of Theorem
\ref{thm9.1} give
\begin{align*}
\big|J^{\Reg,\RNum{2}}_{\Geo,\bi}(f,\mathbf{0})\big|\ll_{\varepsilon,\pi_1',\pi_2',\pi_{p_*}}M^{n-1}T^{\frac{n}{2}-\frac{1}{2}+\varepsilon}\|W_{\infty}'\|_2^2,
\end{align*}
so it remains to evaluate the small-cell term. Set
$S:=\{\infty\}\bigsqcup\,\{p:p\mid p_*M_1'M_2'\}$. For $\mathbf s=(s,0)$ and
$\Re(s)>0$, unfolding gives
\begin{align*}
J^{\Reg}_{\Geo,\sm}(f,\mathbf{s})
=\frac{L^{(S)}(1+s,\pi_1'\times\widetilde\pi_2')}
{\Vol(K_0(M))}\prod_{v\in S}\mathcal I_v(f,s).
\end{align*}
All local integrals converge at $s=0$, and the Rankin--Selberg factor is
holomorphic there because $\pi_1'\not\simeq\pi_2'$. At $p\mid M_1'M_2'$,
the definition of $f_p$ and \cite[Theorem 2.1.1]{Kim10} give $\mathcal I_p(f,0)=L_p(1,\pi_{1,p}'\times\widetilde\pi_{2,p}')$.

At $p=p_*$, Fourier inversion and local Plancherel give
\begin{align*}
\mathcal I_{p_*}(f,0)
=c_{n,p_*}\mathcal P_{1,p_*}\overline{\mathcal P_{2,p_*}},
\qquad
\big|\mathcal I_{p_*}(f,0)\big|
\gg_{\pi_{p_*},\pi_{1,p_*}',\pi_{2,p_*}'}1,
\end{align*}
by \eqref{12.1}. At the Archimedean place, the same unfolding gives
\begin{align*}
\mathcal I_\infty(f,0)
=c_{n,\infty}\int_{G'(\mathbb R)}f_\infty(\iota(y))
\overline{\langle\pi_\infty'(y)W_\infty',W_\infty'\rangle}\,dy
\gg T^{\frac n2-\varepsilon}\|W_\infty'\|_2^2.
\end{align*}
Since $\Vol(K_0(M))^{-1}\asymp M^n$ and the remaining local factors are
fixed and nonzero, the small-cell term satisfies
\begin{align*}
\big|J_{\Geo,\sm}^{\Reg}(f,\mathbf0)\big|
\gg T^{\frac n2-\varepsilon}M^n\|W_\infty'\|_2^2.
\end{align*}
The triangle inequality and the preceding bound for the regular contribution
now prove \eqref{12.44}.
\end{proof}
\subsection{Nonvanishing Problems}\label{12.4}
Write
\begin{align*}
\mathscr N(T,M):=\#\big\{\pi\in
\mathcal A_0(T,M;\pi_\infty,\pi_{p_*},\pi_1',\pi_2'):
L(1/2,\pi\times\pi_1')L(1/2,\pi\times\pi_2')\neq0\big\}.
\end{align*}
Theorem \ref{thmC} follows from the next result.
\begin{thmx}\label{G}
Retain the notation of \textsection\ref{sec12.1}. Then
\begin{equation}\label{12.5.}
\mathscr N(T,M)\gg_\varepsilon (TM)^{-\varepsilon}\begin{cases}
M^{\frac{n}{3n^2+n-1}}T^{\frac{1}{2(3n^2+n-1)}}, & \text{if $M\leq T^{\frac{n+1}{n^2-n-1}}$,}\\
M^{\frac{1}{n^2+n+1}}T^{\frac{1}{2(n^2+n+1)}}, & \text{if $M>T^{\frac{n+1}{n^2-n-1}}$.}
\end{cases}
\end{equation}
In particular,
\begin{equation}\label{12.5}
\mathscr N(T,M)\gg_\varepsilon
|\mathcal A_0(T,M;\pi_\infty,\pi_{p_*},\pi_1',\pi_2')|^
{\frac{1}{n(n+1)(3n^2+n-1)}-\varepsilon}.
\end{equation}
The implied constants depend on $\varepsilon$, $c_\infty$, $C_\infty$,
$\pi_1'$, $\pi_2'$, and $\pi_{p_*}$.
\end{thmx}
\begin{proof}
Theorem \ref{C}, Lemma \ref{lemma12.1}, and Proposition
\ref{proposition12.2} give
\begin{align*}
\sum_{\pi\in\mathcal A_0}
\frac{|L(1/2,\pi\times\pi_1')L(1/2,\pi\times\pi_2')|}
{L(1,\pi,\Ad)\|W_\infty'\|_2^2}
\gg T^{\frac{n(n+1)}2-\varepsilon}M^n.
\end{align*}
Here $\mathcal A_0=\mathcal A_0(T,M;\pi_\infty,\pi_{p_*},\pi_1',\pi_2')$.
Applying Theorem \ref{E} to both central values and using
$\langle\phi_j',\phi_j'\rangle\asymp\|W_\infty'\|_2^2$ yields
\begin{align*}
T^{\frac{n(n+1)}2+\varepsilon}M^{n+\varepsilon}\mathbf L^{-2}
\mathscr N(T,M)
\gg T^{\frac{n(n+1)}2-\varepsilon}M^n.
\end{align*}
Thus $\mathscr N(T,M)\gg_\varepsilon(TM)^{-\varepsilon}\mathbf L^2$,
which is \eqref{12.5.}; combining this with \eqref{12.2.2} gives
\eqref{12.5}.
\end{proof}

\bibliographystyle{alpha}

\bibliography{SC}

\end{document}